# Efficient Alternative Finite Difference WENO Schemes for Hyperbolic Conservation Laws

By

Dinshaw S. Balsara[1,2], Deepak Bhoriya[1], Chi-Wang Shu[3] and Harish Kumar[4]

[1]Physics Department, [2]ACMS Department, University of Notre Dame

[3]Division of Applied Mathematics, Brown University

[4]Department of Mathematics, Indian Institute of Technology, Delhi

**Abstract**

Higher order finite difference Weighted Essentially Non-Oscillatory (WENO) schemes for conservation laws are extremely popular because, for multidimensional problems, they offer high order accuracy at a fraction of the cost of finite volume WENO or DG schemes. Such schemes come in two formulations. The very popular classical finite difference WENO (FD-WENO) method (Shu and Osher, J. Comput. Phys., 83 (1989) 32-78) relies two reconstruction steps applied to two split fluxes. However, the method cannot accommodate different types of Riemann solvers and cannot preserve free stream boundary conditions on curvilinear meshes. This limits its utility. The alternative finite difference WENO (AFD-WENO) method can overcome these deficiencies, however, much less work has been done on this method. The reasons are three-fold. First, it is difficult for the casual reader to understand the intricate logic that requires higher order derivatives of the fluxes to be evaluated at zone boundaries. The analytical methods for deriving the update equation for AFD-WENO schemes are somewhat recondite. To overcome that difficulty, we provide an easily accessible script that is based on a computer algebra system in Appendix A of this paper. Second, the method relies on interpolation rather than reconstruction, and WENO interpolation formulae have not been documented in the literature as thoroughly as WENO reconstruction formulae. In this paper, we explicitly provide all necessary WENO interpolation formulae that are needed for implementing AFD-WENO up to ninth order. The third reason is that AFD-WENO requires higher order derivatives of the fluxes to be available at zone boundaries.

Since those derivatives are usually obtained by finite differencing the zone-centered fluxes, they become susceptible to a Gibbs phenomenon when the solution is non-smooth. The inclusion of those fluxes is also crucially important for preserving the order property when the solution is smooth. This has limited the utility of AFD-WENO in the past even though the method *per se* has many desirable features. Some efforts to mitigate the effect of finite differencing of the fluxes have been tried, but so far they have been done on a case by case basis for the PDE being considered. In this paper we find a general-purpose strategy that is based on a different type of WENO interpolation. This new WENO interpolation takes the first derivatives of the fluxes at zone centers as its inputs and returns the requisite non-linearly hybridized higher order derivatives of flux-like terms at the zone boundaries as its output.

With these three advances, we find that AFD-WENO becomes a robust and general-purpose solution strategy for large classes of conservation laws. It allows any Riemann solver to be used. AFD-WENO has a computational complexity that is entirely comparable to classical FD-WENO, because it relies on two interpolation steps which cost the same as the two reconstruction steps in classical FD-WENO. We apply the method to several stringent test problems drawn from Euler flow, relativistic hydrodynamics and ten-moment equations. The method meets its design accuracy for smooth flow and can handle stringent problems in one and multiple dimensions.

**I) Introduction**

Essentially Non-Oscillatory (ENO) methods for the high order accurate numerical simulation of conservation laws in a finite volume setting were invented in the pioneering work of Harten *et al.* [28]. Shu and Osher [53], [54] soon followed up on this work by introducing finite difference ENO methods. The finite difference versions of higher order schemes are much more computationally efficient compared to their finite volume counterparts. Early ENO schemes suffered from the deficiency that certain problems could cause very rapid switching of the stencil, resulting in a loss of accuracy. Weighted Essentially Non-Oscillatory (WENO) schemes were invented to overcome this deficiency (Liu, Osher & Chan [43], Jiang & Shu [31]). The methods were extended to seventh, ninth and eleventh orders by Balsara & Shu [4] and much later to seventeenth order by Gerolymos, Sénéchal & Vallet [26]. Some of the early deficiencies of WENO schemes stemmed from a loss of accuracy at critical points, and a way out of this problem was presented in Henrick, Aslam & Powers [30], Borges et al. [18] and Castro et al. [20]. A finite difference WENO formulation that applies to systems with non-conservative products has recently been formulated by Balsara *et al.* [12]. For a comprehensive review of WENO schemes, see Shu [55] [56].

The original papers by Shu and Osher [53], [54] contained not one but rather two styles of thinking about finite difference ENO/WENO schemes. One of those styles of thinking, in Shu and Osher [54], has become extremely popular. It consists of splitting the locally Lax-Friedrichs (LLF) flux into two parts – a left-going and a right-going part. The two parts were suitably upwinded via an appropriate choice of stencils. Using a trick stemming from the fundamental theorem of integral calculus (later referred to as the Shu-Osher Lemma), it was shown that a one-dimensional finite volume reconstruction from the point values of the upwinded fluxes would yield a high order finite difference scheme. This then allowed the usage of the same one-dimensional finite volume reconstruction subroutine to approximate multi-dimensional conservation laws dimension by dimension in the high order conservative finite difference schemes. We refer to this algorithm as the classical FD-WENO algorithm. This algorithm has become highly popular to the point where most papers that cite Shu and Osher [53], [54] do so because of this algorithm. However, there was another algorithm that was developed first in Shu and Osher [53]. For a very long time, almost nobody paid attention to that algorithm. A paper by Merriman [46] made some headway in understanding that algorithm. We refer to that algorithm as the alternative formulation of finite difference WENO (AFD-WENO) in this paper, following the terminology first used in Jiang, Shu and Zhang [32]. Subsequent interest in AFD-WENO has emerged sporadically (Jiang, Shu and Zhang [32], [33], Zheng, Shu, Qiu [61], Gao *et al*. [25]) but it is becoming clearer that AFD-WENO represents a strong alternative algorithm to the classical FD-WENO algorithm. Furthermore, AFD-WENO has many strong points relative to the classical FD-WENO. In the next paragraph we explain the reasons why this is so.

One of the reasons for AFD-WENO's slower acceptance is that the it is difficult for the casual reader to understand the intricate logic that requires higher order derivatives of the fluxes to be evaluated at zone boundaries. The analytical methods that give rise to the AFD-WENO update equation are difficult to understand. Because our *first goal* is to make the method more accessible, we provide a script based on a computer algebra system in Appendix A of this paper which shows that the AFD-WENO update equation can be easily derived with modern computational tools.

Classical FD-WENO relies on the availability of a smooth flux. But this restricts it to an LLF Riemann solver or a variant of a Roe-type Riemann solver. This restriction is fundamentally a consequence of the flux reconstruction. In recent years we have seen many different Riemann solvers emerge which have special attributes that make them very useful in various application areas. Those Riemann solvers do not fit well into the strictures of classical FD-WENO. The AFD-WENO algorithm is free of such strictures – any type of Riemann solver can be invoked in a pointwise fashion at the zone boundaries. This makes a well-designed AFD-WENO very broadly applicable to many application areas. Classical FD-WENO also does not take well to preserving the free stream condition on curvilinear meshes; whereas AFD-WENO can indeed take well to curvilinear meshes (Jiang, Shu and Zhang [32], [33]). In the discussions that are contained in this paper we document other potential advantages of AFD-WENO.

For all its advantages, AFD-WENO has also proven to be a little harder to work with, perhaps because it is not as well-developed as classical FD-WENO. To begin with, it relies on WENO interpolation rather than WENO reconstruction. Since the latter is much better known than the former, the widespread acceptance of AFD-WENO has suffered. There is a further reason for its hitherto fore lack of acceptance. It stems from the fact that one has to evaluate higher order derivatives of the flux at the zone boundaries if one wants to use AFD-WENO. These higher derivatives of the flux can cause spurious oscillations when the solution is non-smooth on the mesh. There has been very little work to control these oscillations and the few efforts that have been made are very specific to the PDE being considered (Zheng, Shu, Qiu [61], Gao *et al*. [25]). As a result, it has not been possible to develop AFD-WENO schemes as a general-purpose tool for numerically solving large classes of hyperbolic conservation laws. The *second goal* of this paper is to overcome this limitation by showing that our well-designed AFD-WENO scheme is indeed a general-purpose solver for conservation laws. Consequently, as part of our numerical results, we will show that the same AFD-WENO algorithm can be applied very generally to large classes of hyperbolic conservation laws.

All the alternative finite difference WENO (AFD-WENO) schemes are based on interpolation. It is therefore worthwhile to make a distinction between reconstruction and interpolation. Reconstruction is used in all finite volume schemes and also in many popular finite difference WENO (FD-WENO) reconstruction schemes where the fluxes are reconstructed. It consists of starting with the zone averages in a given stencil and obtaining therefrom the high degree polynomial whose integration over each of the zones of the stencil matches the original zone averages. Interpolation is used less often in the numerical solution of conservation laws, however it is the approach that will be used in this paper. It consists of starting with the point values at each of the zone centers of a stencil and obtaining therefrom the high degree polynomial that matches those point values. Therefore, the two words, reconstruction and interpolation, carry different connotations. When applied to the same stencil, reconstruction and interpolation produce polynomials with the same degree. However, the underlying polynomial coefficients that are produced by invoking reconstruction or interpolation on a given stencil can indeed be different. Standard WENO concepts like linear weights, smoothness indicators, normalized non-linear weights etc. are often the same for reconstruction and interpolation; so there is indeed a beneficial transference of knowledge between them. WENO reconstruction, especially as it relates to FD-WENO, has been amply documented in the literature starting from Jiang and Shu [31], Balsara and Shu [4] and continuing through Balsara, Garain and Shu [8], where it was presented in its most polished form. It turns out that the WENO interpolation, as it relates to AFD-WENO, has not been thoroughly documented in the literature. This may be one reason why the method has not been widely embraced by practitioners. The *third goal* of this paper is to thoroughly document the different WENO interpolation methods that are useful in AFD-WENO. To facilitate easy adoption of AFD-WENO by the greater community, some of the sections of this paper have been designed so that they constitute a one-stop-shop for WENO interpolation formulae.

The paper is divided into the following sections. In Section II we quickly present the AFD-WENO formulation and identify the advances in WENO interpolation that are needed to turn AFD-WENO into a general-purpose algorithm. In Section III we document some of the WENO interpolation algorithms as they apply to zone centered variables. In Section IV we document new types of WENO interpolation algorithms that are indeed invoked at zone boundaries. Section V gives pointwise implementation-related details. Section VI presents an accuracy analysis for a variety of hyperbolic conservation laws. Section VII presents results involving various hyperbolic conservation laws in one dimension. Section VIII does the same in multiple dimensions. Section IX draws some conclusions.

## II) An AFD-WENO Algorithm – Description of Philosophy and Formulation

This Section is split into two Sub-sections. The first Sub-section clearly describes the philosophy behind the AFD-WENO algorithm and should be useful to any reader. The second Sub-section documents a class of AFD-WENO schemes up to ninth order where the Gibbs phenomenon is suppressed via a redesigned WENO algorithm.

### II.a) AFD-WENO Scheme Design Philosophy

We now describe the philosophy behind the Alternative Finite Difference WENO (AFD-WENO) scheme. We focus on the solution of a one-dimensional PDE system given by

$$\partial_t \mathbf{U} + \partial_x \mathbf{F} = 0 \tag{1}$$

Let us establish some notation. Please see Fig. 1. It shows a small sub-section of the mesh function in a few adjacent zones. The zones are labeled by "$i-1, i, i+1$" etc. and their zone centers are denoted by "$x_{i-1}, x_i, x_{i+1}$" etc. The zone boundaries of each zone "$i$" are denoted by "$x_{i-1/2}$" and "$x_{i+1/2}$" with $\Delta x = x_{i+1/2} - x_{i-1/2}$ being a constant because we have assumed a uniform mesh. The associated mesh functions are specified in pointwise fashion at the zone centers and are labeled by "$\mathbf{U}_{i-1}, \mathbf{U}_i, \mathbf{U}_{i+1}$" etc. Here "$\mathbf{U}$" is a vector of primal variables for the hyperbolic PDE that we are considering. For the moment we consider a one-dimensional mesh, but because this is a finite difference scheme, the method can be extended dimension-by-dimension to multiple dimensions. Using any WENO interpolation strategy, as applied to the point values of the mesh function, we can obtain a suitably high order interpolation within each zone. (In the next two sections we will provide a lot of detail on WENO interpolation.) The interpolation within zone "$i$" gives us interpolated $\hat{\mathbf{U}}^-_{i+1/2}$ and $\hat{\mathbf{U}}^+_{i-1/2}$ at the right and left boundaries of the zone being considered, see Fig. 1. We will use a caret to denote such interpolated variables. Note that $\hat{\mathbf{U}}^-_{i+1/2}$ is available at the left side of the zone boundary $x_{i+1/2}$ and $\hat{\mathbf{U}}^+_{i-1/2}$ is available at the right side of the zone boundary $x_{i-1/2}$. Similarly, using our WENO interpolation within the zone "$i+1$", we obtain $\hat{\mathbf{U}}^-_{i+3/2}$ and $\hat{\mathbf{U}}^+_{i+1/2}$ at the right and left boundaries of that zone. Note that $\hat{\mathbf{U}}^+_{i+1/2}$ is available at the right side of the zone boundary $x_{i+1/2}$. Likewise, using our WENO interpolation within the zone "$i-1$", we obtain $\hat{\mathbf{U}}^-_{i-1/2}$ and $\hat{\mathbf{U}}^+_{i-3/2}$ at the right and left boundaries of that zone. Consequently, $\hat{\mathbf{U}}^-_{i-1/2}$ is available at the left side of the zone boundary $x_{i-1/2}$. We assume that a Riemann solver with left and right states given by $\hat{\mathbf{U}}^-_{i+1/2}$ and $\hat{\mathbf{U}}^+_{i+1/2}$ is applied at zone boundary $x_{i+1/2}$ and it yields a resolved state $\mathbf{U}^*_{i+1/2}$ that overlies the zone boundary, as seen in Fig. 1. The resolved state, as well as the structure of the Riemann fan, can be used to obtain a resolved flux given by $\mathbf{F}^*\left(\hat{\mathbf{U}}^-_{i+1/2}, \hat{\mathbf{U}}^+_{i+1/2}\right)$. Likewise, we assume that a Riemann solver with left and right states given by $\hat{\mathbf{U}}^-_{i-1/2}$ and $\hat{\mathbf{U}}^+_{i-1/2}$ is applied at zone boundary $x_{i-1/2}$ and it produces a resolved flux within the Riemann fan given by $\mathbf{F}^*\left(\hat{\mathbf{U}}^-_{i-1/2}, \hat{\mathbf{U}}^+_{i-1/2}\right)$. Since we are hoping to produce a very light-weight scheme, we assume that the Riemann solver will be something very simple like HLL, or HLLI, or HLLC, or LLF with a single resolved state within the Riemann fan. One of the advantages of the AFD-WENO algorithm is that one it is agnostic to the type of Riemann solver that is used. This has the result that one may even choose a Riemann solver that maintains stationary linearly degenerate fields on the mesh, if the application would benefit from it. Similarly, if the application would benefit from a positivity preserving Riemann solver, then such a Riemann solver can also be used. The upshot is that by being agnostic to the kind of Riemann solver that is used, the AFD-WENO algorithm gives the user considerably greater flexibility compared to the classical FD-WENO algorithm (which has in the past only been formulated in its LLF and RF variants). This completes our description of Fig. 1. The Riemann solver that we will use in this entire work will be the HLLI Riemann solver (and its LLFI variant) from Dumbser and Balsara [24] where the philosophy of such Riemann solvers is explained. In Section V of Balsara *et al*. [12] we also abstract results from Dumbser and Balsara [24] in a notation that is more suited for the WENO-AO interpolation shown in Fig. 1. We illustrate Fig. 1 within the context of WENO with Adaptive Order (WENO-AO) interpolation (Balsara, Garain and Shu [8]); however, any form of WENO interpolation would suffice for the purposes of the discussion in this Section.

Let us now assume that we have a high order pointwise WENO interpolation strategy which gives us high order interpolation polynomials within each zone. Say too that at each zone boundary $x_{i+1/2}$ we can use these interpolation polynomials to obtain $\hat{\mathbf{U}}^-_{i+1/2}$ and $\hat{\mathbf{U}}^+_{i+1/2}$. Say also that we invoke the Riemann solver at each zone boundary to obtain $\mathbf{F}^*\left(\hat{\mathbf{U}}^-_{i+1/2}, \hat{\mathbf{U}}^+_{i+1/2}\right)$. Say we naively assert a discrete in space but continuous in time update in the zone "*i*" of the form

$$\partial_t \mathbf{U}_i = -\frac{1}{\Delta x}\left(\mathbf{F}^*\left(\hat{\mathbf{U}}^-_{i+1/2}, \hat{\mathbf{U}}^+_{i+1/2}\right) - \mathbf{F}^*\left(\hat{\mathbf{U}}^-_{i-1/2}, \hat{\mathbf{U}}^+_{i-1/2}\right)\right) \qquad (2)$$

Even if the interpolation has been a very high order pointwise WENO interpolation, eqn. (2) will only yield a FD-WENO scheme that is at best second order accurate! The reason is that the finite difference on the right hand side of eqn. (2) only yields a second order accurate representation of the flux gradient $\partial_x \mathbf{F}$ that is evaluated at the zone center. In other words, the two individual fluxes in eqn. (2), when evaluated at the zone boundaries, will indeed be pointwise high order accurate up to the accuracy of the interpolation. However, their finite differencing in the right hand side of eqn. (2) will only yield a second order accurate finite difference approximation at the zone center! For readers who are more attuned to finite volume schemes, this might seem like a counter-intuitive claim; but it is nevertheless valid.

Let us illustrate the somewhat counter-intuitive claim made above with a simple example. For our simple example, we focus only on the task of achieving third order accuracy. We assume for the purposes of this paragraph that we are dealing with scalar fluxes that are as smooth and differentiable as we desire. At the zone center, taken at $x = 0$, we can make a Taylor series expansion for the flux as

$$\mathbf{F}(x) = f_0 + x(\partial_x f)_0 + \frac{x^2}{2}(\partial_x^2 f)_0 + \frac{x^3}{6}(\partial_x^3 f)_0 + \ldots \tag{3}$$

In the above equation all the terms of the Taylor series, $f_0$, $(\partial_x f)_0$, $(\partial_x^2 f)_0$, $(\partial_x^3 f)_0$ are all evaluated at $x = 0$. We can evaluate eqn. (1) and its higher derivatives at $x = \pm \Delta x / 2$ to get

$$\begin{aligned}
\mathbf{F}(x)\big|_{x=\pm\Delta x/2} &= f_0 \pm \frac{\Delta x}{2}(\partial_x f)_0 + \frac{\Delta x^2}{8}(\partial_x^2 f)_0 \pm \frac{\Delta x^3}{48}(\partial_x^3 f)_0 + \ldots \quad ; \\
\partial_x \mathbf{F}(x)\big|_{x=\pm\Delta x/2} &= (\partial_x f)_0 \pm \frac{\Delta x}{2}(\partial_x^2 f)_0 + \frac{\Delta x^2}{8}(\partial_x^3 f)_0 \pm \ldots \quad ; \\
\partial_x^2 \mathbf{F}(x)\big|_{x=\pm\Delta x/2} &= (\partial_x^2 f)_0 \pm \frac{\Delta x}{2}(\partial_x^3 f)_0 + \ldots \quad .
\end{aligned} \tag{4}$$

By finite differencing the point values of the fluxes at $x = \pm \Delta x / 2$ we can now easily see that

$$\frac{1}{\Delta x}\left[\mathbf{F}(x)\big|_{x=\Delta x/2} - \mathbf{F}(x)\big|_{x=-\Delta x/2}\right] = (\partial_x f)_0 + \frac{\Delta x^2}{24}(\partial_x^3 f)_0 \tag{5}$$

The above expression makes it easy to see that the presence of $\Delta x^2 (\partial_x^3 f)_0 / 24$ prevents eqn. (5) from being a third order accurate expression. The way out of this conundrum is also easy to see. Say we take our numerical fluxes at $x = \pm \Delta x / 2$ to be

$$\mathbf{F}^{num}_{x=\pm\Delta x/2} = \left[\mathbf{F}(x)\big|_{x=\pm\Delta x/2}\right] - \frac{\Delta x^2}{24}\left[\partial_x^2 \mathbf{F}(x)\big|_{x=\pm\Delta x/2}\right] \tag{6}$$

By finite differencing the numerical fluxes from the above equation we can now easily see that

$$\frac{1}{\Delta x}\left[\mathbf{F}^{num}_{x=\Delta x/2} - \mathbf{F}^{num}_{x=-\Delta x/2}\right] = (\partial_x f)_0 + O(\Delta x^4) \tag{7}$$

In other words, with the correction term in eqn. (6), the numerical flux becomes fourth order accurate *in a pointwise finite difference sense*! The concept in the present paragraph was presented without much detail in Shu and Osher [53], causing Merriman [46] to write a paper explaining such nuances of the Alternative FD-WENO (AFD-WENO) scheme. In Appendix A we provide a script based on a computer algebra system that explains how the math from the previous paragraph can be extended to fifth order. The logic of the script can be extended to all higher orders.

If we circle back to eqn. (6) we see that the leading term should come from a suitable Riemann solver and be capable of entropy enforcement. In previous works (Jiang, Shu and Zhang [32], [33]) the second derivative term in eqn. (6) was obtained by a suitably high order accurate finite difference approximation of the second derivative of the pointwise, zone-centered fluxes. Now realize that, since we are using a finite difference approximation for the second derivative, it can introduce spurious oscillations into the scheme if the solution is non-smooth on the mesh. Designing an automatic strategy that switches off the higher order terms when the solution is non-smooth is not easy. In the past, different switches had to be tailor-made for different PDEs (Zheng, Shu, Qiu [61]). While positivity-based switches can be used for the Euler equations to accomplish a suppression of spurious oscillations (Gao *et al*. [25]), such switches have not even been designed for other more complicated PDEs. Furthermore, for certain classes of PDEs that don't have a scalar pressure, it may prove very difficult or impossible to design such switches. These considerations have limited the broad-based adoption of AFD-WENO schemes. It is the fundamental reason why AFD-WENO has not gained much traction in the literature. In this paper we realize that a different type of WENO interpolation strategy can be designed and applied at zone boundaries which can naturally be used to control the spurious oscillations when the solution is non-smooth. It is our hope that this invention helps in the widespread adoption of AFD-WENO schemes.

**II.b) AFD-WENO Formulated Up To Ninth Order of Accuracy**

We are now ready to describe the AFD-WENO fluxes at various orders in the ensuing paragraphs. At each order, we will obtain a numerical flux at each zone boundary which will consist of the resolved flux from the Riemann solver plus a correction stemming from the previously developed insight. Denoting the numerical flux at $x_{i+1/2}$ by $\mathbf{F}_{i+1/2}^{num}$, the discrete in space but continuous in time update takes the simple form

$$\partial_t \mathbf{U}_i = -\frac{1}{\Delta x}\left(\mathbf{F}_{i+1/2}^{num} - \mathbf{F}_{i-1/2}^{num}\right) \quad . \tag{8}$$

Because this is a finite difference scheme, the multidimensional extension is just done dimension-by-dimension.

At third order, the AFD-WENO scheme has a numerical flux given by

$$\mathbf{F}^{num}_{i+1/2} = \mathbf{F}^*\left(\hat{\mathbf{U}}^-_{i+1/2}, \hat{\mathbf{U}}^+_{i+1/2}\right) - \left\{\frac{1}{24}(\Delta x)^2 \left[\partial_x^2 \mathbf{F}\right]_{i+1/2}\right\}. \quad (9)$$

At fifth order, the AFD-WENO scheme has a numerical flux given by

$$\mathbf{F}^{num}_{i+1/2} = \mathbf{F}^*\left(\hat{\mathbf{U}}^-_{i+1/2}, \hat{\mathbf{U}}^+_{i+1/2}\right) - \left\{\frac{1}{24}(\Delta x)^2 \left[\partial_x^2 \mathbf{F}\right]_{i+1/2} - \frac{7}{5760}(\Delta x)^4 \left[\partial_x^4 \mathbf{F}\right]_{i+1/2}\right\}. \quad (10)$$

At seventh order, the AFD-WENO scheme has a numerical flux given by

$$\mathbf{F}^{num}_{i+1/2} = \mathbf{F}^*\left(\hat{\mathbf{U}}^-_{i+1/2}, \hat{\mathbf{U}}^+_{i+1/2}\right) - \left\{\begin{array}{l} \dfrac{1}{24}(\Delta x)^2 \left[\partial_x^2 \mathbf{F}\right]_{i+1/2} - \dfrac{7}{5760}(\Delta x)^4 \left[\partial_x^4 \mathbf{F}\right]_{i+1/2} \\ + \dfrac{31}{967680}(\Delta x)^6 \left[\partial_x^6 \mathbf{F}\right]_{i+1/2} \end{array}\right\}. \quad (11)$$

At ninth order, the AFD-WENO scheme has a numerical flux given by

$$\mathbf{F}^{num}_{i+1/2} = \mathbf{F}^*\left(\hat{\mathbf{U}}^-_{i+1/2}, \hat{\mathbf{U}}^+_{i+1/2}\right) - \left\{\begin{array}{l} \dfrac{1}{24}(\Delta x)^2 \left[\partial_x^2 \mathbf{F}\right]_{i+1/2} - \dfrac{7}{5760}(\Delta x)^4 \left[\partial_x^4 \mathbf{F}\right]_{i+1/2} \\ + \dfrac{31}{967680}(\Delta x)^6 \left[\partial_x^6 \mathbf{F}\right]_{i+1/2} - \dfrac{127}{154828800}(\Delta x)^8 \left[\partial_x^8 \mathbf{F}\right]_{i+1/2} \end{array}\right\}. \quad (12)$$

Now notice that all the higher order derivatives of the flux in eqns. (9) to (12) involve even derivatives of the flux evaluated at the zone boundaries. It is certainly technically possible to design a WENO interpolation that takes as its input the fluxes evaluated at the zone centers and returns as its output the even derivatives of the fluxes at the zone boundaries. However, it would lead to smaller stencils if one could pre-compute some of the derivatives so as to reduce the order of the derivatives at the zone boundaries. To that end, we can realize that

$$\partial_x \mathbf{F} = \mathbf{A}\partial_x \mathbf{U} \;\; ; \;\; \partial_x^2 \mathbf{F} = \partial_x(\mathbf{A}\partial_x \mathbf{U}) \;\; ; \;\; \partial_x^4 \mathbf{F} = \partial_x^3(\mathbf{A}\partial_x \mathbf{U}) \text{ and so on with } \mathbf{A} \equiv \frac{\partial \mathbf{F}}{\partial \mathbf{U}}. \quad (13)$$

Since the characteristic matrix "$\mathbf{A}$" can easily be evaluated at zone centers, and since $\partial_x \mathbf{U}$ can be evaluated at the zone centers as a byproduct of the WENO interpolation, it is easy to evaluate $(\mathbf{A}\partial_x \mathbf{U})$ at zone centers. These can be incorporated into eqns. (9) to (12). The result is a discrete in space but continuous in time update that can be written as

$$\partial_t \mathbf{U}_i = -\frac{1}{\Delta x}\left\{\mathbf{F}^*\left(\hat{\mathbf{U}}_{i+1/2}^-, \hat{\mathbf{U}}_{i+1/2}^+\right) - \mathbf{F}^*\left(\hat{\mathbf{U}}_{i-1/2}^-, \hat{\mathbf{U}}_{i-1/2}^+\right)\right\}$$

$$-\frac{1}{\Delta x}\left\{\begin{bmatrix} -\frac{1}{24}(\Delta x)^2 \left[\partial_x(\mathbf{A}\partial_x \mathbf{U})\right]_{i+1/2} + \frac{7}{5760}(\Delta x)^4 \left[\partial_x^3(\mathbf{A}\partial_x \mathbf{U})\right]_{i+1/2} \\ -\frac{31}{967680}(\Delta x)^6 \left[\partial_x^5(\mathbf{A}\partial_x \mathbf{U})\right]_{i+1/2} + \frac{127}{154828800}(\Delta x)^8 \left[\partial_x^7(\mathbf{A}\partial_x \mathbf{U})\right]_{i+1/2} \end{bmatrix} \\ -\begin{bmatrix} -\frac{1}{24}(\Delta x)^2 \left[\partial_x(\mathbf{A}\partial_x \mathbf{U})\right]_{i-1/2} + \frac{7}{5760}(\Delta x)^4 \left[\partial_x^3(\mathbf{A}\partial_x \mathbf{U})\right]_{i-1/2} \\ -\frac{31}{967680}(\Delta x)^6 \left[\partial_x^5(\mathbf{A}\partial_x \mathbf{U})\right]_{i-1/2} + \frac{127}{154828800}(\Delta x)^8 \left[\partial_x^7(\mathbf{A}\partial_x \mathbf{U})\right]_{i-1/2} \end{bmatrix}\right\} \quad (14)$$

Notice that eqn. (14) is still in conservation form and therefore it should be able to capture shock locations accurately. The black terms in the above equation yield a third order scheme. In that case, the derivatives of $(\mathbf{A}\partial_x \mathbf{U})$ have to be evaluated at the zone boundaries with a WENO scheme that is at least second order accurate. If the red terms are also included, in addition to the black terms, the scheme becomes fifth order accurate. In that case, all the derivatives of $(\mathbf{A}\partial_x \mathbf{U})$ have to be evaluated at the zone boundaries with a WENO scheme that is at least fourth order accurate. If the blue terms are also included, in addition to the black and red terms, we get a seventh order scheme. In that case, all the derivatives of $(\mathbf{A}\partial_x \mathbf{U})$ have to be evaluated at the zone boundaries with a WENO scheme that is at least sixth order accurate. If the magenta terms are included, in addition to the black, red and blue terms, we get a ninth order scheme. In that case, all the derivatives of $(\mathbf{A}\partial_x \mathbf{U})$ have to be evaluated at the zone boundaries with a WENO scheme that is at least eighth order accurate. Such WENO schemes are documented in Section IV. Notice, therefore, that our strategy of first using the original WENO interpolant to obtain $(\mathbf{A}\partial_x \mathbf{U})$ at zone centers allows us to decrease the order of the WENO interpolation that is needed at zone boundaries, and that can be a useful cost savings.

We finish this Section by enumerating some observations about the uses of this new class of schemes. We obviously cannot develop and demonstrate all the uses of the AFD-WENO schemes in this one paper, but our observations can serve as a signpost for future work:-

**1)** Certain problems can have extreme physics. In such circumstances it is valuable to make physics-based shock detectors (see Balsara [3], [5], [7]). Such discontinuity indicators could be directly responsive to the presence of very strong shocks, in which case the higher order finite differencing of the $\left(\mathbf{A}\partial_x\mathbf{U}\right)$-dependent terms in eqn. (14) can be strongly suppressed in the vicinity of very strong shocks.

**2)** Unlike classical FD-WENO schemes, a smooth flux like the LLF flux is not needed; even monotone fluxes can be used in an AFD-WENO scheme (see Jiang, Shu and Zhang [32], [33]). The ADF-WENO schemes presented here can be used with any Riemann solver, including positivity preserving ones. Proofs of positivity preservation usually rely on hybridizing a lower order scheme which is provably positivity preserving with a higher order scheme that may not be provably positivity preserving. The current schemes all reveal that there is a clear split in the numerical flux between the part $\mathbf{F}^*\left(\hat{\mathbf{U}}_{i+1/2}^-, \hat{\mathbf{U}}_{i+1/2}^+\right)$ which comes from the Riemann solver and the part that comes from the higher order derivatives of the fluxes. This could have advantages in advancing proofs for positivity preservation for FD-WENO schemes.

**3)** It is possible to design Riemann solvers that preserve stationary linearly degenerate discontinuities. Schemes that can preserve stationary linearly degenerate discontinuities have clear uses in certain circumstances, like well-balancing in the presence of gravitational forces (Käppeli and Mishra [34], Berberich *et al*. [13], Grosheintz-Laval and Käppeli [27], Käppeli [35]). When a stationary contact discontinuity is present, the discontinuity indicators will all become zero so that we have $\mathbf{F}_{i+1/2}^{num} \to \mathbf{F}^*\left(\hat{\mathbf{U}}_{i+1/2}^-, \hat{\mathbf{U}}_{i+1/2}^+\right)$. It is, therefore, easy to see that if the underlying Riemann solver preserves stationary contact discontinuities, the entire scheme will do the same.

**4)** The AFD-WENO scheme is especially adept at preserving the free stream condition on curvilinear meshes, as shown by (Jiang, Shu and Zhang [32], [33]). These ideas should be very useful in extending the present methods to curvilinear meshes.

**5)** If analytical procedures are used, there are limits to the kinds of curvilinear meshes that can be accessed by an AFD-WENO scheme, as shown by Merriman [46]. Scripts like the one in Appendix A could be used to extend the range of curvilinear meshes that can be handled by AFD-WENO schemes.

**6)** In Balsara *et al*. [12] we have presented high order FD-WENO schemes that can handle hyperbolic PDE systems with non-conservative products. However, the solution vector in many such systems has some flux conservative components and some components that are non-conservative. The present AFD-WENO will be extended to yield AFD-WENO approaches that are flux conservative when they need to be conservative and yet handle non-conservative products. This will be a very novel contribution that dramatically extends the applicability of AFD-WENO-based methodologies, and it will be developed in a subsequent paper that is fully dedicated to systems with non-conservative products.

**7)** AFD-WENO requires one WENO-based interpolation of the zone-centered variables. It also requires another WENO interpolation at the zone boundaries. While finite difference WENO schemes are inherently designed to have low computational complexity, the expensive part of the algorithm is in the WENO reconstruction/interpolation. Classical FD-WENO also involves two reconstruction steps applied to the left-going and right-going fluxes. Consequently, it is fair to say that AFD-WENO and FD-WENO have similar computational complexities. Therefore, the selling point of a well-designed AFD-WENO algorithm relative to its classical FD-WENO counterpart is that it offers greater flexibility at the same computational cost.

This completes a broad-based formulation and discussion of the AFD-WENO schemes. In the next two sections we will document details of the WENO-based interpolation that enable us to obtain efficient implementation of AFD-WENO schemes.

**III) WENO-AO Interpolation at Several Orders for AFD-WENO Schemes**

Weighted Essentially Non-Oscillatory with Adaptive Order (WENO-AO) (Balsara, Garain and Shu [8], hereafter BGS16, Arbogast *et al*. [2], Kumar and Chandrashekar [36], [37], Balsara *et al*. [10], Boscheri and Balsara [17]) is a multiresolution strategy for carrying out WENO reconstruction just like WENO-ZQ (Zhu and Qiu [62]) and subsequent multiresolution WENO by Zhu and Shu [63]. The above authors realized that one can make a non-linear hybridization between a large, centered, very high accuracy stencil (or stencils) and a lower order WENO scheme that is, nevertheless, very stable. This yields a class of adaptive order WENO schemes, which we call WENO-AO (for Adaptive Order). The large, centered, higher order WENO stencil(s) give the reconstruction strategy its high order accuracy when a smooth solution is present on that large stencil. The lower order WENO scheme is meant to stabilize the method when the larger, higher order stencil has a non-smooth solution on it. The lower order scheme can be a third order CWENO scheme (Levy, Puppo and Russo [40], Cravero and Semplice [22], Semplice Coco and Russo [51]) which uses three piecewise quadratic polynomials (BGS16). CWENO makes a good lower order scheme because it is very robust and can, nevertheless, capture extrema. The lower order scheme can also be comprised of two piecewise linear polynomials (Zhu and Qiu [62]). One can even make allowance for the lowest order scheme to go all the way down to first order accuracy (Zhu and Shu [63]) when all stencils have a non-smooth solution. All the intuitive insights for WENO-AO reconstruction described in this paragraph will also go over for WENO-AO interpolation described in the rest of this section.

In this Section we provide sufficient background on *pointwise* WENO-AO interpolation. Please note that *pointwise* WENO-AO interpolation seeks to match the pointwise zone-centered values. This is quite different from the classical WENO reconstruction which seeks to match zone averages. This Section is intended to make this paper self-contained so that anyone can implement the schemes that we will describe in this paper. The insights developed here are crucial to the formulation of the interpolated WENO-AO scheme that is described in subsequent sections. In BGS16 tremendous detail was provided on efficient WENO reconstruction when Legendre polynomials were used for the reconstruction process. In that paper we also showed that our use of Legendre polynomials gives us the very desirable advantage that the smoothness indicators at all orders can be written compactly as the sum of perfect squares. Here we will only restrict ourselves to the WENO-AO interpolation techniques that are relevant to this particular paper. We use Legendre polynomials that span the interval $[-1/2, 1/2]$. (In Balsara *et al*. [6] and Balsara, Samantaray and Subramanian [11] we showed that it is very favorable to cast WENO interpolation and reconstruction in terms of Legendre polynomials.) They are given by

$$L_0(x) = 1 \;;\; L_1(x) = x \;;\; L_2(x) = x^2 - \frac{1}{12} \;;\; L_3(x) = x^3 - \frac{3}{20}x \;;$$

$$L_4(x) = x^4 - \frac{3}{14}x^2 + \frac{3}{560} \;;\; L_5(x) = x^5 - \frac{5}{18}x^3 + \frac{5}{336}x \;;$$

$$L_6(x) = x^6 - \frac{15}{44}x^4 + \frac{5}{176}x^2 - \frac{5}{14784} \;; \quad (15)$$

$$L_7(x) = x^7 - \frac{21}{52}x^5 + \frac{105}{2288}x^3 - \frac{35}{27456}x \;;$$

$$L_8(x) = x^8 - \frac{7}{15}x^6 + \frac{7}{104}x^4 - \frac{7}{2288}x^2 + \frac{7}{329472}.$$

We let "$r$" denote the order of accuracy of the interpolation; for example, an interpolation that is only based on $L_0(x) = 1$, $L_1(x) = x$ and $L_2(x) = x^2 - 1/12$ corresponds to $r = 3$. Our use of Legendre polynomials for the interpolation ensures that our smoothness indicators can still be expressed as the sum of perfect squares. For the rest of this Section we assume that we are dealing with a uniform mesh with zones that are mapped to the unit interval $[-1/2, 1/2]$. Because we are describing a finite difference scheme, we assume that the zone-centered variables are collocated pointwise at the zone centers, please see Fig. 1.

With these preliminaries in place, we describe the *interpolation* strategies for WENO-AO(3), WENO-AO(5,3), WENO-AO(7,3), WENO-AO(7,5,3) and WENO-AO(9,3) schemes in Sub-Sections 3.1, 3.2, 3.3, 3.4 and 3.5 respectively. The reader who wants corresponding narrative for WENO reconstruction may consult BGS16.

### III.a) WENO-AO(3) Interpolation

We will draw upon the $r = 3$ CWENO as our lowest order scheme because it is robust, inexpensive and can capture extrema. We focus on the interpolation problem in a zone labeled by a subscript "0". Consider the neighboring zone-centered variables $\{u_{-2}, u_{-1}, u_0, u_1, u_2\}$. A third order interpolation over the zone labeled "0" can be carried out by using the left-biased $r = 3$ stencil $S_1^{r3}$, the centered $r = 3$ stencil $S_2^{r3}$ and the right-biased $r = 3$ stencil $S_3^{r3}$ that rely on the variables $\{u_{-2}, u_{-1}, u_0\}$, $\{u_{-1}, u_0, u_1\}$ and $\{u_0, u_1, u_2\}$ respectively. These stencils are shown as the magenta, green and blue stencils in Fig. 1. At each of the zone centers we want the *pointwise* WENO interpolation to match the *pointwise* values at the zone centers. In other words, please note that this is not the traditional WENO reconstruction that is used in Jiang and Shu [31] or Balsara and Shu [8]! In this paper we label our stencils with a superscript that denotes the $r^{th}$ order of the polynomial. The subscripts 1, 2, 3 will denote the three stencils under consideration. The third order polynomial resulting from the WENO interpolation also carries the same subscripting convention. The $i^{th}$ interpolated polynomial $P_i^{r3}(x)$ corresponding to stencil $S_i^{r3}$ is then expressed as

$$P_i^{r3}(x) = u_{pt} + u_x \, L_1(x) + u_{x2} \, L_2(x) \ . \tag{16}$$

The left-biased $r = 3$ stencil $S_1^{r3}$ gives

$$u_{pt} = (25 \, u_0 - 2 \, u_{-1} + u_{-2})/24 \ ; \ u_x = (3 \, u_0 - 4 \, u_{-1} + u_{-2})/2 \ ;$$
$$u_{x2} = (u_0 - 2 \, u_{-1} + u_{-2})/2. \tag{17}$$

The centered $r = 3$ stencil $S_2^{r3}$ gives

$$u_{pt} = (22 \, u_0 + u_{-1} + u_1)/24 \ ; \ u_x = (u_1 - u_{-1})/2 \ ; \ u_{x2} = (-2 \, u_0 + u_{-1} + u_1)/2. \tag{18}$$

The right-biased $r = 3$ stencil $S_3^{r3}$ gives

$$u_{pt} = (25 \, u_0 - 2 \, u_1 + u_2)/24 \ ; \ u_x = (-3 \, u_0 + 4 \, u_1 - u_2)/2 \ ;$$
$$u_{x2} = (u_0 - 2 \, u_1 + u_2)/2 \tag{19}$$

Please compare eqns. (17), (18) and (19) in this paper to eqns. (17), (18) and (19) in BGS16 to appreciate the differences between interpolation and reconstruction. Such pointwise WENO interpolations have been explored at lower orders in Sebastian and Shu [50], Carlini, Ferretti and Russo [19] and Shu [55] and in what follows we give many higher order extensions. The smoothness indicator for each of the three stencils is unchanged and can then be written in a very compact form which is a sum of two squares as

$$\beta^{r3} = \left(u_x\right)^2 + \frac{13}{3}\left(u_{x2}\right)^2. \tag{20}$$

In designing a WENO-AO scheme at third order, which can graciously degrade to second or even first order at discontinuities, we wish to make a non-linearly hybridized interpolation using stencils $S_1^{r3}$, $S_2^{r3}$ and $S_3^{r3}$. The $r = 3$ WENO-AO interpolation is described by one parameter $\gamma_{Lo} \in (0,1)$. The linear weights for the stencils $S_1^{r3}$, $S_2^{r3}$ and $S_3^{r3}$ are given by

$$\gamma_1^{r3} = (1-\gamma_{Lo})/2 \quad ; \quad \gamma_2^{r3} = \gamma_{Lo} \quad ; \quad \gamma_3^{r3} = (1-\gamma_{Lo})/2. \tag{21}$$

Notice that for the linear weights we have, $\gamma_1^{r3} + \gamma_2^{r3} + \gamma_3^{r3} = 1$. Typically, we set $\gamma_{Lo} \in [0.8, 0.95]$. We see that when the central stencil $S_2^{r3}$ is smooth we want most of our interpolation to come from the central stencil because it is the most stable choice when the solution is smooth. However, when a suitable comparison of the smoothness indicators shows that the central stencil is non-smooth, we want most (or all) of our interpolation to be weighted towards either the left-biased stencil or the right-biased stencil depending on which one of the two is the smoothest one.

To avoid loss of accuracy at critical points (Borges *et al.*, [18]) we use the smoothness indicators to define

$$\tau = \frac{1}{2}\left( \left|\beta_2^{r3} - \beta_1^{r3}\right| + \left|\beta_2^{r3} - \beta_3^{r3}\right| \right) \tag{22}$$

where $\beta_1^{r3}$, $\beta_2^{r3}$ and $\beta_3^{r3}$ are the smoothness indicators for the three third order stencils. Here $\beta_2^{r3}$ comes from the centered third order stencil. The unnormalized non-linear weights are given by

$$w_1^{r3} = \gamma_1^{r3}\left(1 + \tau^2/\left(\beta_1^{r3} + \varepsilon\right)^2\right) \quad ; \quad w_2^{r3} = \gamma_2^{r3}\left(1 + \tau^2/\left(\beta_2^{r3} + \varepsilon\right)^2\right) \quad ; \\ w_3^{r3} = \gamma_3^{r3}\left(1 + \tau^2/\left(\beta_3^{r3} + \varepsilon\right)^2\right). \tag{23}$$

The normalization of the non-linear weights is given by

$$\bar{w}_1^{r3} = w_1^{r3}/\left(w_1^{r3} + w_2^{r3} + w_3^{r3}\right) \quad ; \quad \bar{w}_2^{r3} = w_2^{r3}/\left(w_1^{r3} + w_2^{r3} + w_3^{r3}\right) \quad ; \quad \bar{w}_3^{r3} = w_3^{r3}/\left(w_1^{r3} + w_2^{r3} + w_3^{r3}\right). \tag{24}$$

The non-linearly hybridized third order accurate WENO-AO interpolation is given by

$$P^{AO(3)}(x) = \bar{w}_1^{r3} P_1^{r3}(x) + \bar{w}_2^{r3} P_2^{r3}(x) + \bar{w}_3^{r3} P_3^{r3}(x). \tag{25}$$

**III.b) WENO-AO(5,3) Interpolation**

Recall that WENO-AO(5,3) from BGS16 consists of a non-linear hybridization between a large, centered, fifth order stencil denoted by $S_3^{r5}$ that relies on the variables $\{u_{-2}, u_{-1}, u_0, u_1, u_2\}$ and the three smaller stencils described above. The three smaller stencils are still shown as the magenta, green and blue stencils in Fig. 1; whereas the larger stencil is shown as the red stencil in Fig. 1. The fifth order accurate polynomial is given by

$$P_3^{r5}(x) = u_{pt} + u_x L_1(x) + u_{x2} L_2(x) + u_{x3} L_3(x) + u_{x4} L_4(x) \tag{26}$$

where the coefficients of the above polynomial are given by

$$\begin{aligned}
u_{pt} &= (5178\, u_0 + 308\, u_{-1} - 17\, u_{-2} + 308\, u_1 - 17\, u_2)/5760, \\
u_x &= (-154\, u_{-1} + 17\, u_{-2} + 154\, u_1 - 17\, u_2)/240, \\
u_{x2} &= (-402\, u_0 + 212\, u_{-1} - 11\, u_{-2} + 212\, u_1 - 11\, u_2)/336, \\
u_{x3} &= (2\, u_{-1} - u_{-2} - 2\, u_1 + u_2)/12, \\
u_{x4} &= (6\, u_0 - 4\, u_{-1} + u_{-2} - 4\, u_1 + u_2)/24.
\end{aligned} \tag{27}$$

Please compare eqn. (27) in this paper to eqn. (29) in BGS16 to appreciate the differences between interpolation and reconstruction. The corresponding smoothness indicator is unchanged and it is given by

$$\beta_3^{r5} = (u_x + u_{x3}/10)^2 + \frac{13}{3}\left(u_{x2} + \frac{123}{455} u_{x4}\right)^2 + \frac{781}{20}(u_{x3})^2 + \frac{1421461}{2275}(u_{x4})^2 . \tag{28}$$

Our task is to eventually make a non-linear hybridization between the larger stencil $S_3^{r5}$ and the smaller stencils $S_1^{r3}$, $S_2^{r3}$ and $S_3^{r3}$. The method is described by two parameters $\gamma_{Hi} \in (0,1)$ and $\gamma_{Lo} \in (0,1)$. The linear weights for the centered $r=5$ stencil $S_3^{r5}$ and the three $r=3$ stencils $S_1^{r3}$, $S_2^{r3}$ and $S_3^{r3}$ are given by

$$\gamma_3^{r5} = \gamma_{Hi} \;;\; \gamma_1^{r3} = (1-\gamma_{Hi})(1-\gamma_{Lo})/2 \;;\; \gamma_2^{r3} = (1-\gamma_{Hi})\gamma_{Lo} \;;\; \gamma_3^{r3} = (1-\gamma_{Hi})(1-\gamma_{Lo})/2 . \tag{29}$$

Notice that for the linear weights we have, $\gamma_3^{r5} + \gamma_1^{r3} + \gamma_2^{r3} + \gamma_3^{r3} = 1$. Typically, we set $\gamma_{Hi} \in [0.8, 0.95]$ and $\gamma_{Lo} \in [0.8, 0.95]$. These numbers themselves give us a glimpse of what is afoot. When a suitable comparison of the smoothness indicators shows that the large central stencil $S_3^{r5}$ is smooth, we want most (or all) of our interpolation to come from the large central stencil. However, when a suitable comparison of the smoothness indicators shows that the large central stencil is non-smooth, we want most (or all) of our interpolation to be weighted towards a very stable, third order accurate, extrema-preserving $r=3$ CWENO-type interpolation.

We now describe the process of obtaining the non-linearly hybridized weights. To avoid loss of accuracy at critical points (Borges *et al.*, [18]) we use the smoothness indicators to define

$$\tau = \frac{1}{3}\left(\left|\beta_3^{r5} - \beta_1^{r3}\right| + \left|\beta_3^{r5} - \beta_2^{r3}\right| + \left|\beta_3^{r5} - \beta_3^{r3}\right|\right) \tag{30}$$

where $\beta_3^{r5}$ is the smoothness indicator for the large centered 5$^{th}$ order stencil and $\beta_1^{r3}$, $\beta_2^{r3}$ and $\beta_3^{r3}$ are the smoothness indicators for the three smaller third order stencils. Here $\beta_2^{r3}$ comes from the centered third order stencil. The unnormalized non-linear weights are given by

$$w_3^{r5} = \gamma_3^{r5}\left(1 + \tau^2\big/\left(\beta_3^{r5} + \varepsilon\right)^2\right) \quad ; \quad w_1^{r3} = \gamma_1^{r3}\left(1 + \tau^2\big/\left(\beta_1^{r3} + \varepsilon\right)^2\right) \quad ;$$
$$w_2^{r3} = \gamma_2^{r3}\left(1 + \tau^2\big/\left(\beta_2^{r3} + \varepsilon\right)^2\right) \quad ; \quad w_3^{r3} = \gamma_3^{r3}\left(1 + \tau^2\big/\left(\beta_3^{r3} + \varepsilon\right)^2\right) \quad . \tag{31}$$

The normalization of the non-linear weights is given by

$$\overline{w}_3^{r5} = w_3^{r5}\big/\left(w_3^{r5} + w_1^{r3} + w_2^{r3} + w_3^{r3}\right) \quad ; \quad \overline{w}_1^{r3} = w_1^{r3}\big/\left(w_3^{r5} + w_1^{r3} + w_2^{r3} + w_3^{r3}\right) \quad ;$$
$$\overline{w}_2^{r3} = w_2^{r3}\big/\left(w_3^{r5} + w_1^{r3} + w_2^{r3} + w_3^{r3}\right) \quad ; \quad \overline{w}_3^{r3} = w_3^{r3}\big/\left(w_3^{r5} + w_1^{r3} + w_2^{r3} + w_3^{r3}\right) \quad . \tag{32}$$

We denote the interpolated polynomial for WENO-AO(5,3) as $P^{AO(5,3)}(x)$. Our task in this paragraph is to describe the construction of the order-preserving, non-linearly hybridized, fifth order polynomial $P^{AO(5,3)}(x)$. The non-linear weights should be combined in such a way that when all the smoothness indicators seem to have almost similar values then only the higher order scheme is obtained. From eqns. (31) and (32) realize that when the four smoothness measures associated with these four stencils have closely similar values, we have $\overline{w}_3^{r5} \to \gamma_3^{r5}$, $\overline{w}_1^{r3} \to \gamma_1^{r3}$, $\overline{w}_2^{r3} \to \gamma_2^{r3}$ and $\overline{w}_3^{r3} \to \gamma_3^{r3}$. We then require that when the limits specified by the previous sentence are attained, we have $P^{AO(5,3)}(x) \to P_3^{r5}(x)$. This is achieved by the following definition

$$\begin{aligned}\mathrm{P}^{\mathrm{AO}(5,3)}(x) &= \frac{\overline{w}_3^{r5}}{\gamma_3^{r5}}\left(\mathrm{P}_3^{r5}(x) - \gamma_1^{r3}\,\mathrm{P}_1^{r3}(x) - \gamma_2^{r3}\,\mathrm{P}_2^{r3}(x) - \gamma_3^{r3}\,\mathrm{P}_3^{r3}(x)\right) \\ &\quad + \overline{w}_1^{r3}\,\mathrm{P}_1^{r3}(x) + \overline{w}_2^{r3}\,\mathrm{P}_2^{r3}(x) + \overline{w}_3^{r3}\,\mathrm{P}_3^{r3}(x)\end{aligned} \qquad (33)$$

In the limit where the larger stencil has a non-smooth solution but the smaller stencils have smooth solutions, we have $\overline{w}_3^{r5} \ll \overline{w}_1^{r3}$ or $\overline{w}_3^{r5} \ll \overline{w}_2^{r3}$ or $\overline{w}_3^{r5} \ll \overline{w}_3^{r3}$ and also $\overline{w}_1^{r3} \to \gamma_1^{r3}$, $\overline{w}_2^{r3} \to \gamma_2^{r3}$ and $\overline{w}_3^{r3} \to \gamma_3^{r3}$. This ensures that the smoothest of the $r=3$ CWENO-type stencils will be sought out by the interpolation polynomial. Notice that the non-linear hybridization that we sought at the beginning of this paragraph has been found via eqn. (33). This completes our description of $\mathrm{P}^{\mathrm{AO}(5,3)}(x)$.

### III.c) WENO-AO(7,3) Interpolation

This Sub-section is a small variation on the previous one. Recall that WENO-AO(7,3) from BGS16 consists of a non-linear hybridization between a large, centered, seventh order stencil denoted by $S_4^{r7}$ that relies on the variables $\{u_{-3}, u_{-2}, u_{-1}, u_0, u_1, u_2, u_3\}$ and the three smaller CWENO-type stencils. The seventh order accurate polynomial is given by

$$\mathrm{P}_4^{r7}(x) = \mathrm{u}_{pt} + \mathrm{u}_x\,L_1(x) + \mathrm{u}_{x2}\,L_2(x) + \mathrm{u}_{x3}\,L_3(x) + \mathrm{u}_{x4}\,L_4(x) + \mathrm{u}_{x5}\,L_5(x) + \mathrm{u}_{x6}\,L_6(x). \qquad (34)$$

The large central stencil gives

$$\begin{aligned}\mathrm{u}_{pt} &= (\,862564\,u_0 + 57249\,u_{-1} - 5058\,u_{-2} + 367\,u_{-3} + 57249\,u_1 - 5058\,u_2 + 367\,u_3\,)/967680\,; \\ \mathrm{u}_x &= (\,-19083\,u_{-1} + 3372\,u_{-2} - 367\,u_{-3} + 19083\,u_1 - 3372\,u_2 + 367\,u_3\,)/26880\,; \\ \mathrm{u}_{x2} &= (\,-34380\,u_0 + 18625\,u_{-1} - 1546\,u_{-2} + 111\,u_{-3} + 18625\,u_1 - 1546\,u_2 + 111\,u_3\,)/26880\,; \\ \mathrm{u}_{x3} &= (\,229\,u_{-1} - 140\,u_{-2} + 17\,u_{-3} - 229\,u_1 + 140\,u_2 - 17\,u_3\,)/864\,; \\ \mathrm{u}_{x4} &= (\,2404\,u_0 - 1671\,u_{-1} + 510\,u_{-2} - 41\,u_{-3} - 1671\,u_1 + 510\,u_2 - 41\,u_3\,)/6336\,; \\ \mathrm{u}_{x5} &= (\,-5\,u_{-1} + 4\,u_{-2} - u_{-3} + 5\,u_1 - 4\,u_2 + u_3\,)/240\,; \\ \mathrm{u}_{x6} &= (\,-20\,u_0 + 15\,u_{-1} - 6\,u_{-2} + u_{-3} + 15\,u_1 - 6\,u_2 + u_3\,)/720.\end{aligned}$$

$$(35)$$

Please compare eqn. (35) in this paper to eqn. (40) in BGS16 to appreciate the differences between interpolation and reconstruction. The corresponding smoothness indicator is unchanged and is given by

$$\beta_4'^{r7} = (u_x + u_{x3}/10 + u_{x5}/126)^2 + \frac{13}{3}\left(u_{x2} + \frac{123}{455}u_{x4} + \frac{85}{2002}u_{x6}\right)^2$$
$$+ \frac{781}{20}\left(u_{x3} + \frac{26045}{49203}u_{x5}\right)^2 + \frac{1421461}{2275}\left(u_{x4} + \frac{81596225}{93816426}u_{x6}\right)^2 \quad (36)$$
$$+ \frac{21520059541}{1377684}(u_{x5})^2 + \frac{15510384942580921}{27582029244}(u_{x6})^2 .$$

The linear weights are analogous to eqn. (29), and are given by replacing the superscript "r5" in eqn. (29) with "r7". We define "$\tau$" analogously to eqn. (30), but we replace the superscript "r5" with "r7". The unnormalized non-linear weights are given by eqn. (31) with a replacement of the superscript "r5" with "r7". The normalization of the non-linear weights is given by eqn. (32) with a replacement of the superscript "r5" with "r7". The interpolated polynomial $P^{AO(7,3)}(x)$ is then given by an expression that is very analogous to eqn. (33) and is given by

$$P^{AO(7,3)}(x) = \frac{\overline{w}_4^{r7}}{\gamma_4^{r7}}\left(P_4^{r7}(x) - \gamma_1^{r3} P_1^{r3}(x) - \gamma_2^{r3} P_2^{r3}(x) - \gamma_3^{r3} P_3^{r3}(x)\right) \quad (37)$$
$$+ \overline{w}_1^{r3} P_1^{r3}(x) + \overline{w}_2^{r3} P_2^{r3}(x) + \overline{w}_3^{r3} P_3^{r3}(x) .$$

This completes our description of $P^{AO(7,3)}(x)$.

### III.d) WENO-AO(7,5,3) Interpolation

Notice that WENO-AO(7,3) from BGS16 produces an abrupt transition from seventh to third order. Because FD-WENO schemes have odd orders of accuracy, we would like to have a scheme that produces a smoother transition from seventh to fifth order and then from fifth to third order. WENO-AO(7,5,3) was designed in BGS16 to rectify that fact. Arbogast et al. [2] and Kumar and Chandrashekar [37] suggested a slightly better arrangement of the linear weights than the original BGS16, but the intent is the same. Here we describe the WENO-AO(7,5,3) interpolation. Notice that it will use the large size, seventh order stencil described in eqns. (34), (35) and (36) but it will also use the intermediate sized, fifth order stencil described in eqns. (26), (27) and (28). In addition, it will of course use the small third order stencils in eqns. (17) to (19).

The linear weights are given by:-

$$\gamma_4^{r7} = \gamma_{Hi} \quad ; \quad \gamma_3^{r5} = (1-\gamma_{Hi})\gamma_{Avg} \quad ;$$

$$\gamma_1^{r3} = (1-\gamma_{Hi})(1-\gamma_{Avg})(1-\gamma_{Lo})/2 \quad ; \quad \gamma_2^{r3} = (1-\gamma_{Hi})(1-\gamma_{Avg})\gamma_{Lo} \quad ; \quad \gamma_3^{r3} = (1-\gamma_{Hi})(1-\gamma_{Avg})(1-\gamma_{Lo})/2$$

(38)

with $\gamma_{Hi} \in [0.8, 0.95]$, $\gamma_{Avg} \in [0.85, 0.95]$ and $\gamma_{Lo} \in [0.85, 0.95]$ being the usual choices. Notice that for the linear weights we have, $\gamma_4^{r7} + \gamma_3^{r5} + \gamma_1^{r3} + \gamma_2^{r3} + \gamma_3^{r3} = 1$. To avoid loss of accuracy at critical points we use

$$\tau = \frac{1}{4}\left(\left|\beta_4^{r7} - \beta_3^{r5}\right| + \left|\beta_4^{r7} - \beta_1^{r3}\right| + \left|\beta_4^{r7} - \beta_2^{r3}\right| + \left|\beta_4^{r7} - \beta_3^{r3}\right|\right) \tag{39}$$

Where $\beta_4^{r7}$ is the smoothness indicator for the large centered 7$^{th}$ order stencil, $\beta_3^{r5}$ is the smoothness indicator for the somewhat smaller centered 5$^{th}$ order stencil and $\beta_1^{r3}$, $\beta_2^{r3}$ and $\beta_3^{r3}$ are the smoothness indicators for the three smaller third order stencils. Here $\beta_2^{r3}$ comes from the centered third order stencil. The unnormalized non-linear weights are

$$w_4^{r7} = \gamma_4^{r7}\left(1 + \tau^3/(\beta_4^{r7} + \varepsilon)^2\right) \quad ; \quad w_3^{r5} = \gamma_3^{r5}\left(1 + \tau^3/(\beta_3^{r5} + \varepsilon)^2\right) \quad ;$$

$$w_1^{r3} = \gamma_1^{r3}\left(1 + \tau^3/(\beta_1^{r3} + \varepsilon)^2\right) \quad ; \quad w_2^{r3} = \gamma_2^{r3}\left(1 + \tau^3/(\beta_2^{r3} + \varepsilon)^2\right) \quad ; \quad w_3^{r3} = \gamma_3^{r3}\left(1 + \tau^3/(\beta_3^{r3} + \varepsilon)^2\right)$$

(40)

The normalization of the non-linear weights is given by

$$\bar{w}_4^{r7} = w_4^{r7}/\left(w_4^{r7} + w_3^{r5} + w_1^{r3} + w_2^{r3} + w_3^{r3}\right) \quad ; \quad \bar{w}_3^{r5} = w_3^{r5}/\left(w_4^{r7} + w_3^{r5} + w_1^{r3} + w_2^{r3} + w_3^{r3}\right) \quad ;$$

$$\bar{w}_1^{r3} = w_1^{r3}/\left(w_4^{r7} + w_3^{r5} + w_1^{r3} + w_2^{r3} + w_3^{r3}\right) \quad ; \quad \bar{w}_2^{r3} = w_2^{r3}/\left(w_4^{r7} + w_3^{r5} + w_1^{r3} + w_2^{r3} + w_3^{r3}\right) \quad ;$$

$$\bar{w}_3^{r3} = w_3^{r3}/\left(w_4^{r7} + w_3^{r5} + w_1^{r3} + w_2^{r3} + w_3^{r3}\right)$$

(41)

We finally get the non-linearly hybridized interpolation polynomial $P^{AO(7,5,3)}(x)$ as

$$P^{AO(7,5,3)}(x) = \frac{\bar{w}_4^{r7}}{\gamma_4^{r7}}\left(P_4^{r7}(x) - \gamma_3^{r5}P_3^{r5}(x) - \gamma_1^{r3}P_1^{r3}(x) - \gamma_2^{r3}P_2^{r3}(x) - \gamma_3^{r3}P_3^{r3}(x)\right)$$
$$+ \bar{w}_3^{r5}P_3^{r5}(x) + \bar{w}_1^{r3}P_1^{r3}(x) + \bar{w}_2^{r3}P_2^{r3}(x) + \bar{w}_3^{r3}P_3^{r3}(x)$$

(42)

This completes our description of $P^{AO(7,5,3)}(x)$.

### III.e) WENO-AO(9,3) Interpolation

A WENO-AO(9,5,3) system can be just as easily created from the style of thinking presented in this and in the previous Sub-section. Therefore, we just document WENO-AO(9,3) here. The ninth order accurate polynomial is given by

$$P_5^{r9}(x) = u_{pt} + u_x\, L_1(x) + u_{x2}\, L_2(x) + u_{x3}\, L_3(x) + u_{x4}\, L_4(x) + u_{x5}\, L_5(x) + u_{x6}\, L_6(x) \\ + u_{x7}\, L_7(x) + u_{x8}\, L_8(x) \quad . \tag{43}$$

The large central stencil gives

$$\begin{aligned}
u_{pt} &= (\,412080590\, u_0 + 29039624\, u_{-1} - 3207892\, u_{-2} + 399032\, u_{-3} - 27859\, u_{-4} + 29039624\, u_1 \\
&\quad - 3207892\, u_2 + 399032\, u_3 - 27859\, u_4\,)\,/\,464486400\;; \\
u_x &= (\,-7259906\, u_{-1} + 1603946\, u_{-2} - 299274\, u_{-3} + 27859\, u_{-4} + 7259906\, u_1 - 1603946\, u_2 \\
&\quad + 299274\, u_3 - 27859\, u_4\,)\,/\,9676800\;; \\
u_{x2} &= (\,-28194190\, u_0 + 15523184\, u_{-1} - 1610524\, u_{-2} + 198224\, u_{-3} - 13789\, u_{-4} + 15523184\, u_1 \\
&\quad - 1610524\, u_2 + 198224\, u_3 - 13789\, u_4\,)\,/\,21288960\;; \\
u_{x3} &= (\,747682\, u_{-1} - 512722\, u_{-2} + 106218\, u_{-3} - 10223\, u_{-4} - 747682\, u_1 + 512722\, u_2 \\
&\quad - 106218\, u_3 + 10223\, u_4\,)\,/\,2280960\;; \\
u_{x4} &= (\,3007170\, u_0 - 2143448\, u_{-1} + 733204\, u_{-2} - 100584\, u_{-3} + 7243\, u_{-4} - 2143448\, u_1 \\
&\quad + 733204\, u_2 - 100584\, u_3 + 7243\, u_4\,)\,/\,6589440\;; \\
u_{x5} &= (\,-2974\, u_{-1} + 2662\, u_{-2} - 918\, u_{-3} + 101\, u_{-4} + 2974\, u_1 - 2662\, u_2 + 918\, u_3 - 101\, u_4\,) \\
&\quad /\,74880\;; \\
u_{x6} &= (\,-4430\, u_0 + 3424\, u_{-1} - 1532\, u_{-2} + 352\, u_{-3} - 29\, u_{-4} + 3424\, u_1 - 1532\, u_2 + 352\, u_3 \\
&\quad - 29\, u_4\,)\,/\,86400\;; \\
u_{x7} &= (\,14\, u_{-1} - 14\, u_{-2} + 6\, u_{-3} - u_{-4} - 14\, u_1 + 14\, u_2 - 6\, u_3 + u_4\,)\,/\,10080\;; \\
u_{x8} &= (\,70\, u_0 - 56\, u_{-1} + 28\, u_{-2} - 8\, u_{-3} + u_{-4} - 56\, u_1 + 28\, u_2 - 8\, u_3 + u_4\,)\,/\,40320.
\end{aligned}$$

(44)

The corresponding smoothness indicator is given by

$$\beta_5^{r9} = (u_x + u_{x3}/10 + u_{x5}/126 + u_{x7}/1716)^2 + \frac{13}{3}\left(u_{x2} + \frac{123}{455}u_{x4} + \frac{85}{2002}u_{x6} + \frac{29}{5577}u_{x8}\right)^2$$

$$+ \frac{781}{20}\left(u_{x3} + \frac{26045}{49203}u_{x5} + \frac{8395}{60918}u_{x7}\right)^2$$

$$+ \frac{1421461}{2275}\left(u_{x4} + \frac{81596225}{93816426}u_{x6} + \frac{618438835}{1829420307}u_{x8}\right)^2$$

$$+ \frac{21520059541}{1377684}\left(u_{x5} + \frac{722379670131}{559521548066}u_{x7}\right)^2$$

$$+ \frac{15510384942580921}{27582029244}\left(u_{x6} + \frac{5423630339859998294}{3024525063803279595}u_{x8}\right)^2$$

$$+ \frac{12210527897166191835083}{443141066068272}(u_{x7})^2$$

$$+ \frac{7550936809810378933608373 1407561}{42818201328263029226415}(u_{x8})^2.$$

(45)

The linear weights are analogous to eqn. (29), and are given by replacing the superscript "r5" in eqn. (29) with "r9". We define "$\tau$" analogously to eqn. (30), but we replace the superscript "r5" with "r9". The unnormalized non-linear weights are given by eqn. (31) with a replacement of the superscript "r5" with "r9". The normalization of the non-linear weights is given by eqn. (32) with a replacement of the superscript "r5" with "r9". The interpolated polynomial $P^{AO(9,3)}(x)$ is then given by

$$P^{AO(9,3)}(x) = \frac{\overline{w}_5^{r9}}{\gamma_5^{r9}}\left(P_5^{r9}(x) - \gamma_1^{r3}P_1^{r3}(x) - \gamma_2^{r3}P_2^{r3}(x) - \gamma_3^{r3}P_3^{r3}(x)\right)$$
$$+ \overline{w}_1^{r3}P_1^{r3}(x) + \overline{w}_2^{r3}P_2^{r3}(x) + \overline{w}_3^{r3}P_3^{r3}(x).$$

(46)

This completes our description of $P^{AO(9,3)}(x)$; and it also completes this Section. Analogously to the $P^{AO(7,5,3)}(x)$ in eqn. (42), we can also construct $P^{AO(9,5,3)}(x)$ in order to obtain a scheme that graciously relinquishes order of accuracy as the solution loses smoothness.

**IV) A New Type of WENO-AO Interpolation that is Applicable to Zone Boundaries**

In the course of designing AFD-WENO algorithms we will need strategies for obtaining high order derivatives of the $(\mathbf{A}\partial_x\mathbf{U})$ variable at the zone boundaries. The FD-WENO algorithm uses zone-centered point values of the conserved variables, and the same interpolant can be used to obtain $\partial_x\mathbf{U}$ at the zone centers. The characteristic matrix "$\mathbf{A}$" can be obtained at the zone centers by making a pointwise evaluation. The term $(\mathbf{A}\partial_x\mathbf{U})$ can, therefore, be evaluated at the zone centers using these point values of the solution vector and its derivative at the zone center. These zone-centered values of $(\mathbf{A}\partial_x\mathbf{U})$ serve as inputs for a new type of WENO interpolation that we describe in this Section. This new type of WENO interpolation takes those zone-centered values of $(\mathbf{A}\partial_x\mathbf{U})$ as inputs and provides nonlinearly hybridized higher order derivatives of $(\mathbf{A}\partial_x\mathbf{U})$ at the zone boundaries as its outputs. This is shown in Fig. 2. The ensuing Sub-sections show how this is done for AFD-WENO schemes of increasing order of accuracy.

## IV.a) Third Order WENO-AO Interpolation of Higher Order Derivatives at Zone Boundaries

We focus on the third order WENO-AO interpolation problem at a zone boundary labeled by a subscript "1/2". Let the origin be centered at this zone boundary. From eqn. (14) we see that we will need a first derivative of $(\mathbf{A}\partial_x\mathbf{U})$ at each zone boundary. Consider the symmetrically placed neighboring zone-centered $(\mathbf{A}\partial_x\mathbf{U})$ variables which we denote in the formulae that follow by $\{f_{-1}, f_0, f_1, f_2\}$. A third order interpolation at the zone boundary labeled "1/2" can be carried out by using the left-biased $S_1^{r3}$ stencil and the right-biased $S_2^{r3}$ stencil that rely on the variables $\{f_{-1}, f_0, f_1\}$ and $\{f_0, f_1, f_2\}$ respectively. These are shown by the magenta and blue stencils in Fig. 2 and we denote them as $S_{zb;1}^{r3}$ and $S_{zb;2}^{r3}$ respectively. Please also see Fig. 2 for the indexing convention. At each of the zone centers we want the *pointwise* WENO interpolation to match the *pointwise* values of $(\mathbf{A}\partial_x\mathbf{U})$ at the zone centers. The $i^{\text{th}}$ interpolated polynomial $P_{zb;i}^{r3}(x)$ corresponding to stencil $S_{zb;i}^{r3}$ is then expressed as

$$P_{zb;i}^{r3}(x) = f_{pt} + f_x\, L_1(x) + f_{x2}\, L_2(x) \ . \tag{47}$$

The left-biased $r = 3$ stencil $S_{zb;1}^{r3}$ gives

$$f_{pt} = (8 f_0 - f_{-1} + 5 f_1) / 12 \ ; \ f_x = (f_1 - f_0) \ ; \ f_{x2} = (f_1 - 2 f_0 + f_{-1}) / 2 \ . \tag{48}$$

The right-biased $r = 3$ stencil $S_{zb;2}^{r3}$ gives

$$f_{pt} = (8 f_1 - f_2 + 5 f_0) / 12 \ ; \ f_x = (f_1 - f_0) \ ; \ f_{x2} = (f_2 - 2 f_1 + f_0) / 2 \ . \tag{49}$$

The form of the smoothness indicators is unchanged, and we denote the two smoothness indicators for the two stencils as $\beta_1^{r3}$ and $\beta_2^{r3}$. The stencils can be given equal linear weights and non-linearly combined as follows

$$\tau = \left| \beta_2^{r3} - \beta_1^{r3} \right| \ ; \ w_1^{r3} = 0.5 \left( 1 + \tau^2 / \left( \beta_1^{r3} + \varepsilon \right)^2 \right) \ ; \ w_2^{r3} = 0.5 \left( 1 + \tau^2 / \left( \beta_2^{r3} + \varepsilon \right)^2 \right) \tag{50}$$

The non-linear weights are normalized as

$$\overline{w}_1^{r3} = w_1^{r3} / \left( w_1^{r3} + w_2^{r3} \right) \ ; \ \overline{w}_2^{r3} = w_2^{r3} / \left( w_1^{r3} + w_2^{r3} \right) \ . \tag{51}$$

The non-linearly hybridized third order accurate WENO interpolation of $(\mathbf{A} \partial_x \mathbf{U})$ at the zone boundary is given by

$$P_{zb}^{AO(3)}(x) = \overline{w}_1^{r3} \ P_1^{r3}(x) + \overline{w}_2^{r3} \ P_2^{r3}(x) \ . \tag{52}$$

It is worth noting that eqn. (52) will be used in a third order scheme only for the purpose of obtaining a first derivative of $(\mathbf{A} \partial_x \mathbf{U})$ at the zone boundary of interest. By examining eqns. (48) and (49) it is easy to intuit that if the two stencils have first derivatives of the same sign then the interpolation in eqn. (52) biases us towards the stencil with the smaller first derivative. If the two stencils have first derivatives with opposite signs, then they will try to cancel one another and the first derivative with the smaller value will be favored. From this discussion we get the essential insight that our novel WENO-AO interpolation is naturally stabilizing and can itself act like a discontinuity indicator for choosing first derivatives of the zone-centered $(\mathbf{A} \partial_x \mathbf{U})$ variable at the zone boundaries.

## IV.b) Fourth Order WENO-AO Interpolation of Higher Order Derivatives at Zone Boundaries

Realize from Fig. 2 and eqn. (14) that a fifth order AFD-WENO scheme will need to obtain first and third derivatives of $\left(\mathbf{A}\partial_x \mathbf{U}\right)$ at the zone boundaries. The smallest stencil that is symmetrical around zone boundary "$i+1/2$" and contains third derivatives is the large fourth order stencil shown in Fig. 2; we call it stencil $S_{zb;c}^{r4}$. We focus on the fourth order WENO-AO interpolation problem in a zone boundary labeled by a subscript "1/2". Let the origin be centered at this zone boundary. Consider the symmetrically placed neighboring zone-centered flux variables $\{f_{-1}, f_0, f_1, f_2\}$. The small left-biased stencil, the small right-biased stencil and the large centered stencil are shown by the magenta, blue and red stencils respectively in Fig. 2 and we denote them by $S_{zb;1}^{r3}$, $S_{zb;2}^{r3}$ and $S_{zb;c}^{r4}$. The smaller stencils were described in the previous Sub-section. A non-linear hybridization at the zone boundary labeled "1/2" can be carried out by using the small left-biased $S_1^{r3}$ stencil, the small right-biased $S_2^{r3}$ stencil and the large centered $S_{zb;c}^{r4}$ stencil that rely on the variables $\{f_{-1}, f_0, f_1\}$, $\{f_0, f_1, f_2\}$ and $\{f_{-1}, f_0, f_1, f_2\}$ respectively. The interpolated polynomial $P_{zb;c}^{r4}(x)$ that corresponds to $S_{zb;c}^{r4}$ is given by

$$P_{zb;c}^{r4}(x) = f_{pt} + f_x\ L_1(x) + f_{x2}\ L_2(x) + f_{x3}\ L_3(x) \tag{53}$$

where we obtain

$$\begin{aligned}
f_{pt} &= (13\ f_0 - f_{-1} + 13\ f_1 - f_2)/24\ ; \\
f_x &= (-63\ f_0 + f_{-1} + 63\ f_1 - f_2)/60\ ; \\
f_{x2} &= (-f_0 + f_{-1} - f_1 + f_2)/4\ ; \\
f_{x3} &= (3\ f_0 - f_{-1} - 3\ f_1 + f_2)/6.
\end{aligned} \tag{54}$$

The corresponding smoothness indicator can be obtained from BGS16 and we denote it by $\beta_c^{r4}$. The three stencils that we are considering can be non-linearly combined as

$$\gamma_c^{r4} = \gamma_{Hi} \quad ; \quad \gamma_1^{r3} = (1-\gamma_{Hi})0.5 \quad ; \quad \gamma_2^{r3} = (1-\gamma_{Hi})0.5 \quad ; \quad \tau = \frac{1}{2}\left(\left|\beta_c^{r4} - \beta_1^{r3}\right| + \left|\beta_c^{r4} - \beta_2^{r3}\right|\right) \quad ;$$

$$w_c^{r4} = \gamma_c^{r4}\left(1 + \tau^2 / \left(\beta_c^{r4} + \varepsilon\right)^2\right) \quad ; \quad w_1^{r3} = \gamma_1^{r3}\left(1 + \tau^2 / \left(\beta_1^{r3} + \varepsilon\right)^2\right) \quad ; \quad w_2^{r3} = \gamma_2^{r3}\left(1 + \tau^2 / \left(\beta_2^{r3} + \varepsilon\right)^2\right) .$$

(55)

The normalization of the non-linear weights is given by

$$\overline{w}_c^{r4} = w_c^{r4} / \left(w_c^{r4} + w_1^{r3} + w_2^{r3}\right) \quad ; \quad \overline{w}_1^{r3} = w_1^{r3} / \left(w_c^{r4} + w_1^{r3} + w_2^{r3}\right) \quad ; \quad \overline{w}_2^{r3} = w_2^{r3} / \left(w_c^{r4} + w_1^{r3} + w_2^{r3}\right) .$$

(56)

The non-linearly hybridized fourth order accurate WENO interpolation of the fluxes at the zone boundary is given by

$$P_{zb}^{AO(4,3)}(x) = \frac{\overline{w}_c^{r4}}{\gamma_c^{r4}}\left(P_c^{r4}(x) - \gamma_1^{r3}P_1^{r3}(x) - \gamma_2^{r3}P_2^{r3}(x)\right) + \overline{w}_1^{r3} P_1^{r3}(x) + \overline{w}_2^{r3} P_2^{r3}(x) . \qquad (57)$$

It is now easy to see that when $\left(\overline{w}_c^{r4}/\gamma_c^{r4}\right) \to 1$, the first and third derivatives of the above polynomial will be obtained exclusively from the fourth order polynomial $P_c^{r4}(x)$. As a result, they will be evaluated with the desired level of accuracy. When $\left(\overline{w}_c^{r4}/\gamma_c^{r4}\right) \to 0$, the third derivative will be zero and the first derivative will also be of lower order of accuracy since it will be obtained exclusively from $\overline{w}_1^{r3} P_1^{r3}(x) + \overline{w}_2^{r3} P_2^{r3}(x)$. When even the first derivatives from the two smaller stencils have opposite signs, the evaluation of the first derivative will be even further suppressed. Therefore, just as in the previous Sub-section, we see here that our novel WENO-AO interpolation at zone boundaries is naturally stabilizing and can itself act like a discontinuity indicator for choosing first and third derivatives of the zone-centered $(\mathbf{A}\partial_x\mathbf{U})$ variable at the zone boundaries.

### IV.c) Sixth Order WENO-AO Interpolation of Higher Order Derivatives at Zone Boundaries

Realize from eqn. (14) that a seventh order AFD-WENO scheme will need to obtain first, third and fifth derivatives of $(\mathbf{A}\partial_x\mathbf{U})$ at the zone boundaries. The smallest stencil that is symmetrical around zone boundary "$i+1/2$" and contains fifth derivatives is a large sixth order stencil; let us denote this as $S^{r6}_{zb;c}$. We focus on the sixth order WENO-AO interpolation problem in a zone boundary labeled by a subscript "1/2". Let the origin be centered at this zone boundary. Consider the symmetrically placed neighboring zone-centered variables $\{f_{-2}, f_{-1}, f_0, f_1, f_2, f_3\}$. Realize, first off, that if all the points are used, this will result in a sixth order accurate interpolation. As before, we will non-linearly hybridize between the small left-biased $S^{r3}_1$ stencil, the small right-biased $S^{r3}_2$ stencil and the large centered $S^{r6}_{zb;c}$ stencil that rely on the variables $\{f_{-1}, f_0, f_1\}$, $\{f_0, f_1, f_2\}$ and $\{f_{-2}, f_{-1}, f_0, f_1, f_2, f_3\}$ respectively. The smaller stencils were described in the previous Sub-section. The interpolated polynomial $P^{r6}_{zb;c}(x)$ that corresponds to the $S^{r6}_{zb;c}$ stencil is given by

$$P^{r6}_{zb;c}(x) = f_{pt} + f_x\ L_1(x) + f_{x2}\ L_2(x) + f_{x3}\ L_3(x) + f_{x4}\ L_4(x) + f_{x5}\ L_5(x) \tag{58}$$

where we obtain

$$\begin{aligned}
f_{pt} &= (802\ f_0 - 93.0\ f_{-1} + 11\ f_{-2} + 802\ f_1 - 93\ f_2 + 11\ f_3)/1440\ ; \\
f_x &= (-1794\ f_0 + 43\ f_{-1} - 3\ f_{-2} + 1794\ f_1 - 43\ f_2 + 3\ f_3)/1680\ ; \\
f_{x2} &= (-29\ f_0 + 33\ f_{-1} - 4.0\ f_{-2} - 29\ f_1 + 33\ f_2 - 4\ f_3)/84\ ; \\
f_{x3} &= (37\ f_0 - 14\ f_{-1} + f_{-2} - 37\ f_1 + 14\ f_2 - f_3)/54\ ; \\
f_{x4} &= (2\ f_0 - 3\ f_{-1} + f_{-2} + 2\ f_1 - 3\ f_2 + f_3)/48\ ; \\
f_{x5} &= (-10\ f_0 + 5\ f_{-1} - f_{-2} + 10\ f_1 - 5\ f_2 + f_3)/120\ .
\end{aligned} \tag{59}$$

The corresponding smoothness indicator can be obtained from BGS16 and we denote it by $\beta^{r6}_c$. The three stencils that we are considering can be non-linearly combined in a fashion that is analogous to eqn. (55) where we replace the superscript "r4" with "r6". The normalization of the non-linear weights is given by an equation that is analogous to eqn. (56) with the superscript "r4" replaced with "r6". The non-linearly hybridized sixth order accurate WENO interpolation of the fluxes at the zone boundary is given by

$$P^{AO(6,3)}_{zb}(x) = \frac{\overline{w}^{r6}_c}{\gamma^{r6}_c}\left(P^{r6}_c(x) - \gamma^{r3}_1 P^{r3}_1(x) - \gamma^{r3}_2 P^{r3}_2(x)\right) + \overline{w}^{r3}_1\ P^{r3}_1(x) + \overline{w}^{r3}_2\ P^{r3}_2(x)\ . \tag{60}$$

It is now easy to see that when $\left(\bar{w}_c^{r6}/\gamma_c^{r6}\right) \to 1$, the first, third and fifth derivatives of the above polynomial will be obtained exclusively from the sixth order polynomial $P_c^{r6}(x)$. As a result, they will be very accurate. When $\left(\bar{w}_c^{r6}/\gamma_c^{r6}\right) \to 0$, the third and fifth derivatives will be zero and the first derivative will also be of lower order of accuracy since it will be obtained exclusively from $\bar{w}_1^{r3} P_1^{r3}(x) + \bar{w}_2^{r3} P_2^{r3}(x)$. When even the first derivatives from the two smaller stencils have opposite signs, the evaluation of the first derivative will be even further suppressed. Therefore, just as in the previous Sub-sections, we see here that our novel WENO-AO interpolation at zone boundaries is naturally stabilizing and can itself act like a discontinuity indicator for choosing first, third and fifth derivatives of the zone-centered $(\mathbf{A}\partial_x\mathbf{U})$ variable at the zone boundaries.

## IV.d) Eighth Order WENO-AO Interpolation of Higher Order Derivatives at Zone Boundaries

Realize from eqn. (14) that a ninth order AFD-WENO scheme will need to obtain first, third, fifth and seventh derivatives of $(\mathbf{A}\partial_x\mathbf{U})$ at the zone boundaries. We would like to use this section to demonstrate that it is possible to design a multiresolution method at zone boundaries. We focus on the eighth order WENO-AO interpolation problem in a zone boundary labeled by a subscript "1/2". Let the origin be centered at this zone boundary. Consider the symmetrically placed neighboring zone-centered flux variables $\{f_{-3}, f_{-2}, f_{-1}, f_0, f_1, f_2, f_3, f_4\}$; let us denote this as $S_{zb;c}^{r8}$. Realize, first off, that if all the points are used, this will result in an eighth order interpolation. An eighth order interpolation at the zone boundary labeled "1/2" can be carried out by using the small left-biased $S_1^{r3}$ stencil, the small right-biased $S_2^{r3}$ stencil, the sixth order centered $S_{zb;c}^{r6}$ stencil and the eighth order centered $S_{zb;c}^{r8}$ stencil that rely on the variables $\{f_{-1}, f_0, f_1\}$, $\{f_0, f_1, f_2\}$, $\{f_{-2}, f_{-1}, f_0, f_1, f_2, f_3\}$ and $\{f_{-3}, f_{-2}, f_{-1}, f_0, f_1, f_2, f_3, f_4\}$ respectively. The coefficients for the $S_1^{r3}$ stencil, the $S_2^{r3}$ stencil and the $S_{zb;c}^{r6}$ stencil have already been given in eqns. (48), (49) and (59) respectively. The interpolated polynomial $P_{zb;c}^{r8}(x)$ that corresponds to the $S_{zb;c}^{r8}$ stencil is given by

$$P_{zb;c}^{r8}(x) = f_{pt} + f_x\ L_1(x) + f_{x2}\ L_2(x) + f_{x3}\ L_3(x) + f_{x4}\ L_4(x) + f_{x5}\ L_5(x) + f_{x6}\ L_6(x) + f_{x7}\ L_7(x)$$

(61)

where we obtain

$$f_{pt} = (68323\ f_0 - 9531\ f_{-1} + 1879\ f_{-2} - 191\ f_{-3} + 68323\ f_1 - 9531\ f_2 + 1879\ f_3 - 191\ f_4)/120960 ;$$
$$f_x = (-325685\ f_0 + 9399\ f_{-1} - 1093\ f_{-2} + 79\ f_{-3} + 325685\ f_1 - 9399\ f_2 + 1093\ f_3 - 79\ f_4)/302400 ;$$
$$f_{x2} = (-2655\ f_0 + 3243\ f_{-1} - 655\ f_{-2} + 67\ f_{-3} - 2655\ f_1 + 3243\ f_2 - 655\ f_3 + 67\ f_4)/6720 ;$$
$$f_{x3} = (111365\ f_0 - 45171\ f_{-1} + 5377\ f_{-2} - 391\ f_{-3} - 111365\ f_1 + 45171\ f_2 - 5377\ f_3 + 391\ f_4)/142560 ;$$
$$f_{x4} = (449\ f_0 - 729\ f_{-1} + 317\ f_{-2} - 37\ f_{-3} + 449\ f_1 - 729\ f_2 + 317\ f_3 - 37\ f_4)/6336 ;$$
$$f_{x5} = (-2645\ f_0 + 1431\ f_{-1} - 373\ f_{-2} + 31\ f_{-3} + 2645\ f_1 - 1431\ f_2 + 373\ f_3 - 31\ f_4)/18720 ;$$
$$f_{x6} = (-5\ f_0 + 9\ f_{-1} - 5\ f_{-2} + f_{-3} - 5\ f_1 + 9\ f_2 - 5\ f_3 + f_4)/1440 ;$$
$$f_{x7} = (35\ f_0 - 21\ f_{-1} + 7\ f_{-2} - f_{-3} - 35\ f_1 + 21\ f_2 - 7\ f_3 + f_4)/5040 .$$

(62)

The corresponding smoothness indicator can be obtained from BGS16 and we denote it by $\beta_c^{r8}$. The four stencils that we are considering can be non-linearly combined as

$$\gamma_c^{r8} = \gamma_{Hi} \ ; \ \gamma_c^{r6} = (1-\gamma_{Hi})\gamma_{Avg} \ ; \ \gamma_1^{r3} = (1-\gamma_{Hi})(1-\gamma_{Avg})0.5 \ ; \ \gamma_2^{r3} = (1-\gamma_{Hi})(1-\gamma_{Avg})0.5 \ ;$$
$$\tau = \frac{1}{3}\left(\left|\beta_c^{r8} - \beta_c^{r6}\right| + \left|\beta_c^{r8} - \beta_1^{r3}\right| + \left|\beta_c^{r8} - \beta_2^{r3}\right|\right) \ ;$$
$$w_c^{r8} = \gamma_c^{r6}\left(1 + \tau^4/(\beta_c^{r8} + \varepsilon)^2\right) \ ; \ w_c^{r6} = \gamma_c^{r6}\left(1 + \tau^4/(\beta_c^{r6} + \varepsilon)^2\right) \ ;$$
$$w_1^{r3} = \gamma_1^{r3}\left(1 + \tau^4/(\beta_1^{r3} + \varepsilon)^2\right) \ ; \ w_2^{r3} = \gamma_2^{r3}\left(1 + \tau^4/(\beta_2^{r3} + \varepsilon)^2\right) .$$

(63)

The normalization of the non-linear weights is given by

$$\bar{w}_c^{r8} = w_c^{r8}/(w_c^{r8} + w_c^{r6} + w_1^{r3} + w_2^{r3}) \ ; \ \bar{w}_c^{r6} = w_c^{r6}/(w_c^{r8} + w_c^{r6} + w_1^{r3} + w_2^{r3}) \ ;$$
$$\bar{w}_1^{r3} = w_1^{r3}/(w_c^{r8} + w_c^{r6} + w_1^{r3} + w_2^{r3}) \ ; \ \bar{w}_2^{r3} = w_2^{r3}/(w_c^{r8} + w_c^{r6} + w_1^{r3} + w_2^{r3}) .$$

(64)

The non-linearly hybridized eighth order accurate WENO interpolation of the fluxes at the zone boundary is given by

$$P_{zb}^{AO(8,6,3)}(x) = \frac{\bar{w}_c^{r8}}{\gamma_c^{r8}}\left(P_c^{r8}(x) - \gamma_c^{r6}P_c^{r6}(x) - \gamma_1^{r3}P_1^{r3}(x) - \gamma_2^{r3}P_2^{r3}(x)\right)$$
$$+ \bar{w}_c^{r6}P_c^{r6}(x) + \bar{w}_1^{r3}P_1^{r3}(x) + \bar{w}_2^{r3}P_2^{r3}(x) .$$

(65)

We see once again that our novel WENO-AO interpolation at zone boundaries is naturally stabilizing and can itself act like a discontinuity indicator for choosing first, third, fifth and seventh derivatives of the zone-centered $(\mathbf{A}\partial_x \mathbf{U})$ variable at the zone boundaries. Moreover, because of the inclusion of $P_c^{r6}(x)$ in eqn. (65), our WENO method can graciously degrade the order of accuracy of the derivatives that are used.

**IV.e) Multiresolution WENO Interpolation of Higher Order Derivatives at Zone Boundaries**

In this paper we have described WENO-AO interpolation. However, all the formulae developed here can be seamlessly used also for Multiresolution WENO interpolation. When interpolating to zone boundaries, the smallest stencil is a piecewise linear stencil that is formed from $\{f_0, f_1\}$. Eqns. (53) and (54) then provide us with a fourth order strategy for interpolating to zone boundaries. Eqns. (58) and (59) then provide us with a sixth order strategy for interpolating to zone boundaries. Eqns. (61) and (62) then provide us with an eighth order strategy for interpolating to zone boundaries. The non-linear hybridization between stencils is the same as in Zhu and Shu [63].

**V) Pointwise Implementation of Our AFD-WENO Scheme for Conservation Laws**

Now that the above discussions are understood, we provide a pointwise implementation of our AFD-WENO scheme for treating non-linear hyperbolic PDEs that can be written in the form of conservation laws. We realize that the update equation, i.e. eqn. (14), has a lot of terms. The optimal sequence of steps given below is designed so that at the end of each step we catalogue the parts of eqn. (14) that are in hand so that we can finally assemble the entire equation. The pointwise implementation of our AFD-WENO scheme into a numerical code goes according to the following steps:-

**1)** We start with the mesh function as shown in Fig. 1. This means that at each zone center $x_i$ we have a pointwise value for the conserved variable $\mathbf{U}_i$. AFD-WENO schemes are always implemented in dimension-by-dimension fashion, so we only describe one of the dimensional updates here.

**2)** From the conserved variables in each zone, obtain the primitive variables. Use them both to obtain the normalized right and left eigenvectors in the conserved variables.

**3)** As shown in Fig. 1, we use the WENO-AO algorithm from Section III. That Section includes all closed form expressions that are needed for WENO interpolation in one dimension. This consists of making a non-linear hybridization between a large high order accurate stencil and smaller lower order accurate stencils. We use WENO-AO or multiresolution WENO interpolation in characteristic variables. As a result, the neighboring zones around zone "$i$" are projected into the characteristic space of zone "$i$". The extent of these neighboring zones depends on the desired order of the scheme. The fifth order case is explicitly shown in Fig. 1. (Expanding the large stencil by one zone on either side adds two further orders of accuracy.) Once the variables in the neighboring zones around zone "$i$" are projected into the characteristic space of zone "$i$", WENO-AO interpolation is carried out in the characteristic space. Projecting the interpolation characteristic variables back into the space of right eigenvectors gives us high order accurate $\hat{\mathbf{U}}^{-}_{i+1/2}$ and $\hat{\mathbf{U}}^{+}_{i-1/2}$ within each zone "$i$", as shown in Fig. 1. Also evaluate $(\mathbf{A}\partial_x\mathbf{U})$ at the zone centers. Since this step involves projecting all the zones in all the the stencils of interest into the characteristic space of each zone "$i$" using eigenvectors, it is one of the two computationally expensive steps of the algorithm. By the end of this step we should have $\hat{\mathbf{U}}^{-}_{i+1/2}$ and $\hat{\mathbf{U}}^{+}_{i+1/2}$ at each zone boundary and $(\mathbf{A}\partial_x\mathbf{U})_i$ at each zone center.

**4)** At each zone boundary $x_{i+1/2}$, use the left and right states $\hat{\mathbf{U}}^{-}_{i+1/2}$ and $\hat{\mathbf{U}}^{+}_{i+1/2}$ to obtain the left-most and right-most going speeds of the Riemann fan; these are denoted by $S_{L;i+1/2}$ and $S_{R;i+1/2}$. Please note that at this point in the game we are not yet seeking the resolved state within the Riemann fan; that will come later.

**5)** Now, at each zone boundary $x_{i+1/2}$, we hand in the speeds $S_{L;i+1/2}$ and $S_{R;i+1/2}$ as well as the states $\hat{\mathbf{U}}^{-}_{i+1/2}$ and $\hat{\mathbf{U}}^{+}_{i+1/2}$ to the Riemann solver. This gives us the resolved state $\mathbf{U}^{*}_{i+1/2}$ and the resolved flux $\mathbf{F}^{*}_{HLLI;i+1/2}$ at each zone boundary. Please use the formulae in Section III of Balsara *et al.* (2023) if the HLLI Riemann solver is being used. By the end of this step we should have $\mathbf{U}^{*}_{i+1/2}$ and $\mathbf{F}^{*}_{HLLI;i+1/2}$ at each zone boundary.

**6)** If one wants to make a characteristic projection of the $(\mathbf{A}\partial_x\mathbf{U})$ variables, we can do that using $\mathbf{U}^*_{i+1/2}$. This could be useful in the next step. Therefore, we find the matrices of right and left eigenvectors corresponding to the resolved state $\mathbf{U}^*_{i+1/2}$ at each zone boundary "*i+1/2*". Please notice that if the HLLI Riemann solver is used, then we will naturally be constructing the left and right eigenvectors; so it is worthwhile to derive the maximum use from those eigenvectors.

**7)** Use the boundary-centered WENO-AO interpolation results from Section IV and Fig. 2 to interpolate the zone-centered $(\mathbf{A}\partial_x\mathbf{U})$ variables and their higher derivatives to the zone boundaries. (While it is physically meaningful and advantageous to project the zone-centered $(\mathbf{A}\partial_x\mathbf{U})$ variables in the eigenspace of the resolved state $\mathbf{U}^*_{i+1/2}$ at each zone boundary "*i+1/2*", one may also choose to make a component by component WENO-AO interpolation of the same variables. In our experience, this less expensive option works just as well as characteristic projection.) This gives us suitably high order derivatives of the $(\mathbf{A}\partial_x\mathbf{U})$ variables at each zone boundary. Since this step involves projecting all the zones in all the stencils of interest into the characteristic space of each zone boundary "*i+1/2*" using eigenvectors, it is the second of the two computationally expensive steps of the algorithm. By the end of this step we should have higher order derivatives like $\left[\partial_x(\mathbf{A}\partial_x\mathbf{U})\right]_{i+1/2}$, $\left[\partial_x^3(\mathbf{A}\partial_x\mathbf{U})\right]_{i+1/2}$, $\left[\partial_x^5(\mathbf{A}\partial_x\mathbf{U})\right]_{i+1/2}$ and $\left[\partial_x^7(\mathbf{A}\partial_x\mathbf{U})\right]_{i+1/2}$ (as needed) at each of the zone boundaries.

**8)** Now realize from the previous steps that we have acquired all the terms that will contribute to eqn. (14). That gives us one spatially higher order update stage of a multistage RK update strategy.

**9)** The above points have only shown one stage of the scheme. It can be coupled with an SSP-RK update strategy, say from Shu and Osher [53] or Spiteri and Ruuth [58], [59], to achieve higher order in time.

**10)** Some of the PDEs also have stiff source terms; these are usually relaxation terms that enable the system to relax to several useful physical limits. The AFD-WENO method makes it very simple to treat stiff source terms because the source terms are treated pointwise and are collocated at the exact same location as the primal variables. For this reason, when stiff source terms are present, we recommend using the Runge-Kutta IMEX methods from Pareschi and Russo [49]; see also Kupka *et al*. [38].

Notice that Steps 3 and 8 in our AFD-WENO algorithm are indeed the computationally expensive parts of the algorithm because they involve characteristic projections over large stencils. However, please compare this to the reconstruction of the left-going and right-going LLF fluxes in classical FD-WENO. That too counts as two steps where we have to make characteristic projections over large stencils. Therefore, AFD-WENO that is presented here has the same computational complexity as classical FD-WENO! Of course, the AFD-WENO presented here can be used with different types of Riemann solvers and also on curvilinear meshes, thus adding to its versatility. If multiresolution WENO interpolation is used instead of WENO-AO interpolation, there is a slight increase in computational complexity because many more stencils have to be constructed. But the stencil width is the same for both interpolation strategies, so the increase in cost is not substantial. This completes our pointwise implementation-oriented description of our AFD-WENO scheme.

## VI) Accuracy Analysis

In this section, we present several two-dimensional accuracy analyses for the Euler flow, Relativistic hydrodynamics and the 10-moment equations for rarefied gases.

For the time-update, we use a third order accurate SSP-RK scheme from Shu and Osher [53], and a fourth order scheme from Spiteri and Ruuth [58], [59]. The base level grid for all of these accuracy tests was run with a CFL of 0.4. For the spatially third order accurate scheme we use third order SSP-RK scheme, and for the fifth and seventh order accurate scheme we use fourth order time stepping. Consequently, for the fifth and seventh order scheme we reduce the time-step size as the mesh was refined so that the temporal error remain dominated by the spatial error. When a spatially fifth order scheme is used with a temporally fourth order accurate time-stepping strategy, then every doubling of the mesh requires a reduction in the timestep that goes as $\Delta t \rightarrow \Delta t (1/2)^{5/4}$. Similarly, when a spatially seventh order scheme is used with a temporally fourth order time-stepping strategy, then every doubling of the mesh requires a reduction in the timestep that goes as $\Delta t \rightarrow \Delta t (1/2)^{7/4}$.

**VI.a) Two-Dimensional Vortex Problem for Euler Flow**

In this Sub-section, we consider a two-dimensional hydrodynamic vortex. The vortex propagates on a square mesh along with being advected in the diagonal direction. The detailed set-up of the problem is given in Pao and Salas [48], and Balsara and Shu [4]. To minimize the effect of small jumps in the velocity field at the periodic boundaries we double the computational domain and stopping time for the seventh and ninth order schemes. The accuracy results for the density variable is presented in Table I. In the first half of Table I we present the accuracy results when a WENO-AO limiting process is used, where we see that the third through ninth order schemes reach their design accuracies very well. In the second half of the table we show the accuracy results when a Multiresolution WENO limiting process is used. We see that the two forms of non-linear limiting produce comparable accuracies.

**Table I shows the accuracy analysis of two-dimensional hydrodynamic vortex problem. The density variable is shown. WENO-AO results are shown for orders 3, 5, 7, 9. We also show Multiresolution WENO results for orders 3, 5, 7, 9. The lowest order scheme in the multiresolution formulation was a piecewise linear scheme with MC$_\beta$ limiter.**

| Order 3 WENO-AO | $L_1$ Error | $L_1$ Accuracy | $L_\infty$ Error | $L_\infty$ Accuracy |
|---|---|---|---|---|
| $64^2$ | 1.01690E-03 | | 1.65204E-02 | |
| $128^2$ | 1.55293E-04 | 2.71 | 2.38188E-03 | 2.79 |
| $256^2$ | 2.00234E-05 | 2.96 | 3.10946E-04 | 2.94 |
| $512^2$ | 2.50449E-06 | 3.00 | 3.90121E-05 | 2.99 |
| Order 5 WENO-AO | | | | |
| $64^2$ | 7.68939E-05 | | 8.96471E-04 | |
| $128^2$ | 2.44609E-06 | 4.97 | 4.41780E-05 | 4.34 |
| $256^2$ | 7.34858E-08 | 5.06 | 1.85043E-06 | 4.58 |
| $512^2$ | 2.80911E-09 | 4.71 | 9.58492E-08 | 4.27 |
| Order 7 WENO-AO | | | | |
| $64^2$ | 1.57282E-04 | | 2.18673E-02 | |
| $128^2$ | 4.07955E-06 | 5.27 | 4.22366E-04 | 5.69 |
| $256^2$ | 4.30990E-08 | 6.56 | 5.09831E-06 | 6.37 |
| $512^2$ | 3.54277E-10 | 6.93 | 4.05577E-08 | 6.97 |
| Order 9 WENO-AO | | | | |
| $64^2$ | 7.64417E-05 | | 6.60120E-03 | |
| $128^2$ | 8.47453E-07 | 6.50 | 1.47463E-04 | 5.48 |
| $256^2$ | 2.39272E-09 | 8.47 | 2.45157E-07 | 9.23 |
| $384^2$ | 6.95433E-11 | 8.73 | 6.93602E-09 | 8.79 |
| | | | | |
| Order 3 Multires WENO | | | | |
| $64^2$ | 2.07635E-03 | | 3.67584E-02 | |
| $128^2$ | 2.63248E-04 | 2.98 | 6.67737E-03 | 2.46 |
| $256^2$ | 3.30907E-05 | 2.99 | 1.12047E-03 | 2.58 |
| $512^2$ | 5.58111E-06 | 2.57 | 2.29457E-04 | 2.29 |
| Order 5 Multires WENO | | | | |
| $64^2$ | 8.91942E-05 | | 1.59998E-03 | |
| $128^2$ | 3.24049E-06 | 4.78 | 6.66914E-05 | 4.58 |
| $256^2$ | 1.02816E-07 | 4.98 | 2.53590E-06 | 4.72 |
| $512^2$ | 2.55432E-09 | 5.33 | 4.77709E-08 | 5.73 |
| Order 7 Multires WENO | | | | |
| $64^2$ | 1.75104E-04 | | 2.14962E-02 | |
| $128^2$ | 4.38706E-06 | 5.32 | 4.25091E-04 | 5.66 |
| $256^2$ | 4.31078E-08 | 6.67 | 5.09790E-06 | 6.38 |
| $512^2$ | 3.54349E-10 | 6.93 | 4.05499E-08 | 6.97 |
| Order 9 Multires WENO | | | | |
| $64^2$ | 6.50461E-05 | | 5.98146E-03 | |
| $128^2$ | 8.47465E-07 | 6.26 | 1.47477E-04 | 5.34 |
| $256^2$ | 2.39274E-09 | 8.47 | 2.45157E-07 | 9.23 |
| $384^2$ | 6.96344E-11 | 8.72 | 6.93627E-09 | 8.79 |

## VI.b) Two-Dimensional Vortex Problem for Relativistic Hydrodynamical Flow

In this Sub-section we consider the two-dimensional Relativistic hydrodynamics (RHD) equations from [1],[14]. For an ideal fluid, the system can be written in a conservation form for the mass density $D$, momentum density $\mathbf{M} = (M_x, M_y)$ and the total energy density $E$.

$$\frac{\partial}{\partial t}\begin{pmatrix} D \\ M_x \\ M_y \\ E \end{pmatrix} + \frac{\partial}{\partial x}\begin{pmatrix} Dv_x \\ M_x v_x + p \\ M_y v_x \\ M_x \end{pmatrix} + \frac{\partial}{\partial y}\begin{pmatrix} Dv_y \\ M_x v_y \\ M_y v_y + p \\ M_y \end{pmatrix} = 0$$

The primitive quantities $\rho$, $\mathbf{v}=(v_x,v_y)$ and $p$ are the proper mass density, fluid velocity vector and isotropic gas pressure, respectively. The expressions for the conserved quantities are given by

$$D = \Gamma\rho, \qquad M_x = \rho h \Gamma^2 v_x, \, M_y = \rho h \Gamma^2 v_y, \text{ and } E = \rho h \Gamma^2 - p.$$

Here, $\gamma$ is the specific heat constant, $h = 1 + \dfrac{\gamma}{\gamma-1}\dfrac{p}{\rho}$ is the special enthalpy and $\Gamma = \dfrac{1}{\sqrt{1-v_x^2-v_y^2}}$ is the Lorentz factor. The speed of light is assumed to be unity. The complete set of right-eigenvectors can be found in [14], and the set of left-eigenvectors can be obtained by analytically inverting the right-eigenvector matrix.

In this Sub-section we perform the accuracy study for the above system. In 2016, Balsara and Kim [9] constructed an isentropic vortex problem for the Relativistic Magneto-Hydrodynamics system (RMHD). For the construction of initial solution, an ordinary differential equation needed to be integrated numerically. For the RHD system, the analytic solution of the isentropic vortex problem with the explicit expression is derived in [41]. For the cartesian co-ordinates, the vortex problem has been used to study the accuracy of the high-order accurate entropy stable schemes in [23] where the authors have given the explicit expressions in great details in their Section 4.2, therefore we don't describe them here and use the same set up. In Table II we show the accuracy analysis for the density variable. In the first half, WENO-AO results are shown for the orders 3, 5, 7, 9 and in the lower half of the Table II Multiresolution WENO results are shown. As previosuly, we double the computational domain and stopping time for the seventh and ninth order schemes to minimize the effect of small jumps in the velocity field at the periodic boundaries. We observe that both the schemes are able to reach the design accuracy.

**Table II shows the accuracy analysis of two-dimensional Vortex Problem for Relativistic Hydrodynamical Flow. The density variable is shown. WENO-AO results are shown for orders 3, 5, 7, 9. We also show Multiresolution WENO results for orders 3, 5, 7, 9. The lowest order scheme in the multiresolution formulation was a piecewise linear scheme with MC$_\beta$ limiter.**

| Order 3 WENO-AO | $L_1$ Error | $L_1$ Accuracy | $L_\infty$ Error | $L_\infty$ Accuracy |
|---|---|---|---|---|
| $64^2$ | 2.23475E-03 | | 4.53575E-02 | |
| $128^2$ | 3.70004E-04 | 2.59 | 9.30784E-03 | 2.28 |
| $256^2$ | 4.82911E-05 | 2.94 | 1.12257E-03 | 3.05 |
| $512^2$ | 6.05348E-06 | 3.00 | 1.43072E-04 | 2.97 |
| Order 5 WENO-AO | | | | |
| $64^2$ | 2.67573E-04 | | 1.26773E-02 | |
| $128^2$ | 1.01298E-05 | 4.72 | 6.84625E-04 | 4.21 |
| $256^2$ | 2.93989E-07 | 5.11 | 2.72598E-05 | 4.65 |
| $512^2$ | 1.24462E-08 | 4.56 | 8.76169E-07 | 4.96 |
| Order 7 WENO-AO | | | | |
| $64^2$ | 7.61264E-04 | | 4.28187E-02 | |
| $128^2$ | 2.36564E-05 | 5.01 | 7.12349E-03 | 2.59 |
| $256^2$ | 2.93116E-07 | 6.33 | 1.74332E-04 | 5.35 |
| $512^2$ | 2.46174E-09 | 6.90 | 1.54570E-06 | 6.82 |
| Order 9 WENO-AO | | | | |
| $64^2$ | 9.21005E-04 | | 5.89138E-02 | |
| $128^2$ | 7.82577E-06 | 6.88 | 3.60729E-03 | 4.03 |
| $256^2$ | 2.49843E-08 | 8.29 | 1.65346E-05 | 7.77 |
| $384^2$ | 7.05366E-10 | 8.80 | 5.04467E-07 | 8.61 |
| | | | | |
| Order 3 Multires WENO | | | | |
| $64^2$ | 4.51667E-03 | | 1.21497E-01 | |
| $128^2$ | 8.25931E-04 | 2.45 | 2.60012E-02 | 2.22 |
| $256^2$ | 1.30874E-04 | 2.66 | 6.56014E-03 | 1.99 |
| $512^2$ | 2.14674E-05 | 2.61 | 1.60335E-03 | 2.03 |
| Order 5 Multires WENO | | | | |
| $64^2$ | 2.52914E-04 | | 1.00131E-02 | |
| $128^2$ | 1.24974E-05 | 4.34 | 7.96119E-04 | 3.65 |
| $256^2$ | 4.06384E-07 | 4.94 | 3.42472E-05 | 4.54 |
| $512^2$ | 1.27305E-08 | 5.00 | 9.16503E-07 | 5.22 |
| Order 7 Multires WENO | | | | |
| $64^2$ | 8.99424E-04 | | 4.20644E-02 | |
| $128^2$ | 2.45318E-05 | 5.20 | 7.12111E-03 | 2.56 |
| $256^2$ | 2.93295E-07 | 6.39 | 1.74295E-04 | 5.35 |
| $512^2$ | 2.46175E-09 | 6.90 | 1.54570E-06 | 6.82 |
| Order 9 Multires WENO | | | | |
| $64^2$ | 9.09933E-04 | | 7.71162E-02 | |
| $128^2$ | 7.82667E-06 | 6.86 | 3.60735E-03 | 4.42 |
| $256^2$ | 2.49842E-08 | 8.29 | 1.65346E-05 | 7.77 |
| $384^2$ | 7.05372E-10 | 8.80 | 5.04466E-07 | 8.61 |

**VI.c) Two-Dimensional Sinusoidal Problem for 10-Moment Rarefied Gas Flow**

In this Sub-section we consider the two-dimensional ten-moment Gaussian closure model (for the homogeneous case) which is described thoroughly in [45][15]. The system is given by

$$\frac{\partial}{\partial t}\begin{pmatrix} \rho \\ \rho v_x \\ \rho v_y \\ E_{xx} \\ E_{xy} \\ E_{yy} \end{pmatrix} + \frac{\partial}{\partial x}\begin{pmatrix} \rho v_x \\ \rho v_x^2 + p_{xx} \\ \rho v_x v_y + p_{xy} \\ \rho v_x^3 + 3 v_x p_{xx} \\ \rho v_x^2 v_y + 2 v_x p_{xy} + v_y p_{xx} \\ \rho v_x v_y^2 + v_x p_{yy} + 2 v_y p_{xy} \end{pmatrix} + \frac{\partial}{\partial y}\begin{pmatrix} \rho v_y \\ \rho v_x v_y + p_{xy} \\ \rho v_y^2 + p_{yy} \\ \rho v_y v_x^2 + v_y p_{xx} + 2 v_x p_{xy} \\ \rho v_y^2 v_x + 2 v_y p_{xy} + v_x p_{yy} \\ \rho v_y^3 + 3 v_y p_{yy} \end{pmatrix} = 0,$$

where $\rho$ is the fluid density, $\mathbf{v} = (v_x, v_y)$ is the fluid velocity vector, $\mathbf{p} = (p_{xx}, p_{xy}, p_{yy})$ is the symmetric pressure tension and $\mathbf{E} = (E_{xx}, E_{xy}, E_{yy})$ is the symmetric energy tensor. The energy tensor is obtained using the ideal equation of state

$$\mathbf{E} = \rho \mathbf{v} \otimes \mathbf{v} + \mathbf{p}.$$

The complete set of right-eigenvectors can be found in [52], and the set of left-eigenvectors can be obtained by inverting the right-eigenvector matrix.

To perform the accuracy study for the Ten-Moment model, we consider the two-dimensional sinusoidal problem, which is an extension of the one-dimensional sinusoidal problem presented in [15]; see their Section 5.1. The computational domain is given by $[-5, 5]$ with the periodic boundary conditions. The exact solution in terms of the primitive variables is given by

$$\rho(x, y, t) = 2 + \sin(2\pi(x + y - t)), \quad v_x(x, y, 0) = v_y(x, y, 0) = 1,$$
$$p_{xx}(x, y, 0) = p_{yy}(x, y, 0) = 1, \quad p_{xy}(x, y, 0) = 0.$$

We run the simulation till time $t = 0.5$ and compute the $L_1$ and $L_\infty$ error for the density variable in Table III. In the first half of the Table III, we show the accuracy results obtained using the WENO-AO scheme. In the lower half, Multiresolution WENO results are shown. We observe that both the schemes are able to reach the design accuracies for all the orders.

**Table III shows the accuracy analysis of two-dimensional Sinusoidal Problem for 10-Moment Rarefied Gas Flow. The density variable is shown. WENO-AO results are shown for orders 3, 5, 7, 9. We also show Multiresolution WENO results for orders 3, 5, 7, 9. The lowest order scheme in the multiresolution formulation was a piecewise linear scheme with $MC_\beta$ limiter.**

| Order 3 WENO-AO | $L_1$ Error | $L_1$ Accuracy | $L_\infty$ Error | $L_\infty$ Accuracy |
|---|---|---|---|---|
| $16^2$ | 3.62882E-02 | | 6.09802E-02 | |
| $32^2$ | 4.51988E-03 | 3.01 | 8.12445E-03 | 2.91 |
| $64^2$ | 5.51424E-04 | 3.04 | 1.02394E-03 | 2.99 |
| $128^2$ | 6.79346E-05 | 3.02 | 1.27881E-04 | 3.00 |
| Order 5 WENO-AO | | | | |
| $16^2$ | 1.08275E-03 | | 1.83973E-03 | |
| $32^2$ | 3.25056E-05 | 5.06 | 5.85970E-05 | 4.97 |
| $64^2$ | 9.88880E-07 | 5.04 | 1.83499E-06 | 5.00 |
| $128^2$ | 3.04681E-08 | 5.02 | 5.73462E-08 | 5.00 |
| Order 7 WENO-AO | | | | |
| $16^2$ | 3.72538E-05 | | 5.91197E-05 | |
| $32^2$ | 2.69054E-07 | 7.11 | 4.64525E-07 | 6.99 |
| $64^2$ | 2.01648E-09 | 7.06 | 3.66999E-09 | 6.98 |
| $128^2$ | 1.65398E-11 | 6.93 | 2.77980E-11 | 7.04 |
| Order 9 WENO-AO | | | | |
| $16^2$ | 5.83791E-06 | | 9.29804E-06 | |
| $32^2$ | 2.49009E-09 | 11.20 | 4.17883E-09 | 11.12 |
| $64^2$ | 4.93670E-12 | 8.98 | 7.91567E-12 | 9.04 |
| | | | | |
| Order 3 Multires WENO | | | | |
| $16^2$ | 7.15532E-02 | | 2.06001E-01 | |
| $32^2$ | 2.03270E-02 | 1.82 | 6.91685E-02 | 1.57 |
| $64^2$ | 4.39041E-03 | 2.21 | 2.30104E-02 | 1.59 |
| $128^2$ | 8.74362E-04 | 2.33 | 7.66948E-03 | 1.59 |
| Order 5 Multires WENO | | | | |
| $16^2$ | 6.29043E-03 | | 1.83363E-02 | |
| $32^2$ | 1.08591E-03 | 2.53 | 4.45898E-03 | 2.04 |
| $64^2$ | 2.91893E-05 | 5.22 | 1.94605E-04 | 4.52 |
| $128^2$ | 1.51436E-06 | 4.27 | 2.60444E-05 | 2.90 |
| Order 7 Multires WENO | | | | |
| $16^2$ | 7.89142E-05 | | 1.91619E-04 | |
| $32^2$ | 2.72952E-07 | 8.18 | 1.01114E-06 | 7.57 |
| $64^2$ | 3.25457E-09 | 6.39 | 1.95742E-08 | 5.69 |
| $128^2$ | 1.65374E-11 | 7.62 | 3.84652E-11 | 8.99 |
| Order 9 Multires WENO | | | | |
| $16^2$ | 5.84072E-06 | | 9.27601E-06 | |
| $32^2$ | 2.49074E-09 | 11.20 | 4.13097E-09 | 11.13 |
| $64^2$ | 4.93681E-12 | 8.98 | 7.89235E-12 | 9.03 |

## VII) One-Dimensional Test Problems

In this Section, we focus on one-dimensional test problems. In Sub-section 7.1 we present one-dimensional problems for the Euler flow; In Sub-section 7.1 we present one-dimensional test problems for the Relativistic Hydrodynamics flow and in Sub-section 7.3 we consider one-dimensional Riemann problems for the Ten-Moment model. For all the simulations presented in this Section we used a CFL of 0.8 with a third order SSP-RK time-stepping scheme.

## VII.a) One-dimensional test problems for Euler Flow

We present three one-dimensional test cases for the Euler Flow. The initial states, value of the specific heat constant $(\gamma)$ and final times $(t_{end})$ for these three test problems are given in Table IV. All the test problems from Table IV have computational domain $[-0.5, 0.5]$.

Test 1 is the Sod's shock-tube problem from [57], whose solution consist of a right going shock and contact discontinuity, and a left going rarefaction wave. The resulting density, velocity and pressure profiles obtained from the fifth order scheme are shown in Fig. 3. The exact solution is shown in solid lines. We observe that the numerical results for the Test 1 match with the exact solution precisely. The seventh and ninth order schemes also perform well on this problem and are not shown here. Next, we consider the Lax shock-tube problem described in [39]. The initial states for this problem are given by Test 2 in Table IV. Fig. 4 shows the density, velocity and pressure profiles obtained from the seventh order scheme. The exact solution is denoted by the solid lines. The obtained result in Fig. 4 shows that the scheme is able to capture the rarefaction, contact discontinuity and shocks efficiently. The fifth and ninth order schemes also perform well on this problem and are not shown here.

Next, we show that the method performs well for a stringent problem where the initialized profile consists of two strong shocks. The initial states for the problem is given by Test 3 in Table IV and the complete description of the problem is given in Woodward and Colella [60]. The computational domain that spans $[-0.5, 0.5]$, which has been partitioned into 1000 zones, is considered. To simulate this test case, the flattening algorithm introduced in [7] was employed. Fig. 5 shows the density, velocity and pressure profiles obtained from the ninth order scheme. A reference solution was computed using a third order AFD-WENO scheme on a 4000 zone mesh and is shown with the solid lines in Fig. 5. The obtained result shows the precise coincidence of all the profiles with the reference solution. The fifth and seventh order schemes also perform well on this problem and are not shown here.

**Table IV gives the left and right initial states, specific heat constant $(\gamma)$ and final times $(t_{end})$ of the three test problems for the Euler flow.**

|  | $x$ | $\rho$ | $v_x$ | $v_y$ | $p$ | $\gamma$ | $t_{end}$ |
|---|---|---|---|---|---|---|---|
| Test 1 (Sod-Shock) | $x < 0$ | 1 | 0 | 0 | 1 | 1.4 | 0.2 |
|  | $x > 0$ | 0.125 | 0 | 0 | 0.1 |  |  |
| Test 2 (Lax-Shock) | $x < 0$ | 0.445 | 0.698 | 0 | 3.528 | 1.4 | 0.13 |
|  | $x > 0$ | 0.5 | 0 | 0 | 0.571 |  |  |
| Test 3 (Blast wave) | $x < -0.4$ | 1 | 0 | 0 | 1000 | 1.4 | 0.038 |
|  | $-0.4 < x < 0.4$ | 1 | 0 | 0 | 0.01 |  |  |
|  | $x > 0.4$ | 1 | 0 | 0 | 100 |  |  |

**VII.b) One-dimensional test problems for the Relativistic Hydrodynamics Equations**

In this Sub-Section, we present seven one-dimensional test cases for the Relativistic Hydrodynamic equations. The system has been described in Sub-section 6.2. The initial states, value of the specific heat constant $(\gamma)$ and final times $(t_{end})$ for these seven test problems are given in Table V. All the test problems from Table V have computational domain $[-0.5, 0.5]$. Test 1 to Test 4 in Table V are described in [44],[47]. Test 5 to Test 7 are given in [23].

**Table V gives the left and right initial states, specific heat constant ($\gamma$) and final times ($t_{end}$) of the three test problems for the Relativistic Hydrodynamics Equations.**

|  | $x$ | $\rho$ | $v_x$ | $v_y$ | $p$ | $\gamma$ | $t_{end}$ |
|---|---|---|---|---|---|---|---|
| Test 1 | $x < 0$ | 1 | -0.6 | 0 | 10 | 5/3 | 0.4 |
|  | $x > 0$ | 10 | 0.5 | 0 | 20 |  |  |
| Test 2 | $x < 0$ | 10.0 | 0 | 0 | 40/3 | 5/3 | 0.4 |
|  | $x > 0$ | 1 | 0 | 0 | $10^{-6}$ |  |  |
| Test 3 | $x < 0$ | 1 | 0 | 0 | $10^3$ | 5/3 | 0.4 |
|  | $x > 0$ | 1 | 0 | 0 | $10^{-2}$ |  |  |
| Test 4 | $x < 0$ | 1 | 0.9 | 0 | 1 | 4/3 | 0.4 |
|  | $x > 0$ | 1 | 0 | 0 | 10 |  |  |
| Test 5 | $x < 0$ | 1 | -0.7 | 0 | 20 | 5/3 | 0.4 |
|  | $x > 0$ | 1 | 0.7 | 0 | 20 |  |  |
| Test 6 | $x < -0.4$ | 1.0 | 0.0 | 0.0 | 1000.0 | 1.4 | 0.43 |
|  | $-0.4 < x < 0.4$ | 1.0 | 0.0 | 0.0 | 0.01 |  |  |
|  | $x > 0.4$ | 1.0 | 0.0 | 0.0 | 100.0 |  |  |
| Test 7 | $x < 0$ | 5 | 0 | 0 | 50 | 5/3 | 0.35 |
|  | $x > 0$ | $2 + 0.3\sin(50x)$ | 0 | 0 | 5 |  |  |

The first five (Test 1-5) problems in Table V are the Riemann problems for which the exact solution is computed using the exact Riemann solver given in [44]. The exact solution is denoted by solid lines. Fig. 6 shows the density, velocity and pressure profiles for the Test 1 obtained from the fifth order scheme, and we observe that the results are relatively close to the exact solution and resolve both rarefaction and the contact for the density variable. The seventh and ninth order schemes also perform well on this problem and are not shown here. The exact solution for the Test 2 consist of a left-moving rarefaction wave, a contact discontinuity, and a right-moving shock wave in the density variable. The test problem is run using the fifth order scheme. Fig. 7 shows that the obtained results have efficiently captured the shock, the rarefaction wave, and the contact discontinuity in the density variable with a good agreement with the exact solution. The seventh and ninth order schemes also perform well on this problem and are not shown here. The exact solution for the Test 3 contains a narrow shock in the density variable which is difficult to capture. In Fig. 8, we obtain results for Test 3 from the seventh order scheme and observe that the scheme is able to capture the narrow shock. However, the resolution of the simulation can be increased to obtain even closer resemblance with the exact solution. The fifth and ninth order schemes also perform well on this problem and are not shown here. In Fig. 9, we show results for the Test 4 obtained from the seventh order scheme and observe that the obtained profile for the density variable has captured the slowly left-moving shock wave, the contact discontinuity, and a fast right-moving shock wave, and the solution is in good agreement with the exact solution. To avoid small oscillations in the obtained results (oscillations as seen in [23]) for this test case, we have utilized the flattening algorithm described in [9]. The fifth and ninth order schemes also perform well on this problem and are not shown here. In Fig. 10, we show results for the Test 5 obtained from the ninth order scheme. The simulation for this test case was performed utilizing the flattening algorithm presented in [9]. We observe a similar undershoot in the density variable as in [23] (also see [42]). The obtained density result captured the left-moving rarefaction wave and the right-moving rarefaction wave with a good agreement with the exact solution. The fifth and seventh order schemes also perform well on this problem and are not shown here.

Next, we consider the relativistic blast wave problem which is denoted by Test 6 in Table V. Given the extreme nature of the test case and the need to accurately capture the waves, a grid with 4000 zones is used. Similar to the Blast wave problem in Euler flows, we employ the flattening algorithm presented in [9] to simulate this particular test case. The zoomed solutions within the interval $[0.0, 0.03]$, obtained using the ninth order scheme, are shown in Fig. 11. We also plot the reference solution using solid lines, which has been obtained using the third order scheme with 15000 zones. We see a clear resemblance between the obtained result and the reference solution. The fifth and seventh order schemes also perform well on this problem and are not shown here.

Next, we present results for the Test 7 of Table V. This problem is used to test the capabilities of shock-capturing schemes to accurately capture and resolve the small-scale flow features. The exact solution contains sinusoidal density profile with both shocks and rarefaction waves. In Fig. 12, we plot the density, velocity and pressure profiles obtained from the seventh order scheme with 400 zones. We also plot the reference solution in solid black lines obtained using the third order scheme with 4000 zones. We see that the result in Fig. 12 closely matches the reference solution. The fifth and ninth order schemes also perform well on this problem and are not shown here. We re-run the simulation at 140 zones using the third, fifth, seventh and ninth order schemes and plot the density profiles in Fig. 13. We observe that, at such a low resolution, the third order result is much more diffusive in capturing the sinusoidal profile comparatively to the fifth, seventh and ninth order results. This highlights the value of higher order schemes.

### VII.c) One-dimensional test problems for the Ten-Moment rarefied gas flow model

In this Sub-Section, we present three one-dimensional Riemann problems for the Ten-Moment rarefied gas flow. The system has been described in Sub-section 6.3. The initial states and final times $(t_{end})$ for these three test problems are specified in Table VI. A description of the considered Riemann problems is given in [52] (also see [15]).

**Table VI gives the left and right initial states and final times $(t_{end})$ of the three test problems for the Ten-Moment Rarefied Gas Equations.**

|  | $x$ | $\rho$ | $v_x$ | $v_y$ | $p_{xx}$ | $p_{xy}$ | $p_{yy}$ | $t_{end}$ |
|---|---|---|---|---|---|---|---|---|
| Test 1 | $x < 0$ | 1 | 0 | 0 | 2 | 0.05 | 0.6 | 0.125 |
|  | $x > 0$ | 0.125 | 0 | 0 | 0.2 | 0.1 | 0.2 |  |
| Test 2 | $x < 0$ | 1 | 1 | 1 | 1 | 0 | 1 | 0.125 |
|  | $x > 0$ | 1 | -1 | -1 | 1 | 0 | 1 |  |
| Test 3 | $x < 0$ | 2 | -0.5 | -0.5 | 1.5 | 0.5 | 1.5 | 0.15 |
|  | $x > 0$ | 1 | 1 | 1 | 1 | 0 | 1 |  |

Test 1 is the Sod's shock-tube Riemann problem for Ten-Moment, whose solution contains a shock, contact wave and rarefaction wave. The resulting density, x-velocity, y-velocity and pressure tensor profiles obtained from the fifth order scheme are shown in Fig. 14. The exact solution is shown in solid lines. We observe a close match between the numerical result and the exact solution. The seventh and ninth order schemes also perform well on this problem and are not shown here. Next, we consider the Two-Shock wave problem. The initial states for this problem are given by Test 2 in Table VI. Fig. 15 shows the density, x-velocity, y-velocity and pressure tensor profiles obtained from the seventh order scheme. The exact solution is denoted by the solid lines. We see that the scheme is able to capture the two opposite moving shock waves. Similar to the Euler's case in [42], we also see an undershoot at the center of the domain for the density and an overshoot for the pressure component $p_{yy}$. The fifth and ninth order schemes also perform well on this problem and are not shown here. Test 3 in Table VI is known as Two-Rarefaction waves problem. The exact solution of this Riemann problem contains two rarefaction waves. The results for the density, x-velocity, y-velocity and pressure tensor variables obtained from the ninth order scheme are shown in Fig. 16. We observe a good match between the numerical and exact solution. The fifth and seventh order schemes also perform well on this problem and are not shown here.

## VIII) Multidimensional Test Problems

In this Section, we focus on several two-dimensional test problems. We wish to demonstrate that the presented scheme works well for multi-dimensional stringent problems. In Sub-section 8.1, we present multi-dimensional problems for the Euler flow; In Sub-section 8.1, we present multi-dimensional test problems for the Relativistic Hydrodynamics flow, and in Sub-section 8.3, we consider a multi-dimensional problem for the Ten-Moment model. We used a CFL of 0.4 with a third order SSP-RK time-stepping scheme for all the simulations presented in this Section.

## VIII.a) Multi-dimensional test problems for Euler Flow.

In this Sub-section, we present two multi-dimensional problems for Euler flow. The first test problem is the Forward facing step problem. This problem was first introduced by Woodward and Colella [60]. We simulate the problem on a computational domain that spans $[0,3] \times [0,1]$. At the left boundary of the domain, an ideal gas flows in at a speed of Mach 3 with a density of 1.4, a pressure of 1 and a ratio of specific heats of 1.4. At the upper corner, a forward-facing step is set up at the position $(0.6, 0.2)$. At the left and right boundary, we use outflow boundary conditions. At the top and bottom boundary, we apply the reflective boundary conditions, and that includes all parts of the forward-facing step. At the step corner, the singularity has been treated in the same manner as in Woodward and Colella [60]. The problem was run to a final time of 0.4 on a $1440 \times 480$ zone mesh. We use the flattening algorithm presented in [7] to simulate this particular test case. Fig. 17 shows the obtained result for the density variable at a final time of 0.4 using the fifth order accurate HLL-based AFD WENO scheme. The panel shows the sharp profiles and a pronounced roll-up of the vortex sheet can be seen. The seventh and ninth order schemes also perform well on this problem; therefore, they are not shown here.

Next, we consider the Double Mach Reflection (DMR) problem from Woodward and Colella [60] for the Euler flow. We use the same setup here. The problem simulates the multi-dimensional similarity solution that develops when an angled wedge is placed in a supersonic flow. Here we consider a domain that spans $[0,4] \times [0,1]$. A Mach 10 shock, positioned at an angle of $60^o$ to the bottom boundary, is initialized at the boundary point $x = 1/6$. For the value of $x < 1/6$, the post-shock conditions are used at the boundary, which matches the pre-wedge states. For values of $x > 1/6$, the reflective boundary condition is used, which exactly matches the windward face of the wedge. The upper boundary is modified consistently with the motion of the oblique shock. At the left edge, we use the inflow boundary conditions. At the right edge, we use the outflow boundary conditions. The un-shocked region is initialized with a density of 1.4, a pressure of 1, and a ratio of specific heats of 1.4. The problem was run to a final time of 0.2 on a $1920 \times 480$ zone mesh. We use the flattening algorithm presented in [7] to simulate this particular test case. Fig. 18 shows the plot for the density variable obtained from the seventh order accurate HLL-based AFD WENO scheme. It is customary to only image the partial domain $[0,3] \times [0,1]$. We also supply a small panel that shows a zoomed-up region around the Mach stem. We observe that the scheme is able to capture the instabilities (Fig. 18b) that develop around the Mach stem. This shows the significance of higher-order schemes. The fifth and ninth order schemes also perform well on this problem and therefore, they are not shown here.

**VIII.b) Multi-dimensional test problems for the Relativistic Hydrodynamics Equations**

In this Sub-section, we present multiple two-dimensional test problems for the relativistic flows. In Table VII, we consider three two-dimensional Riemann problems from [23]. All the test problems from Table VII have the computational domain $[-0.5, 0.5] \times [-0.5, 0.5]$. We simulate the first problem (2DRP-1) from Table VII using the fifth order accurate LLF-based AFD-WENO scheme on a mesh of $200 \times 200$ zones. For this problem, Fig. 19 shows the 25 equally spaced contours for the density and pressure logarithms. We observe that the four initial vortex sheets mutually interact with each other, resulting in the formation of a spiral pattern characterized by the reduced rest-mass density near the center of the computational domain. The obtained results are precise and closely resemble the results reported in reference [23]. The seventh and ninth order schemes also perform well on this problem and therefore, they are not shown here. Next, we run the second problem (2DRP-2) from Table VII using the seventh order accurate LLF-based AFD-WENO scheme on a mesh of $200 \times 200$ zones. For this problem, Fig. 20 shows the 25 equally spaced contours for the density and pressure logarithms. The solution is characterized by the interactions of planar rarefaction waves, leading to the formation of two symmetrical shock waves, and the results depicted in Fig. 20 demonstrate that the scheme successfully captures these intricate structures. The fifth and ninth order schemes also perform well on this problem and therefore, they are not shown here. Next, we simulate the last problem (2DRP-3) from Table VII using the ninth order accurate LLF-based AFD-WENO scheme on a mesh of $200 \times 200$ zones. For this problem, Fig. 21 shows the 25 equally spaced contours for the density and pressure logarithms. In this particular test case, the initial membranes experience breakdown, leading to the appearance of two contact discontinuities along the left and bottom boundaries of the domain, and at the top-right boundary, the system evolves, giving rise to the formation of two curved front shocks. The results in Fig. 21 demonstrate that the chosen scheme effectively captures these complex structures. The fifth and seventh order schemes also perform well on this problem and therefore, they are not shown here.

**Table VII gives the initial states, specific heat constant $(\gamma)$ and final times $(t_{end})$ of the three multi-dimensional Riemann problems for the Relativistic Hydrodynamics Equations.**

|  | $x, y$ | $\rho$ | $v_x$ | $v_y$ | $p$ | $\gamma$ | $t_{end}$ |
|---|---|---|---|---|---|---|---|
| **2DRP-1** | $x>0, y>0$ | 0.5 | 0.5 | -0.5 | 5 | 5/3 | 0.4 |
|  | $x<0, y>0$ | 1 | 0.5 | 0.5 | 5 |  |  |
|  | $x<0, y<0$ | 3 | -0.5 | 0.5 | 5 |  |  |
|  | $x>0, y<0$ | 1.5 | -0.5 | -0.5 | 5 |  |  |
| **2DRP-2** | $x>0, y>0$ | 1 | 0 | 0 | 1 | 5/3 | 0.4 |
|  | $x<0, y>0$ | 0.5771 | -0.3529 | 0 | 0.4 |  |  |
|  | $x<0, y<0$ | 1 | -0.3529 | -0.3529 | 1 |  |  |
|  | $x>0, y<0$ | 0.5771 | 0 | -0.3529 | 0.4 |  |  |
| **2DRP-3** | $x>0, y>0$ | 0.0351452161 | 0 | 0 | 0.1629310565 | 5/3 | 0.4 |
|  | $x<0, y>0$ | 0.1 | 0.7 | 0 | 1 |  |  |
|  | $x<0, y<0$ | 0.5 | 0 | 0 | 1 |  |  |
|  | $x>0, y<0$ | 0.1 | 0 | 0.7 | 1 |  |  |

Next, we consider the shock-bubble interaction problem from [29]. In [16], authors have simulated the same problem in a discontinuous Galerkin framework. In this test, a shock wave moving to the left interacts with a bubble and generates several interesting wave patterns. The initial states for the shock and bubbles are given in Table VIII. We consider two types of bubbles centered at $(215,0)$ with a radius of 25 units. The first bubble from table VIII is denoted by SB-1 and the second bubble is denoted by SB-2. The computational domain is given by $[0,325]\times[-45,45]$. At the top and bottom boundaries ($y = \pm 45$), we use the reflective boundary conditions. At the left and right boundaries ($x = 0$ and $x = 325$) we use the Dirichlet boundary conditions with the boundary values given by the respective initial left and right shock states specified in Table VIII.

**Table VIII gives the initial states of the shock and the two bubbles for the two-dimensional shock-bubble interaction problem for Relativistic Hydrodynamics Equations. The specific heat constant $(\gamma)$ and final times $(t_{end})$ for both bubbles are also specified in this table.**

|  | $x, y$ | $\rho$ | $v_x$ | $v_y$ | $p$ | $\gamma$ | $t_{end}$ |
|---|---|---|---|---|---|---|---|
| **Shock:** | $x < 265$ | 1 | 0 | 0 | 0.05 | | |
|  | $x > 265$ | 1.86522508063 | -0.19678110737 | 0 | 0.15 | | |
| **Bubble-1: (SB-1)** | $(x-215)^2 + y^2 < 25^2$ | 0.1358 | 0 | 0 | 0.05 | 5/3 | 450 |
| **Bubble-2: (SB-2)** | $(x-215)^2 + y^2 < 25^2$ | 3.1538 | 0 | 0 | 0.05 | 5/3 | 500 |

For the first shock-bubble problem (SB-1), we simulate the problem till the final time $t = 450$ and present the obtained results for the density variable in Fig. 22. Results shown in Fig. 22a are obtained using the third order scheme; Results shown in Fig. 22b are obtained using the fifth order scheme; Results shown in Fig. 22c are obtained using the seventh order scheme and Fig. 22d contains the results obtained using the ninth order AFD-WENO scheme. We observe that all the employed schemes exhibit good precision in capturing the small-scale wave structures. Notably, the higher-order schemes demonstrate a greater level of detail in representing these wave structures, hence showing the value of higher-order schemes.

The setup for the second shock-bubble problem (SB-2) is very similar to the previous one, except that the inside states of the bubble are now heavier. This means that there is more mass inside the bubble compared to bubble-1. The final time for this problem is given by $t = 500$. Fig. 23a, b, c, and d show the density profiles obtained using the third, fifth, seventh, and ninth order accurate AFD-WENO schemes, respectively. Once again, we observe that all the schemes effectively capture the discontinuities and wave structures with remarkable precision. As previously, higher-order schemes yield more accurate results, highlighting their superior performance.

**VIII.c) Multi-dimensional test problems for the Ten-Moment rarefied gas flow model**

In this Sub-section, we present a two-dimensional test problem from [45] (also see [15]) for the Ten-Moment rarefied gas flow model. The initial states for this test contain low-density and low-pressure zones. The computational domain for this test is $[-2, 2] \times [-2, 2]$. We use outflow boundary conditions at all the boundaries. The initial profiles are given by

$$\rho = 1, \qquad p_{xx} = p_{yy} = 2, \, p_{xy} = 0, \qquad (v_x, v_y) = 8(x/r, y/r), \quad \text{where } r = \sqrt{x^2 + y^2}.$$

We simulate the problem using the seventh-order accurate LLF-based AFD-WENO scheme on a mesh of $200 \times 200$ zones. The final time is given by $t = 0.05$. We use the flattening algorithm (with $\kappa = 1$) described in Appendix B to simulate this particular test case. Fig. 24 contains the obtained results. Profile in Fig. 24a shows the density variable, Fig. 24b corresponds to the $p_{xx}$-component of the pressure tensor, Fig. 24c corresponds to the $p_{xy}$-component of the pressure tensor, and Fig. 24d corresponds to the $p_{yy}$-component of the pressure tensor. We observe a good resemblance between the obtained result and reported results in [45][15]. The fifth and ninth order schemes also perform well on this problem and therefore, they are not shown here.

**IX) Conclusions**

In this paper we have striven to make the alternative finite difference WENO (AFD-WENO) as easily accessible and general-purpose as classical finite difference WENO (FD-WENO). A well-developed AFD-WENO offers many extremely desirable advantages over classical FD-WENO. First, it can work with any type of Riemann solver and different fields of study may have their own special Riemann solver that they find beneficial. Such a flexibility is not available in FD-WENO, whereas it is indeed available in AFD-WENO formulations. For instance, some Riemann solvers are better at positivity preservation, which can help with stringent problems.

Other Riemann solvers are better at preserving stationary contact discontinuities, which can help with well-balancing. Yet other Riemann solvers might be better at handling non-conservative products, which could be useful in applications where the governing equation simply cannot be written in conservation form. This flexibility of invoking the Riemann solver at pointwise locations can also be exploited to ensure that the final scheme respects the preservation of free stream conditions on curvilinear meshes. On such curvilinear meshes, the flux reconstruction of classical FD-WENO methods becomes a liability and this enables AFD-WENO methods come to the fore. For all the very good reasons mentioned above, the overarching motivation of this paper was to derive an AFD-WENO formulation that is absolutely general-purpose and easily accessible with all requisite formulae available in one place, i.e. within this paper. The AFD-WENO formulation that we have obtained has a computational complexity that is comparable to classical FD-WENO. This has the happy consequence that there are only up-sides, and no down-sides to switch from classical FD-WENO to the AFD-WENO presented here.

To achieve the broad goal of developing AFD-WENO into a general purpose algorithm for conservation laws, three major barriers had to be overcome. Overcoming these barriers, therefore, provides the three major goals of this paper, which are addressed in each of the next three paragraphs.

One of the reasons why AFD-WENO has not seen much uptake in the community stems from the fact that it is quite difficult to understand the scheme. Indeed, the update equation for AFD-WENO entails many higher derivatives of the fluxes at zone boundaries, and it is difficult to understand where these derivatives come from. In this work, Section II demystifies the derivation of the AFD-WENO update equation by working out the special case of third order AFD-WENO in great detail. In the associated Appendix A we also provide a script based on a computer algebra system that shows how the fifth order scheme can be obtained as a natural outcome of applying the script. Higher order AFD-WENO schemes can also be naturally derived by extending the script. Therefore, in Section II we achieve the first goal of this paper which is to make AFD-WENO very accessible to the greater community.

AFD-WENO also relies on interpolation rather than the more familiar WENO reconstruction. In reconstruction, the volume-averaged entities are reconstructed up to the desired accuracy via a suitable polynomial basis. The same polynomial basis are also used in AFD-WENO.

However, now we require that the interpolating polynomial should match the point values in the zones that make up the interpolating stencil. Section III provides all the explicit formulae for doing this up to ninth order. In that sense, the second goal of this paper is to serve as a one-stop-shop for all the interpolation formulae that are needed in AFD-WENO up to ninth order.

The structure of the AFD-WENO update is such that it relies on the evaluation of the higher derivatives of the fluxes at zone boundaries. When the solution is smooth, these derivatives contribute to higher order accuracy. When the solution is non-smooth, these higher derivatives can be a source of spurious oscillations. In the past, this was one of the major stumbling blocks in the development of AFD-WENO schemes. Some solutions that are specific to the PDE at hand had been developed for controlling these higher order derivatives. In this paper a general-purpose method is developed which relies on a novel WENO interpolation that takes the derivatives of the fluxes at the zone centers as its inputs and returns non-linearly hybridized higher order derivatives of the flux terms at zone boundaries as its output. In the spirit of serving as a one-stop-shop for all the interpolation formulae that are needed in AFD-WENO, our third goal is to document all those interpolation formulae up to eighth order. This is done in Section IV and it is sufficient to yield an AFD-WENO scheme that has up to ninth order of accuracy. Since our method is based on the alternative formulation of WENO finite difference as in [32,33], it can be adapted to maintain free-stream exactly on curvilinear meshes as in [33], which is left for our future work.

A pointwise description of the ADF-WENO algorithm that simplifies its implementation is presented in Section V. Accuracy analysis is presented in Section VI. One dimensional tests are presented in Section VII. Multidimensional tests are presented in Section VIII. In future work we will extend our methods to include hyperbolic PDEs which are completely general in that they may have flux form or non-conservative products.

**Acknowledgements**

DSB acknowledges support via NSF grant, NSF-AST-2009776, NASA-2020-1241 and NASA grant 80NSSC22K0628. DSB and HK acknowledge support from a Vajra award, VJR/2018/00129 and also a travel grant from Notre Dame International. CWS acknowledges support via NSF grant DMS-2309249.

**Ethical Statement**
**i. Compliance with Ethical Standards** : This manuscript complies with all ethical standards for scientific publishing.
**ii. (in case of Funding) Funding** : The funding has been acknowledged. DSB acknowledges support via NSF grants NSF-19-04774, NSF-AST-2009776, NASA-2020-1241 and NASA-80NSSC22K0628. DSB and HK acknowledge support from a Vajra award, VJR/2018/00129. CWS acknowledges support via AFOSR grant FA9550-20-1-0055 and NSF grant DMS-2010107.
**iii. Conflict of Interest** : On behalf of all authors, the corresponding author states that there is no conflict of interest.
**iv. Ethical approval** : N/A
**v. Informed consent** : N/A

**vi. Data Statement** : All data that was used in the generation of the figures has been stored and available for later use.

# Appendix A) Mathematica Script for the derivation of a 5th Order AFD-WENO Scheme

We present the Mathematica script for understanding the formulae in Shu and Osher (1988). The script is extensible to all orders.

```
(* ---------------------------------------------------------------------- *)
(* This is the Mathematica equivalent of the Merriman (2003) derivation.  *)
(* First step is to retrieve the classical derivation. *)
(* But we also examine the non-conservative products. *)
(* ---------------------------------------------------------------------- *)
(* We start with a Taylor series expansion of the flux. *)

    f[x_]:= f0 + x * dxf + x^2 * dx2f / 2 + x^3 * dx3f / 6 + x^4 * dx4f / 24 +
        x^5 * dx5f / 120

(* ------------------------------------------------------------ *)
(* Next, evaluate the flux and its higher derivatives at x = 1/2. Dxn denotes the nth derivative. *)

    Fph = f[x]/.x->1/2
    DxFph = D[f[x],x]/.x->1/2
    Dx2Fph = D[D[f[x],x],x]/.x->1/2
    Dx3Fph = D[D[D[f[x],x],x],x]/.x->1/2
    Dx4Fph = D[D[D[D[f[x],x],x],x],x]/.x->1/2
    Dx5Fph = D[D[D[D[D[f[x],x],x],x],x],x]/.x->1/2

(* ------------------------------------------------------------ *)
(* Similarly, we evaluate the flux and its higher derivatives at x = -1/2 *)

    Fmh = f[x]/.x->-1/2
    DxFmh = D[f[x],x]/.x->-1/2
    Dx2Fmh = D[D[f[x],x],x]/.x->-1/2
    Dx3Fmh = D[D[D[f[x],x],x],x]/.x->-1/2
    Dx4Fmh = D[D[D[D[f[x],x],x],x],x]/.x->-1/2
    Dx5Fmh = D[D[D[D[D[f[x],x],x],x],x],x]/.x->-1/2

(* ---------------------------------------------------------------------- *)
(* Consider the 5th order case for the derivation of AFD-WENO *)
(* ---------------------------------------------------------------------- *)

    Flux5Ph = Fph + a * DxFph + b * Dx2Fph + c * Dx3Fph + d * Dx4Fph + e * Dx5Fph

    Flux5Mh = Fmh + a * DxFmh + b * Dx2Fmh + c * Dx3Fmh + d * Dx4Fmh + e * Dx5Fmh

(* At x = 0 we only want the "dxf" coefficient to be unity and all higher *)
(* order terms to be zero. *)
```

```
Coefficient[ Expand[Simplify[Flux5Ph - Flux5Mh]], dxf]   (* This must evaluate to 1 *)

Solve [ {Coefficient[ Expand[Simplify[Flux5Ph - Flux5Mh]], dx2f] == 0,
     Coefficient[ Expand[Simplify[Flux5Ph - Flux5Mh]], dx3f] == 0,
     Coefficient[ Expand[Simplify[Flux5Ph - Flux5Mh]], dx4f] == 0,
     Coefficient[ Expand[Simplify[Flux5Ph - Flux5Mh]], dx5f] == 0,
     Coefficient[ Expand[Simplify[Flux5Ph - Flux5Mh]], dx6f] == 0},
    {a, b, c, d, e} ]

a = 0
b = - (1/24)
c = 0
d = (7/5760)
e = 0

(* The "derivatives of the flux" that we should put at the zone boundaries are given by *)

Flux5Ph = Fph - (1/24) * Dx2Fph + (7/5760) * Dx4Fph
Flux5Mh = Fmh - (1/24) * Dx2Fmh + (7/5760) * Dx4Fmh

(* The above coefficients give us the higher order derivatives of the fluxes that we need *)
(* at the zone boundaries in order to get the accurate flux gradient at the zone center. *)
(* -------------------------------------------------- *)
```

## Appendix B) Flattener function for the Ten-Moment rarefied gas flow model

The flattening algorithm is employed in simulations to identify regions of strong shocks within the computational domain. The algorithm aims to improve the simulation accuracy by reducing numerical artifacts near discontinuities. The flattener functions rely on comparing the divergence of the velocity field to a characteristic speed associated with the specific problem being solved. For the Euler flows, such functions have been previously defined by Colella and Woodward [21] and Balsara [7]. For the Relativistic Magneto-hydrodynamics equations, a flattener has been presented in Balsara and Kim [9]. Along the same lines, we give a flattener function for the Ten-Moment model. The method begins by calculating the divergence of the velocity, $(\nabla \cdot \mathbf{v})_{i,j}$, and sound-like speed $c_{s;i,j}$ (some approximation of the characteristic speed), within a specific zone $(i, j)$. In the two-dimensional cartesian mesh, the quantities are defined by

$$(\nabla \cdot \mathbf{v})_{i,j} = \frac{v_{x;i+1,j} - v_{x;i-1,j}}{\Delta x} + \frac{v_{y;i,j+1} - v_{y;i,j-1}}{\Delta y},$$

$$c_{s;i,j} = \sqrt{\frac{P_{i,j}}{\rho_{i,j}}}, \quad \text{where} \quad P_{i,j} = \sqrt{p_{xx;i,j} p_{yy;i,j} - p_{xy;i,j}^2}$$

To detect a shock, the undivided divergence of the velocity within a zone must be compared to the minimum sound-like speed in the zone $(i, j)$ and its immediate neighbors. The minimum sound-like speed from all the neighbors is obtained as follows

$$c_{s;i,j}^{\text{min-nbr}} = \min(c_{s;i-1,j-1}, c_{s;i-1,j}, c_{s;i-1,j+1}, c_{s;i,j-1}, c_{s;i,j}, c_{s;i,j+1}, c_{s;i+1,j-1}, c_{s;i+1,j}, c_{s;i+1,j+1})$$

In each zone with an extent of $\Delta x$ and $\Delta y$, the flattener function is defined as

$$\eta_{i,j} = \min\left[1, \max\left[0, \frac{|(\nabla \cdot \mathbf{v})_{i,j}| \max(\Delta x, \Delta y)}{\kappa c_{s;i,j}^{\text{min-nbr}}} - 1\right]\right]$$

The parameter $\kappa$ is set to 0.3, which has been found to work well across different orders and problem types. The flattener function does not modify the reconstruction when the flow is smooth or consists of rarefactions, and in that case, $\eta_{i,j} = 0$. However, the flattener function gradually increases from $\eta_{i,j} = 0$ to $\eta_{i,j} = 1$ when strong shocks are present.

The inclusion of pressure variation in the flattener algorithm allows for a more comprehensive stabilization of the flow simulation. It ensures that not only the zones already influenced by shocks but also the zones on edge receive appropriate flattening treatment. This improvement helps maintain numerical stability and accuracy throughout the simulation, particularly in regions where shocks are forming or propagating. In the x-direction, the flattener can be extended to the neighboring cell if the following conditions are satisfied.

$$\text{if } \left((\eta_{i,j} > 0) \text{ and } (\eta_{i+1,j} = 0) \text{ and } (P_{i,j} > P_{i+1,j})\right) \text{ then } \eta_{i+1,j} = \eta_{i,j}$$
$$\text{if } \left((\eta_{i,j} > 0) \text{ and } (\eta_{i-1,j} = 0) \text{ and } (P_{i,j} > P_{i-1,j})\right) \text{ then } \eta_{i-1,j} = \eta_{i,j}.$$

In situations involving multidimensional problems, the above strategy can be applied to each of the principal directions of the mesh.

**Figure Captions**

*Fig. 1 shows part of the mesh around zone "i". The mesh functions are collocated at the zone centers, as shown by the thick dots. The zone boundaries are shown by the vertical lines. The figure also shows the stencils associated with the zone "i" for the third and fifth order **pointwise** WENO-AO interpolation strategies. We have three smaller third order stencils and a large fifth order stencil. For third order WENO-AO, only the three smaller stencils are used, whereas the larger stencil is also used for fifth order WENO-AO. The interpolated variables at the zone boundaries are shown with a caret. The variables with a superscript star are resolved states obtained by the pointwise application of a Riemann solver at the zone boundaries.*

*Fig. 2 shows part of the mesh around zone boundary "i+1/2". The products of characteristic matrices with the gradients are evaluated pointwise at the zone centers, as shown by the thick dots. The zone boundaries are shown by the vertical lines. The figure also shows the stencils associated with the zone boundary "i+1/2" for the third and fifth order AFD-WENO schemes. We have two smaller third order stencils and a large fourth order stencil. The large stencil is, therefore, capable of providing the first and third derivatives when the smoothness in the solution warrants it. For third order AFD-WENO, only the two smaller stencils are used, whereas the larger stencil is also used for fifth order AFD-WENO. The first and third derivatives of the product of the characteristic matrix with the gradient are shown at the zone boundary of interest. Only the first derivatives are needed at third order, but third derivatives are also needed at fifth order.*

*Fig. 3) Euler Flow: Sod shock-tube problem. Panels a, b and c show the density, velocity and pressure profiles respectively at time t=0.2 obtained using the $5^{th}$ order LLF-based AFD-WENO scheme with 200 zones. Solid lines denote the exact solution. The $7^{th}$ and $9^{th}$ order schemes also perform well on this problem and are not shown here.*

*Fig. 4) Euler Flow: Lax shock-tube problem. Panels a, b and c show the density, velocity and pressure profiles respectively at time t=0.13 obtained using the $7^{th}$ order LLF-based AFD-WENO scheme with 200 zones. Solid lines denote the exact solution. The $5^{th}$ and $9^{th}$ order schemes also perform well on this problem and are not shown here.*

*Fig. 5) Euler Flow: Blast wave interaction problem. Panels a, b and c show the density, velocity and pressure profiles respectively at time t=0.038 obtained using the LLF-based $9^{th}$ order AFD-WENO scheme with 1000 zones. Solid lines denote the reference solution. The $5^{th}$ and $7^{th}$ order schemes also perform well on this problem and are not shown here.*

*Fig. 6) Relativistic Hydrodynamic Flow: Test-1 (Riemann Problem). Panels a, b and c show the density, velocity and pressure profiles respectively at time t=0.4 obtained using the $5^{th}$ order LLF-based AFD-WENO scheme with 200 zones. Solid lines denote the exact solution. The $7^{th}$ and $9^{th}$ order schemes also perform well on this problem and are not shown here.*

*Fig. 7) Relativistic Hydrodynamic Flow: Test-2 (Riemann Problem). Panels a, b and c show the density, velocity and pressure profiles respectively at time t=0.4 obtained using the $5^{th}$ order LLF-based AFD-WENO scheme with 200 zones. Solid lines denote the exact solution. The $7^{th}$ and $9^{th}$ order schemes also perform well on this problem and are not shown here.*

*Fig. 8) Relativistic Hydrodynamic Flow: Test-3 (Riemann Problem). Panels a, b and c show the density, velocity and pressure profiles respectively at time t=0.4 obtained using the $7^{th}$ order LLF-based AFD-WENO scheme with 400 zones. Solid lines denote the exact solution. The $5^{th}$ and $9^{th}$ order schemes also perform well on this problem and are not shown here.*

*Fig. 9) Relativistic Hydrodynamic Flow: Test-4 (Riemann Problem). Panels a, b and c show the density, velocity and pressure profiles respectively at time t=0.4 obtained using the $7^{th}$ order LLF-based AFD-WENO scheme with 200 zones. Solid lines denote the exact solution. The $5^{th}$ and $9^{th}$ order schemes also perform well on this problem and are not shown here.*

*Fig. 10) Relativistic Hydrodynamic Flow: Test-5 (Riemann Problem). Panels a, b and c show the density, velocity and pressure profiles respectively at time t=0.4 obtained using the $9^{th}$ order LLF-based AFD-WENO scheme with 200 zones. Solid lines denote the exact solution. The $5^{th}$ and $7^{th}$ order schemes also perform well on this problem and are not shown here.*

*Fig. 11) Relativistic Hydrodynamic Flow: Test-6 (Blast wave interaction problem). Panels a, b and c show the magnified image of the density, velocity and pressure profiles respectively at time t=0.43 obtained using the $9^{th}$ order LLF-based AFD-WENO scheme with 4000 zones. Solid lines denote the reference solution. The $5^{th}$ and $7^{th}$ order schemes also perform well on this problem and are not shown here.*

*Fig. 12) Relativistic Hydrodynamic Flow: Test-7 (Density perturbation Problem). Panels a, b and c show the density, velocity and pressure profiles respectively at time t=0.35 obtained using the $7^{th}$ order LLF-based AFD-WENO scheme with 400 zones. The converged solution was obtained on a mesh with 4000 zones and a $3^{rd}$ order scheme. Solid lines denote the reference solution.*

*Fig. 13) Relativistic Hydrodynamic Flow: Test-7 (Density perturbation Problem). Panel a shows the density profile at time t=0.35 obtained using the $3^{rd}$, $5^{th}$, $7^{th}$ and $9^{th}$ order LLF-based AFD-WENO schemes with 140 zones. Panel b shows zoomed image of the Panel a. Solid black line denotes the reference solution.*

*Fig. 14) Ten-Moment Rarefied Gas flow: Test-1 (Riemann Problem). Panels a, b and c show the density, x-velocity and y-velocity respectively; and panels d, e and f show the xx, xy and yy-component of the pressure tensor, respectively, at time t=0.125 obtained using the $5^{th}$ order LLF-based AFD-WENO scheme with 200 zones. Solid lines denote the exact solution. The $7^{th}$ and $9^{th}$ order schemes also perform well on this problem and are not shown here.*

*Fig. 15) Ten-Moment Rarefied Gas flow: Test-2 (Riemann Problem). Panels a, b and c show the density, x-velocity and y-velocity respectively; and panels d, e and f show the xx, xy and yy-component of the pressure tensor, respectively, at time t=0.125 obtained using the $7^{th}$ order LLF-based AFD-WENO scheme with 200 zones. Solid lines denote the exact solution. The $5^{th}$ and $9^{th}$ order schemes also perform well on this problem and are not shown here.*

*Fig. 16) Ten-Moment Rarefied Gas flow: Test-3 (Riemann Problem). Panels a, b and c show the density, x-velocity and y-velocity respectively; and panels d, e and f show the xx, xy and yy-component of the pressure tensor, respectively, at time t=0.15 obtained using the $9^{th}$ order LLF-based AFD-WENO scheme with 200 zones. Solid lines denote the exact solution. The $5^{th}$ and $7^{th}$ order schemes also perform well on this problem and are not shown here.*

*Fig. 17) Euler Flow: Forward facing step problem. Panel shows the density contours using the $5^{th}$ order accurate HLLI-based AFD-WENO scheme with $1440 \times 480$ zones. 30 contours were fit between a range of 0.1 and 6.62. The $7^{th}$ and $9^{th}$ order schemes also perform well on this problem and are not shown here.*

*Fig. 18) Euler Flow: Double Mach reflection problem. Fig. 18a shows the density contours using the $7^{th}$ order accurate HLLI-based AFD-WENO scheme with $1920 \times 480$. zones. Fig. 18b shows*

*the detailed view of the density profile. 20 contours were fit between a range of 1.0 and 20.5. The $5^{th}$ and $9^{th}$ order schemes also perform well on this problem and are not shown here.*

*Fig. 19) Relativistic Hydrodynamic Flow: 2DRP-1 (2D Riemann Problem-1). Fig. 19a shows the density logarithm and Fig. 19b shows the pressure logarithm at time t=0.4 obtained using the $5^{th}$ order LLF-based AFD-WENO scheme with 200 ×200 zones. 25 contours were fit in between the range of minimum and maximum value. The $7^{th}$ and $9^{th}$ order schemes also perform well on this problem and are not shown here.*

*Fig. 20) Relativistic Hydrodynamic Flow: 2DRP-2 (2D Riemann Problem-2). Fig. 20a shows the density logarithm and Fig. 20b shows the pressure logarithm at time t=0.4 obtained using the $7^{th}$ order LLF-based AFD-WENO scheme with 200 ×200 zones. 25 contours were fit in between the range of minimum and maximum value. The $5^{th}$ and $9^{th}$ order schemes also perform well on this problem and are not shown here.*

*Fig. 21) Relativistic Hydrodynamic Flow: 2DRP-3 (2D Riemann Problem-3). Fig. 21a shows the density logarithm and Fig. 21b shows the pressure logarithm at time t=0.4 obtained using the $9^{th}$ order LLF-based AFD-WENO scheme with 200 ×200 zones. 25 contours were fit in between the range of minimum and maximum value. The $5^{th}$ and $7^{th}$ order schemes also perform well on this problem and are not shown here.*

*Fig. 22a) Relativistic Hydrodynamic Flow: SB-1 (Shock-Bubble interaction-1). Panel shows the density profile at time t=450 obtained using the $3^{rd}$ order LLF-based AFD-WENO scheme with 650 ×180 zones.*

*Fig. 22b) Relativistic Hydrodynamic Flow: SB-1 (Shock-Bubble interaction-1). Panel shows the density profile at time t=450 obtained using the $5^{th}$ order LLF-based AFD-WENO scheme with 650 ×180 zones.*

*Fig. 22c) Relativistic Hydrodynamic Flow: SB-1 (Shock-Bubble interaction-1). Panel shows the density profile at time t=450 obtained using the $7^{th}$ order LLF-based AFD-WENO scheme with 650 ×180 zones.*

*Fig. 22d) Relativistic Hydrodynamic Flow: SB-1 (Shock-Bubble interaction-1). Panel shows the density profile at time t=450 obtained using the 9$^{th}$ order LLF-based AFD-WENO scheme with 650 ×180 zones.*

*Fig. 23a) Relativistic Hydrodynamic Flow: SB-2 (Shock-Bubble interaction-2). Panel shows the density profile at time t=500 obtained using the 3$^{rd}$ order LLF-based AFD-WENO scheme with 650 ×180 zones.*

*Fig. 23b) Relativistic Hydrodynamic Flow: SB-2 (Shock-Bubble interaction-2). Panel shows the density profile at time t=500 obtained using the 5$^{th}$ order LLF-based AFD-WENO scheme with 650 ×180 zones.*

*Fig. 23c) Relativistic Hydrodynamic Flow: SB-2 (Shock-Bubble interaction-2). Panel shows the density profile at time t=500 obtained using the 7$^{th}$ order LLF-based AFD-WENO scheme with 650 ×180 zones.*

*Fig. 23d) Relativistic Hydrodynamic Flow: SB-2 (Shock-Bubble interaction-2). Panel shows the density profile at time t=500 obtained using the 9$^{th}$ order LLF-based AFD-WENO scheme with 650 ×180 zones.*

*Fig. 24) Ten-Moment Rarefied Gas flow: 2D near vacuum test problem. Fig. 24a shows the density, Fig. 24b shows the $p_{xx}$ component, Fig. 24c Shows the $p_{xy}$ component and Fig. 24d shows the $p_{yy}$ component at time t=0.05 obtained using the 7$^{th}$ order LLF-based AFD-WENO scheme with 200 zones. The 5$^{th}$ and 9$^{th}$ order schemes also perform well on this problem and are not shown here.*

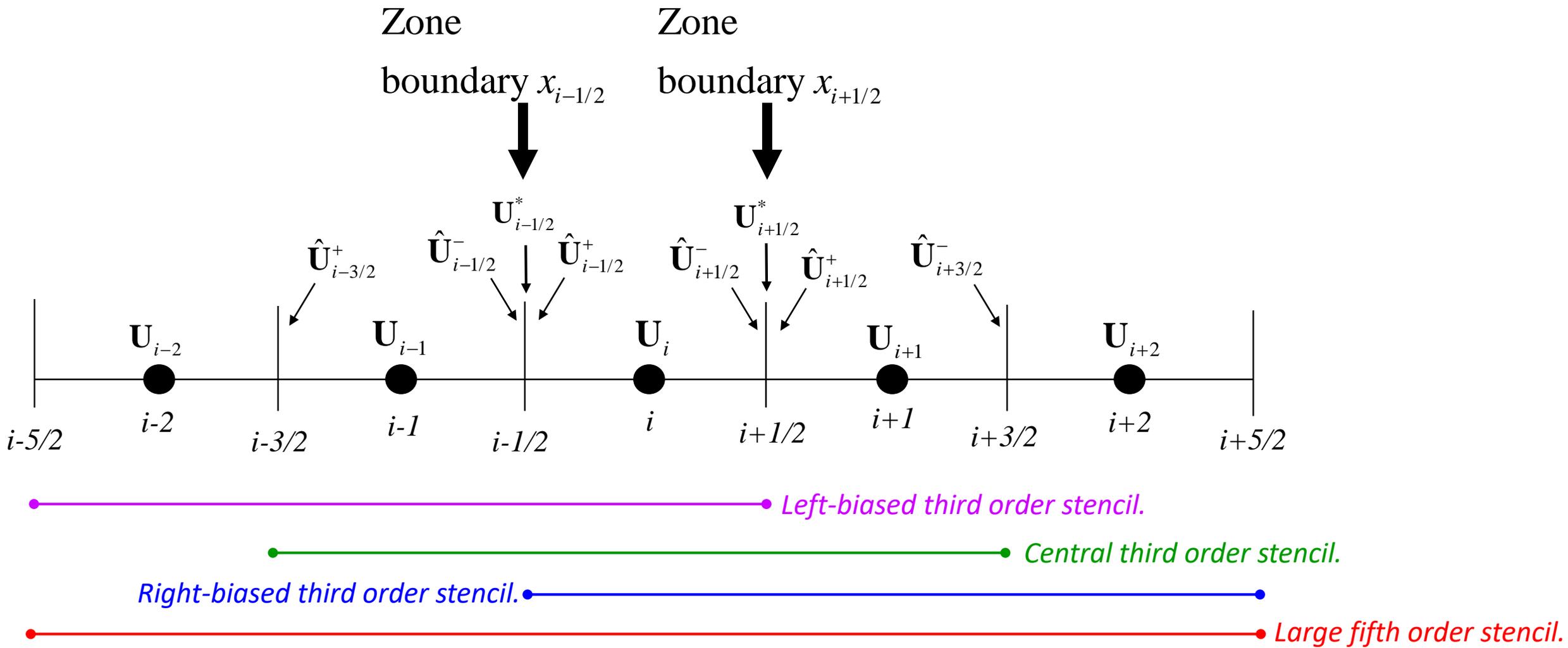

Fig. 1 shows part of the mesh around zone "i". The mesh functions are collocated at the zone centers, as shown by the thick dots. The zone boundaries are shown by the vertical lines. The figure also shows the stencils associated with the zone "i" for the third and fifth order **pointwise** WENO-AO interpolation strategies. We have three smaller third order stencils and a large fifth order stencil. For third order WENO-AO, only the three smaller stencils are used, whereas the larger stencil is also used for fifth order WENO-AO. The interpolated variables at the zone boundaries are shown with a caret. The variables with a superscript star are resolved states obtained by the pointwise application of a Riemann solver at the zone boundaries.

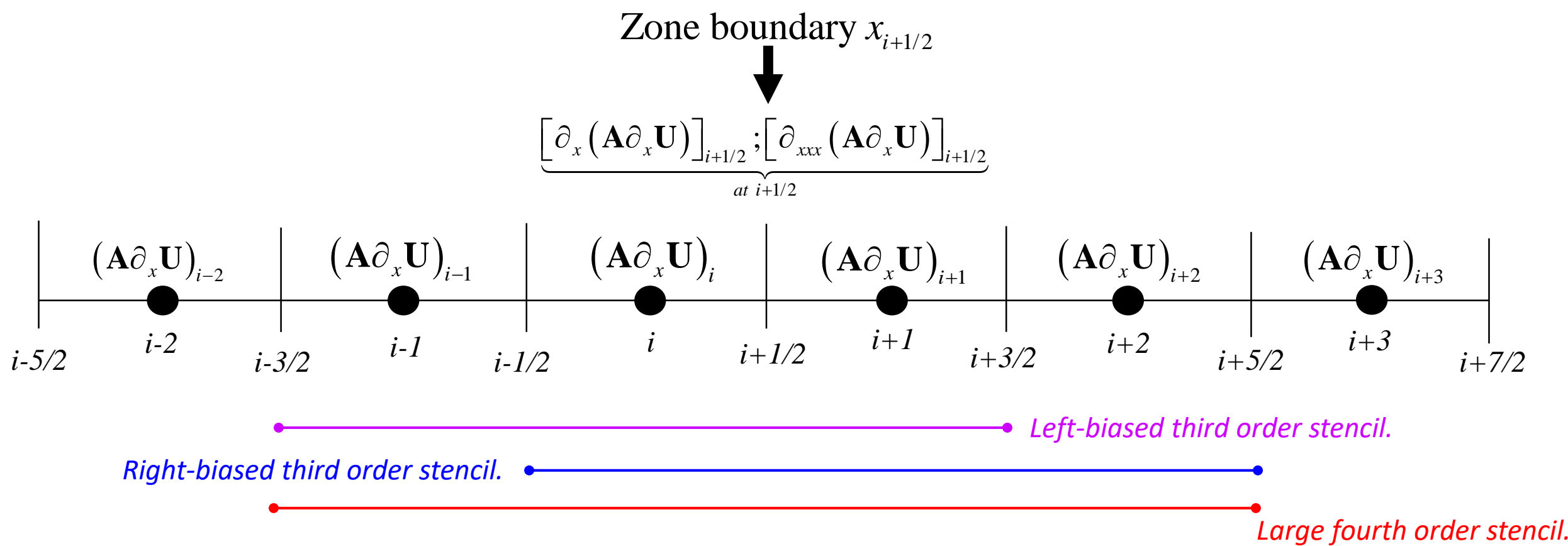

Fig. 2 shows part of the mesh around zone boundary "i+1/2". The products of characteristic matrices with the gradients are evaluated pointwise at the zone centers, as shown by the thick dots. The zone boundaries are shown by the vertical lines. The figure also shows the stencils associated with the zone boundary "i+1/2" for the third and fifth order AFD-WENO schemes. We have two smaller third order stencils and a large fourth order stencil. The large stencil is, therefore, capable of providing the first and third derivatives when the smoothness in the solution warrants it. For third order AFD-WENO, only the two smaller stencils are used, whereas the larger stencil is also used for fifth order AFD-WENO. The first and third derivatives of the product of the characteristic matrix with the gradient are shown at the zone boundary of interest. Only the first derivatives are needed at third order, but third derivatives are also needed at fifth order.

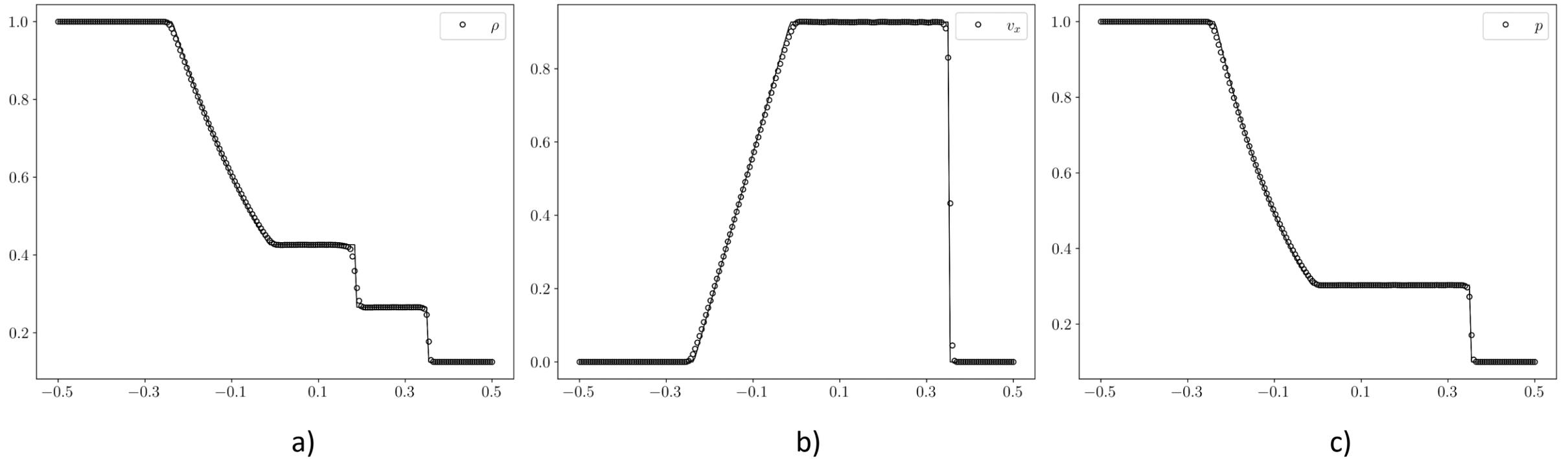

Fig. 3) Euler Flow: Sod shock-tube problem. Panels a, b and c show the density, velocity and pressure profiles respectively at time t=0.2 obtained using the 5$^{th}$ order LLF-based AFD-WENO scheme with 200 zones. Solid lines denote the exact solution. The 7$^{th}$ and 9$^{th}$ order schemes also perform well on this problem and are not shown here.

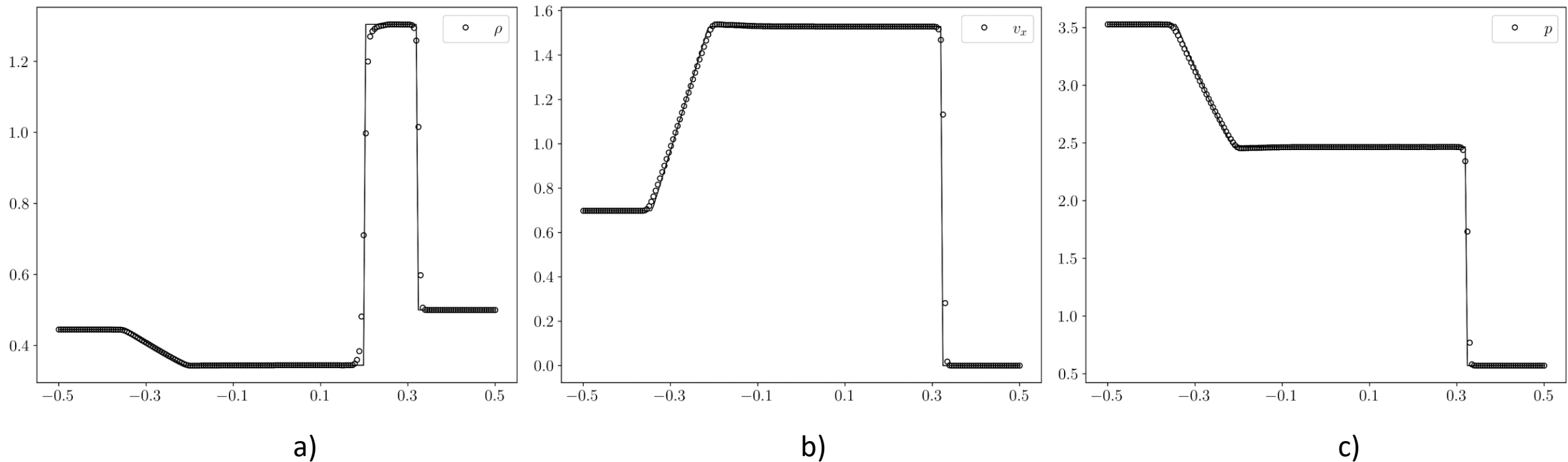

Fig. 4) Euler Flow: Lax shock-tube problem. Panels a, b and c show the density, velocity and pressure profiles respectively at time t=0.13 obtained using the 7th order LLF-based AFD-WENO scheme with 200 zones. Solid lines denote the exact solution. The 5th and 9th order schemes also perform well on this problem and are not shown here.

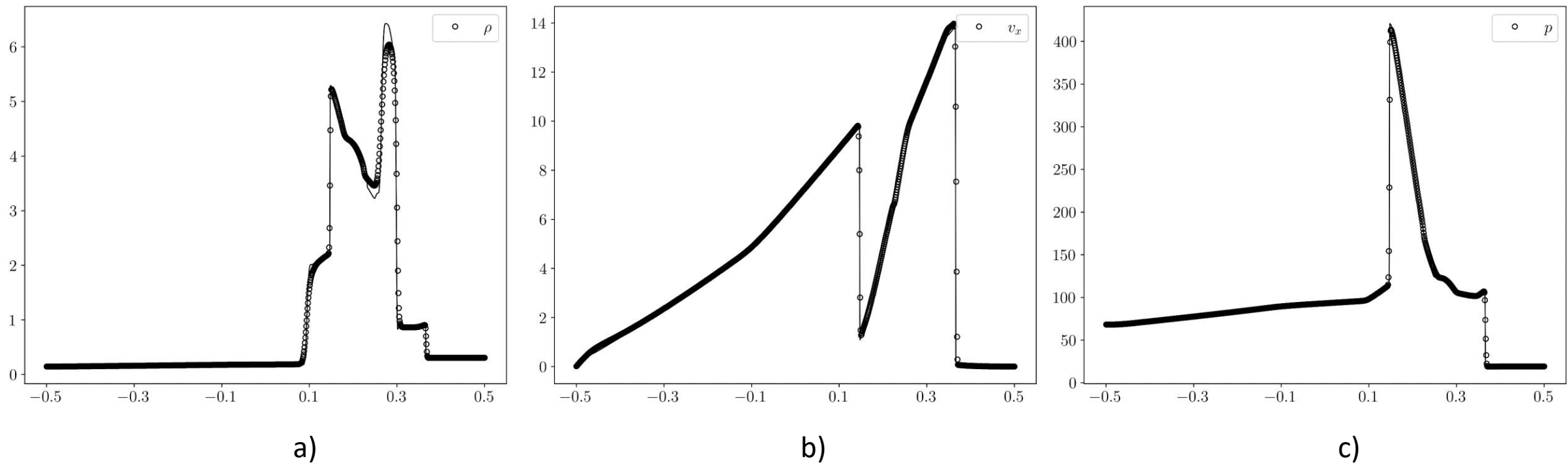

*Fig. 5) Euler Flow: Blast wave interaction problem. Panels a, b and c show the density, velocity and pressure profiles respectively at time t=0.038 obtained using the LLF-based 9$^{th}$ order AFD-WENO scheme with 1000 zones. Solid lines denote the reference solution. The 5$^{th}$ and 7$^{th}$ order schemes also perform well on this problem and are not shown here.*

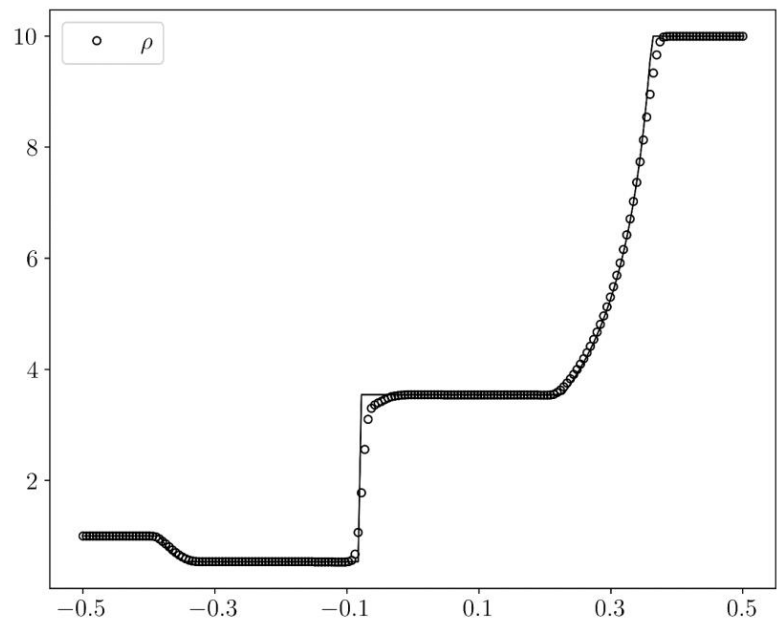 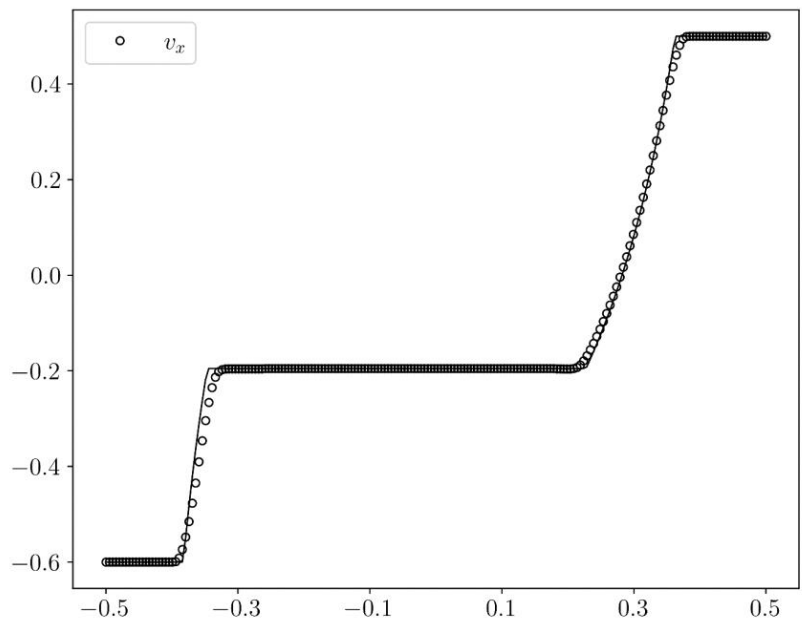 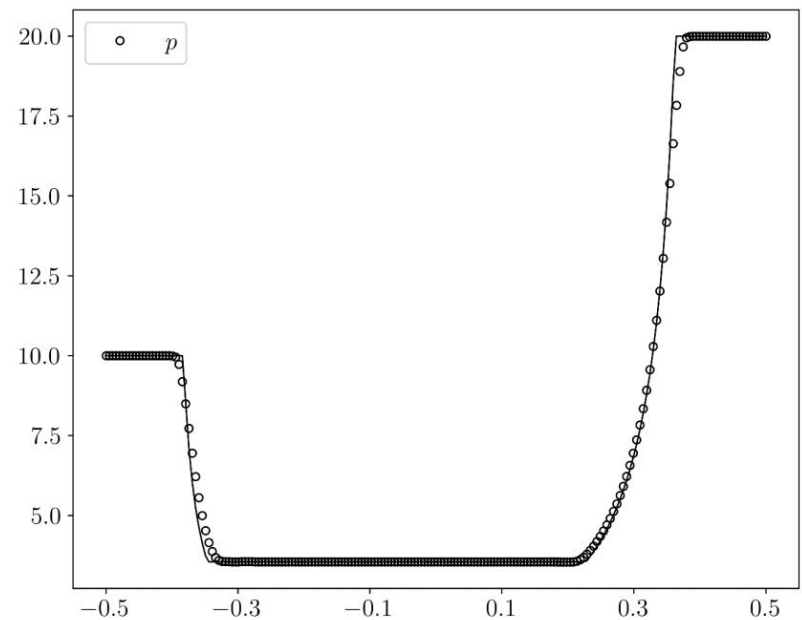

a)              b)              c)

*Fig. 6) Relativistic Hydrodynamic Flow: Test-1 (Riemann Problem). Panels a, b and c show the density, velocity and pressure profiles respectively at time t=0.4 obtained using the $5^{th}$ order LLF-based AFD-WENO scheme with 200 zones. Solid lines denote the exact solution. The $7^{th}$ and $9^{th}$ order schemes also perform well on this problem and are not shown here.*

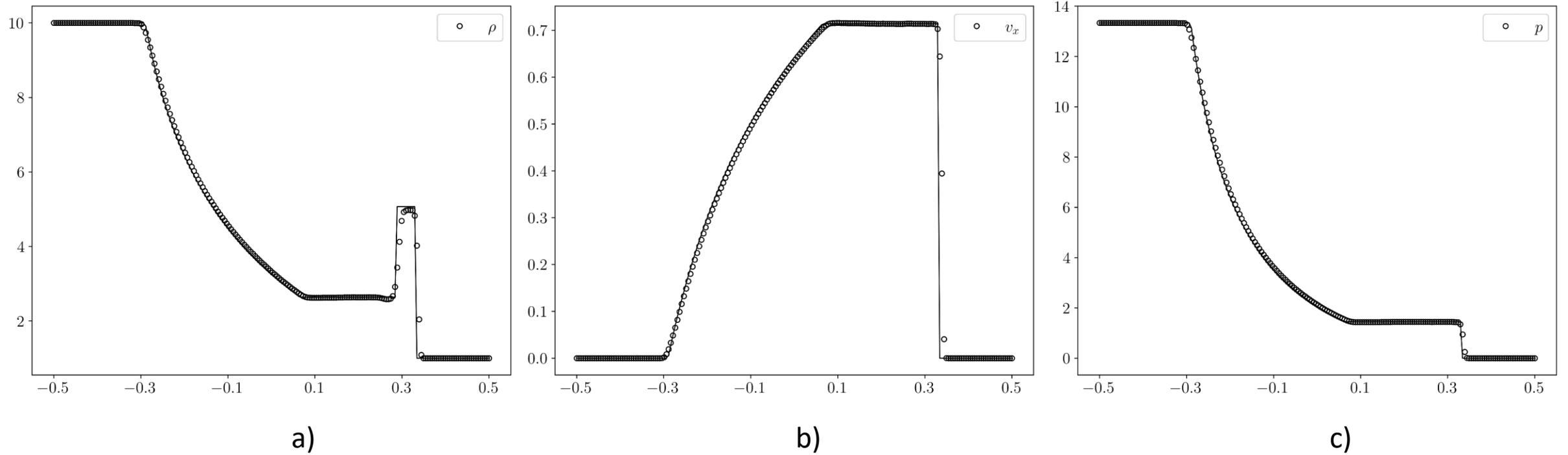

*Fig. 7) Relativistic Hydrodynamic Flow: Test-2 (Riemann Problem). Panels a, b and c show the density, velocity and pressure profiles respectively at time t=0.4 obtained using the 5th order LLF-based AFD-WENO scheme with 200 zones. Solid lines denote the exact solution. The 7th and 9th order schemes also perform well on this problem and are not shown here.*

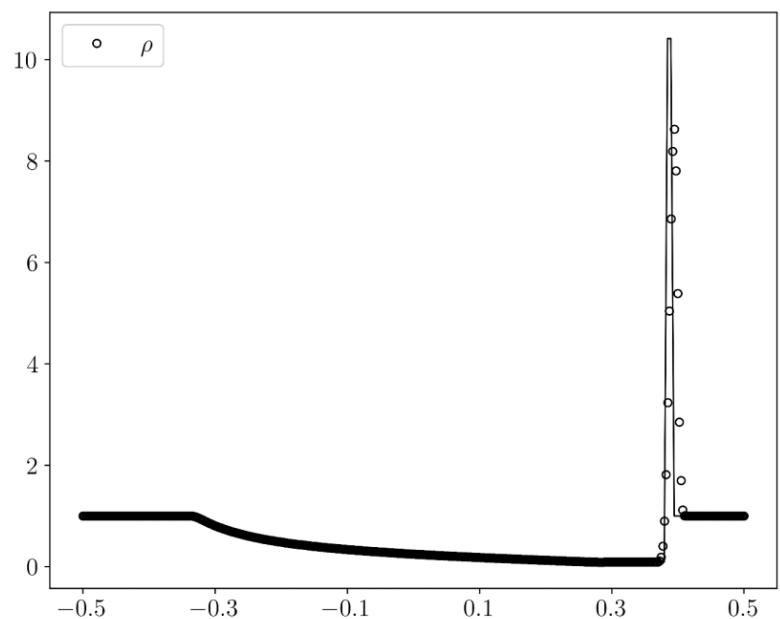 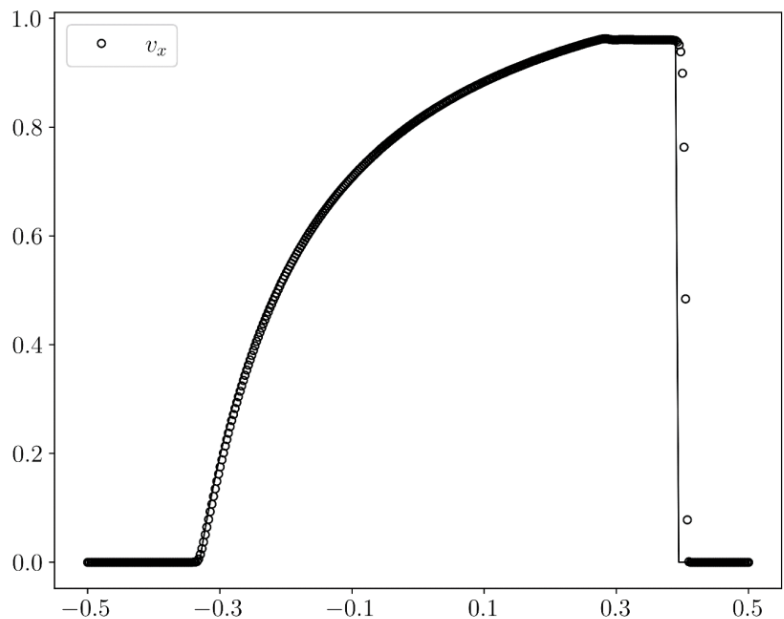 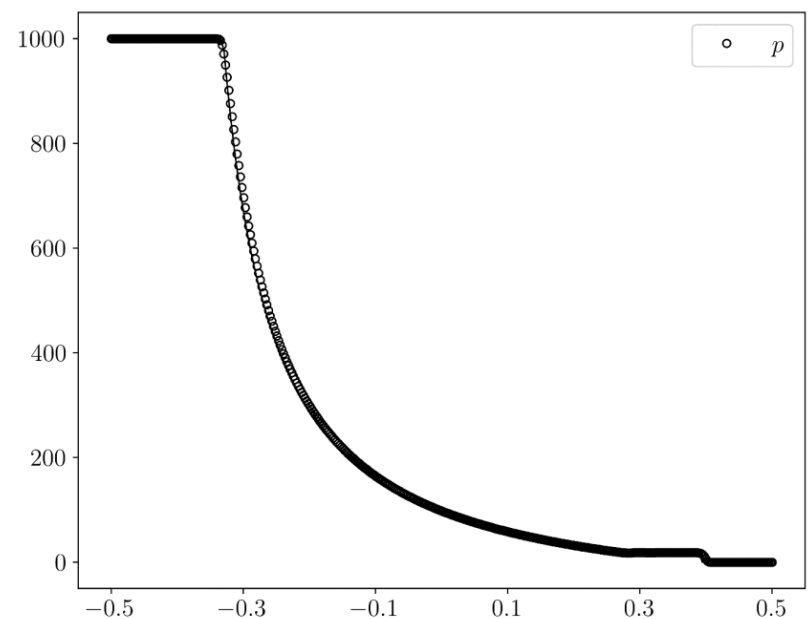

*Fig. 8) Relativistic Hydrodynamic Flow: Test-3 (Riemann Problem). Panels a, b and c show the density, velocity and pressure profiles respectively at time t=0.4 obtained using the 7$^{th}$ order LLF-based AFD-WENO scheme with 400 zones. Solid lines denote the exact solution. The 5$^{th}$ and 9$^{th}$ order schemes also perform well on this problem and are not shown here.*

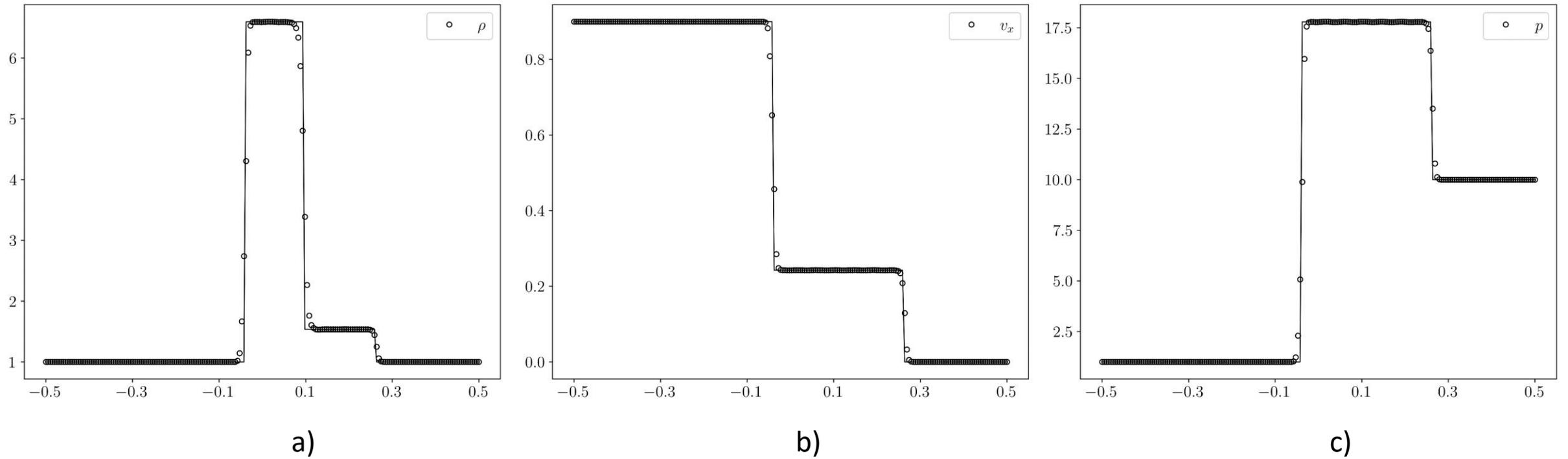

Fig. 9) Relativistic Hydrodynamic Flow: Test-4 (Riemann Problem). Panels a, b and c show the density, velocity and pressure profiles respectively at time t=0.4 obtained using the 7$^{th}$ order LLF-based AFD-WENO scheme with 200 zones. Solid lines denote the exact solution. The 5$^{th}$ and 9$^{th}$ order schemes also perform well on this problem and are not shown here.

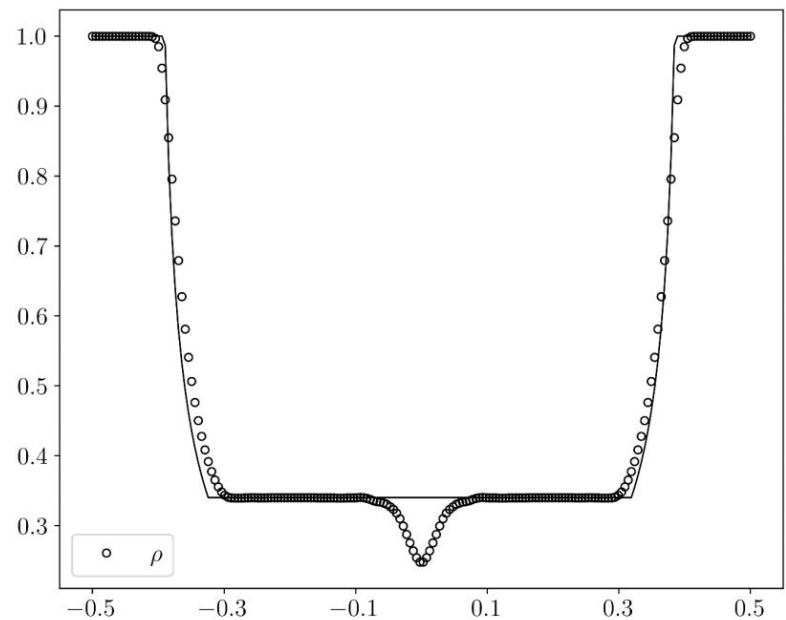 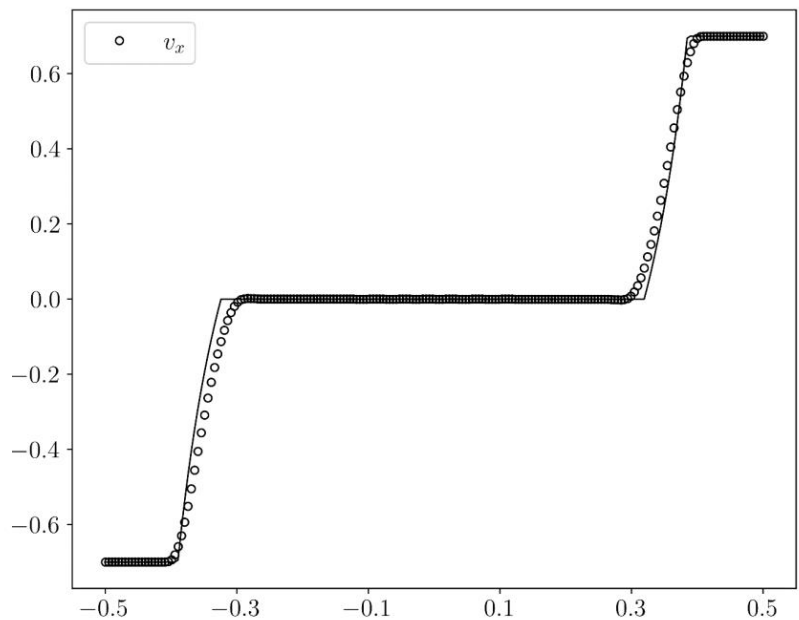 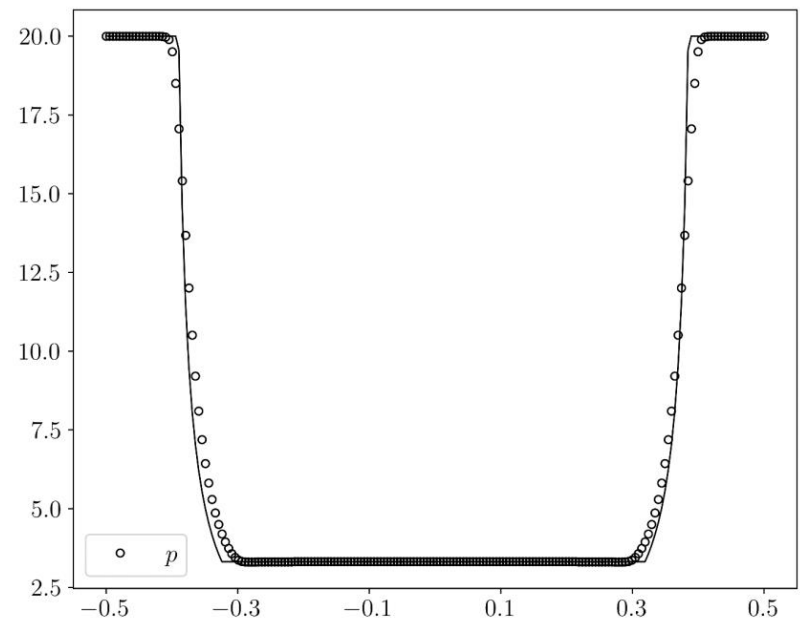

*Fig. 10) Relativistic Hydrodynamic Flow: Test-5 (Riemann Problem). Panels a, b and c show the density, velocity and pressure profiles respectively at time t=0.4 obtained using the 9$^{th}$ order LLF-based AFD-WENO scheme with 200 zones. Solid lines denote the exact solution. The 5$^{th}$ and 7$^{th}$ order schemes also perform well on this problem and are not shown here.*

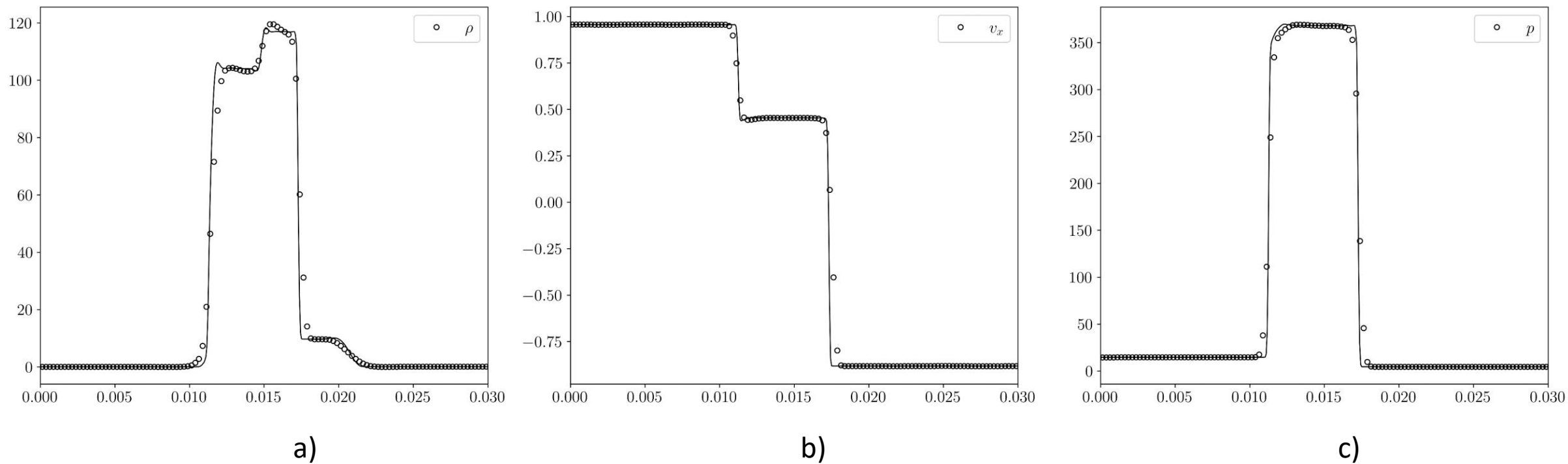

*Fig. 11) Relativistic Hydrodynamic Flow: Test-6 (Blast wave interaction problem). Panels a, b and c show the magnified image of the density, velocity and pressure profiles respectively at time t=0.43 obtained using the 9th order LLF-based AFD-WENO scheme with 4000 zones. Solid lines denote the reference solution. The 5th and 7th order schemes also perform well on this problem and are not shown here.*

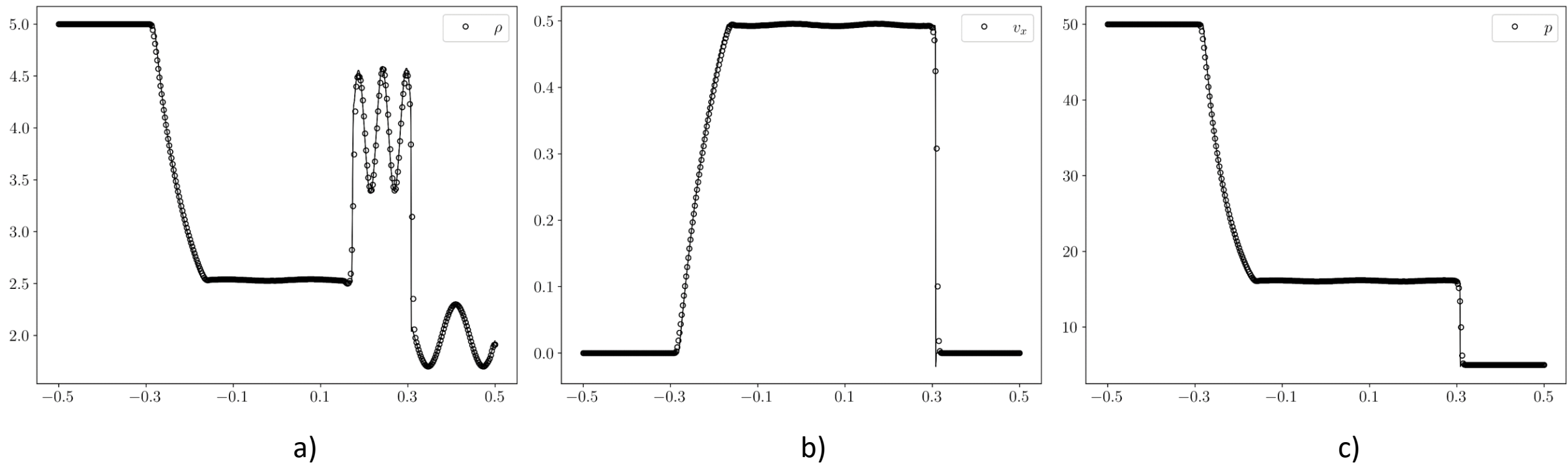

*Fig. 12) Relativistic Hydrodynamic Flow: Test-7 (Density perturbation Problem). Panels a, b and c show the density, velocity and pressure profiles respectively at time t=0.35 obtained using the 7th order LLF-based AFD-WENO scheme with 400 zones. The converged solution was obtained on a mesh with 4000 zones and a 3rd order scheme. Solid lines denote the reference solution.*

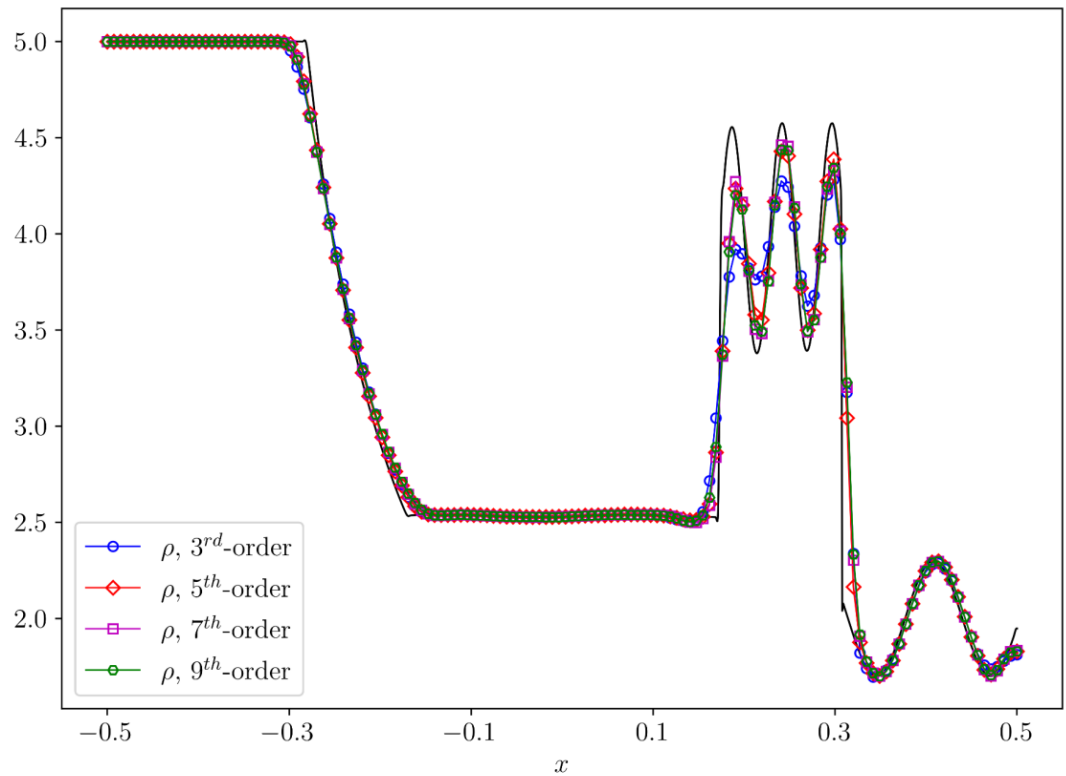 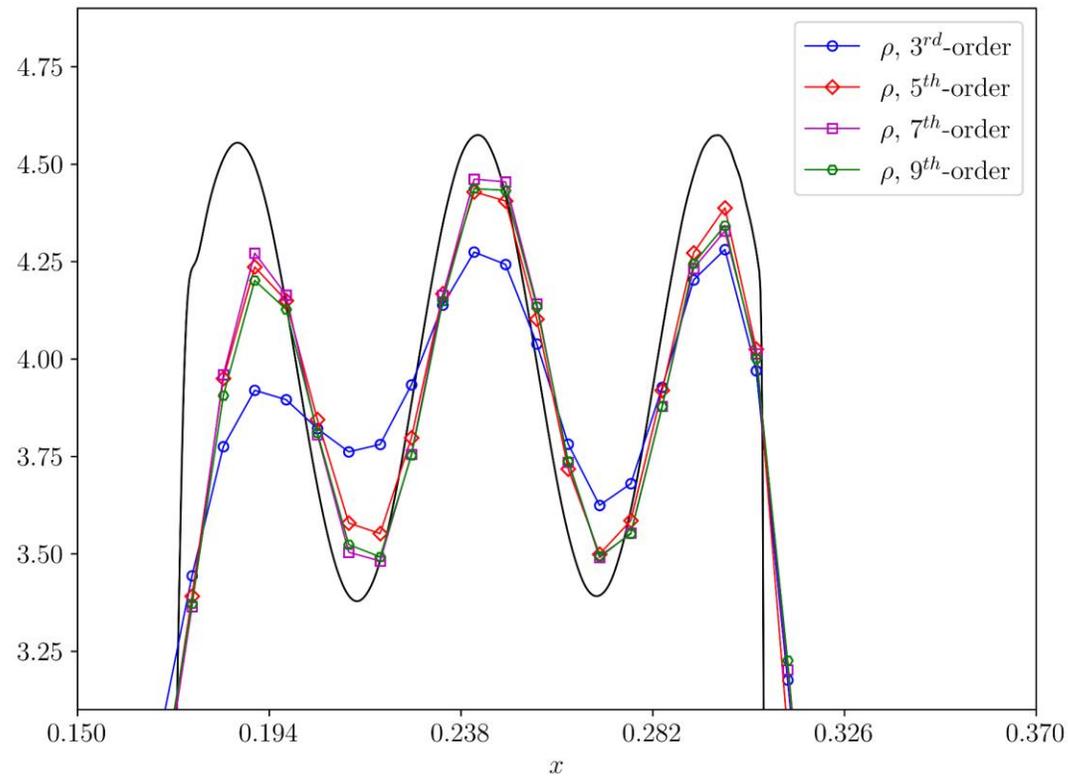

a)

b)

*Fig. 13) Relativistic Hydrodynamic Flow: Test-7 (Density perturbation Problem). Panel a shows the density profile at time t=0.35 obtained using the 3$^{rd}$, 5$^{th}$, 7$^{th}$ and 9$^{th}$ order LLF-based AFD-WENO schemes with 140 zones. Panel b shows zoomed image of the Panel a. Solid black line denotes the reference solution.*

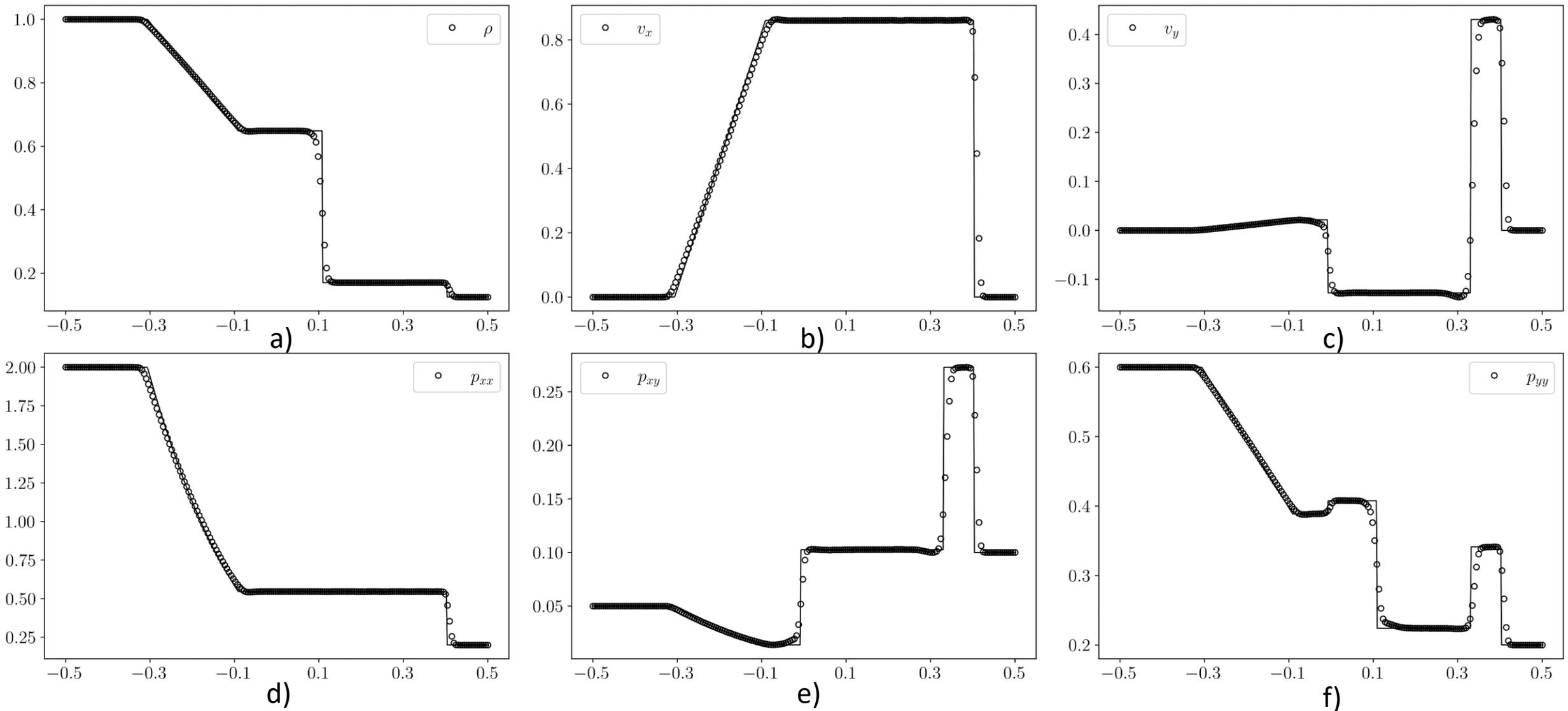

Fig. 14) Ten-Moment Rarefied Gas flow: Test-1 (Riemann Problem). Panels a, b and c show the density, x-velocity and y-velocity respectively; and panels d, e and f show the xx, xy and yy-component of the pressure tensor, respectively, at time t=0.125 obtained using the 5$^{th}$ order LLF-based AFD-WENO scheme with 200 zones. Solid lines denote the exact solution. The 7$^{th}$ and 9$^{th}$ order schemes also perform well on this problem and are not shown here.

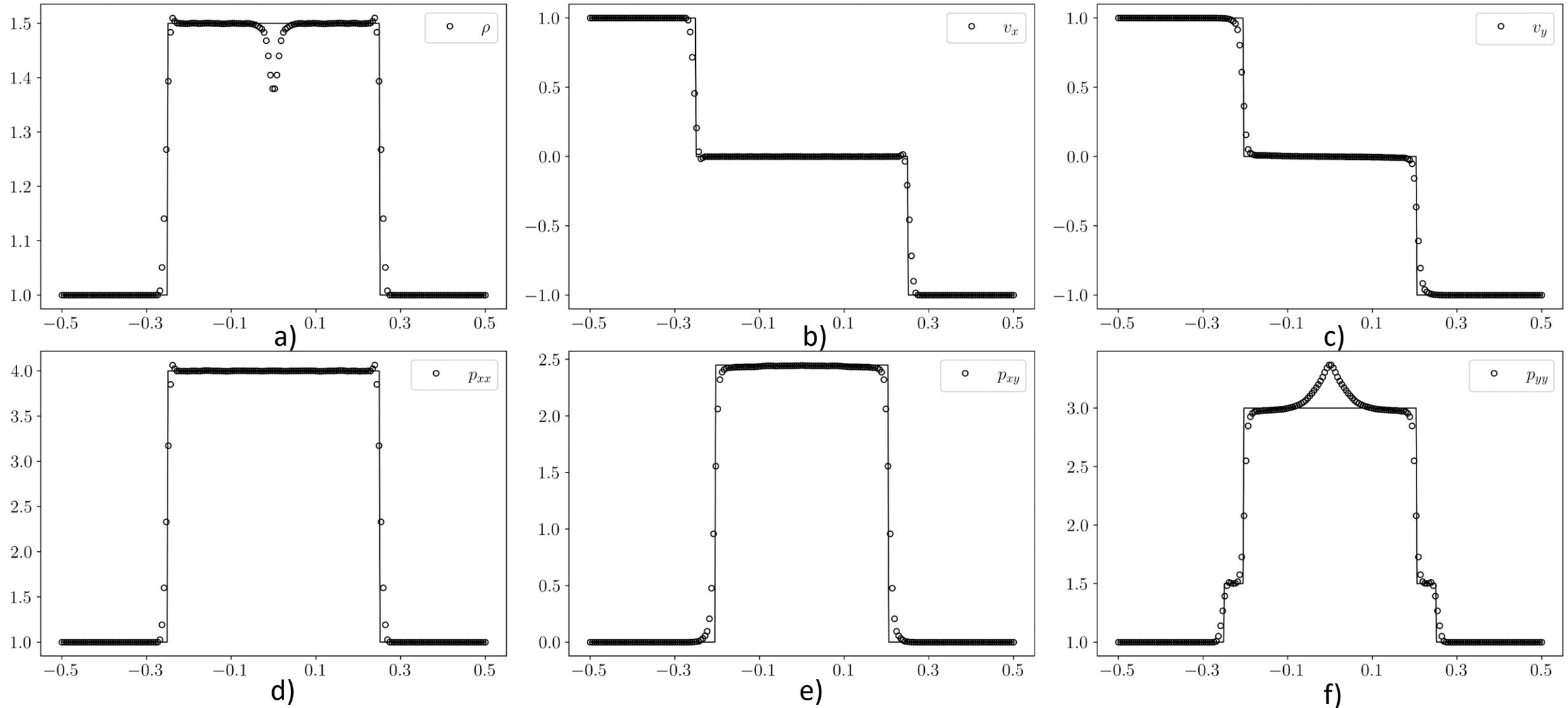

Fig. 15) Ten-Moment Rarefied Gas flow: Test-2 (Riemann Problem). Panels a, b and c show the density, x-velocity and y-velocity respectively; and panels d, e and f show the xx, xy and yy-component of the pressure tensor, respectively, at time t=0.125 obtained using the 7th order LLF-based AFD-WENO scheme with 200 zones. Solid lines denote the exact solution. The 5th and 9th order schemes also perform well on this problem and are not shown here.

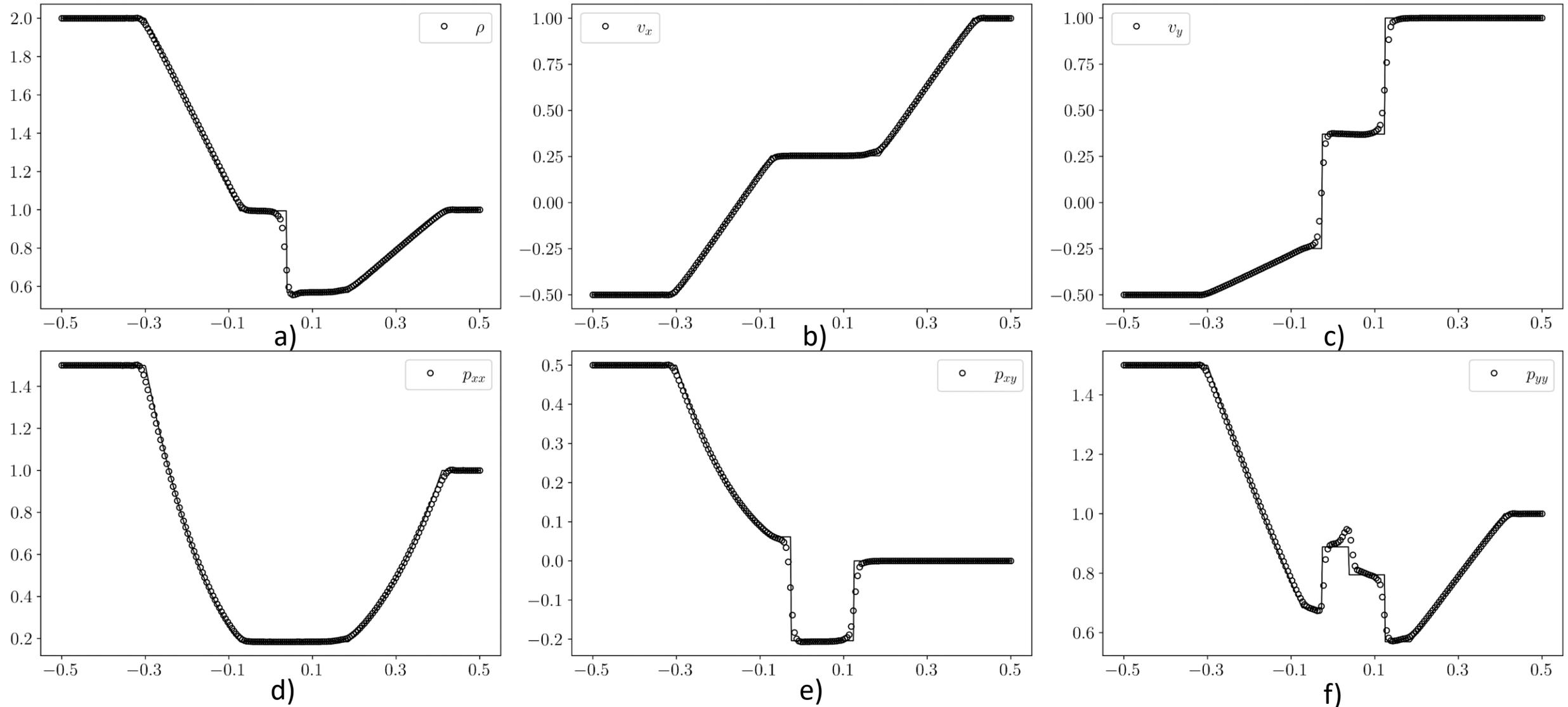

Fig. 16) Ten-Moment Rarefied Gas flow: Test-3 (Riemann Problem). Panels a, b and c show the density, x-velocity and y-velocity respectively; and panels d, e and f show the xx, xy and yy-component of the pressure tensor, respectively, at time t=0.15 obtained using the 9$^{th}$ order LLF-based AFD-WENO scheme with 200 zones. Solid lines denote the exact solution. The 5$^{th}$ and 7$^{th}$ order schemes also perform well on this problem and are not shown here.

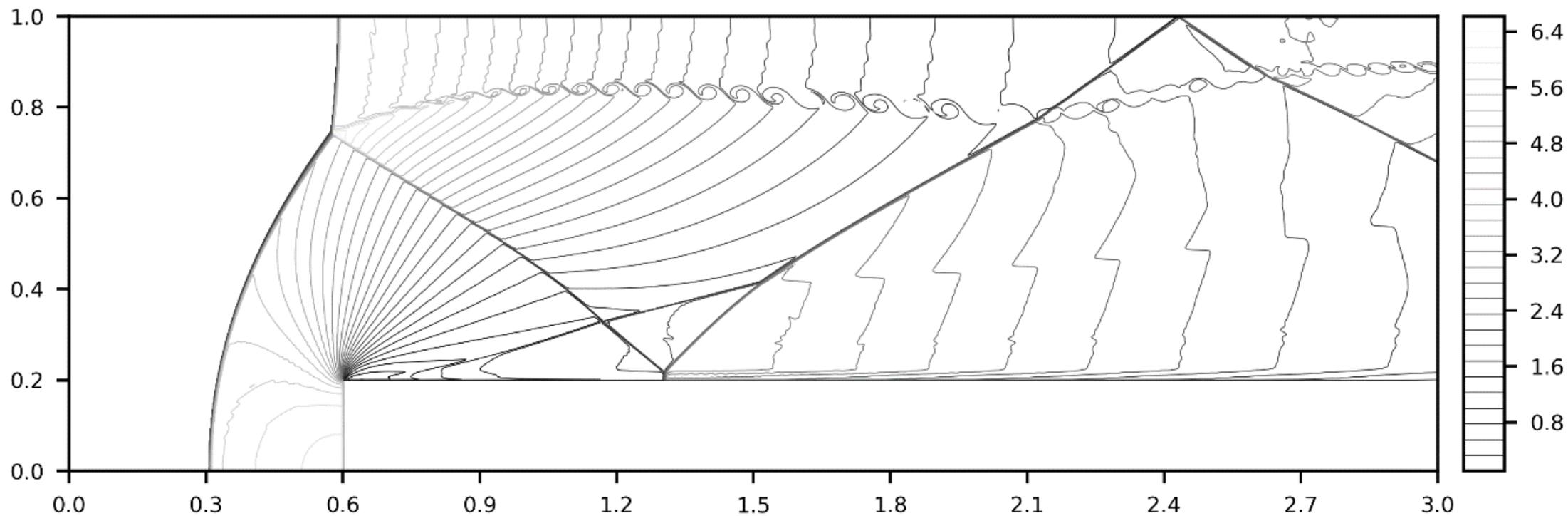

*Fig. 17) Euler Flow: Forward facing step problem. Panel shows the density contours using the 5th order accurate HLLI-based AFD-WENO scheme with 1440 ×480 zones. 30 contours were fit between a range of 0.1 and 6.62. The 7th and 9th order schemes also perform well on this problem and are not shown here.*

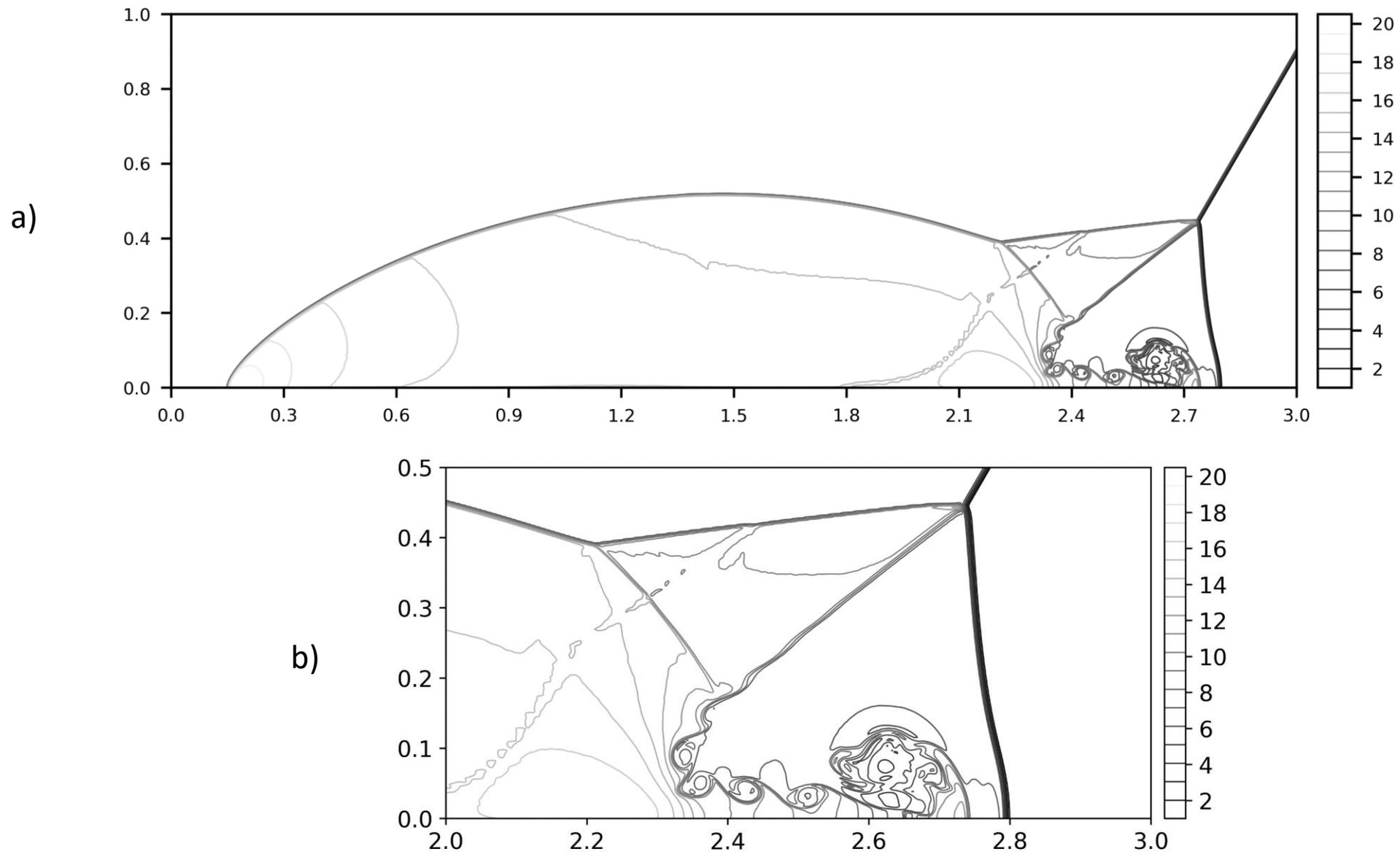

*Fig. 18) Euler Flow: Double Mach reflection problem. Fig. 18a shows the density contours using the 7th order accurate HLLI-based AFD-WENO scheme with 1920×480. zones. Fig. 18b shows the detailed view of the density profile. 20 contours were fit between a range of 1.0 and 20.5. The 5th and 9th order schemes also perform well on this problem and are not shown here.*

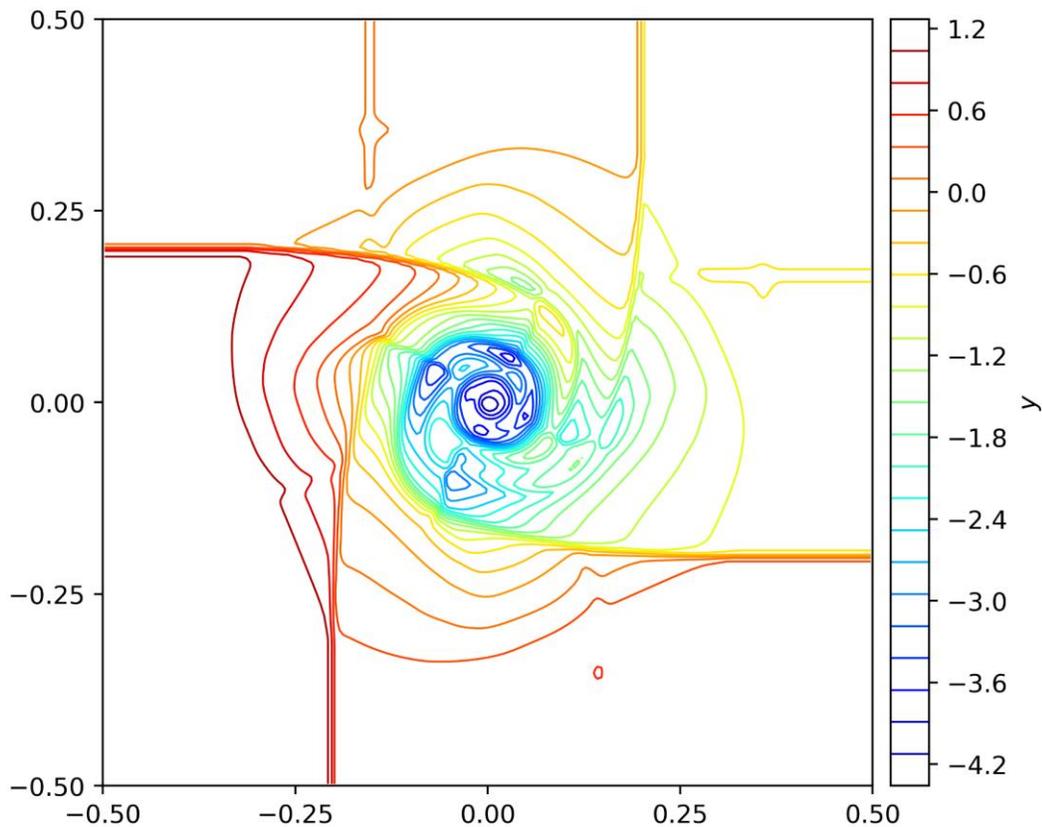 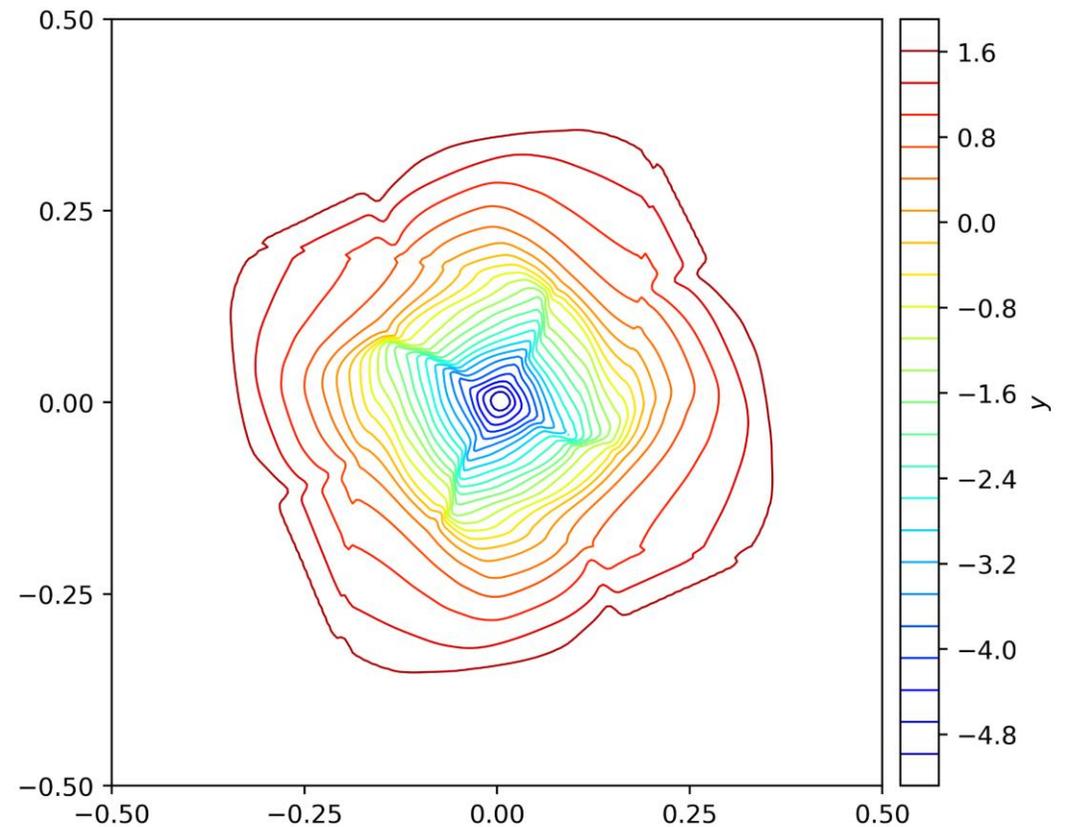

*Fig. 19) Relativistic Hydrodynamic Flow: 2DRP-1 (2D Riemann Problem-1). Fig. 19a shows the density logarithm and Fig. 19b shows the pressure logarithm at time t=0.4 obtained using the 5$^{th}$ order LLF-based AFD-WENO scheme with 200 ×200 zones. 25 contours were fit in between the range of minimum and maximum value. The 7$^{th}$ and 9$^{th}$ order schemes also perform well on this problem and are not shown here.*

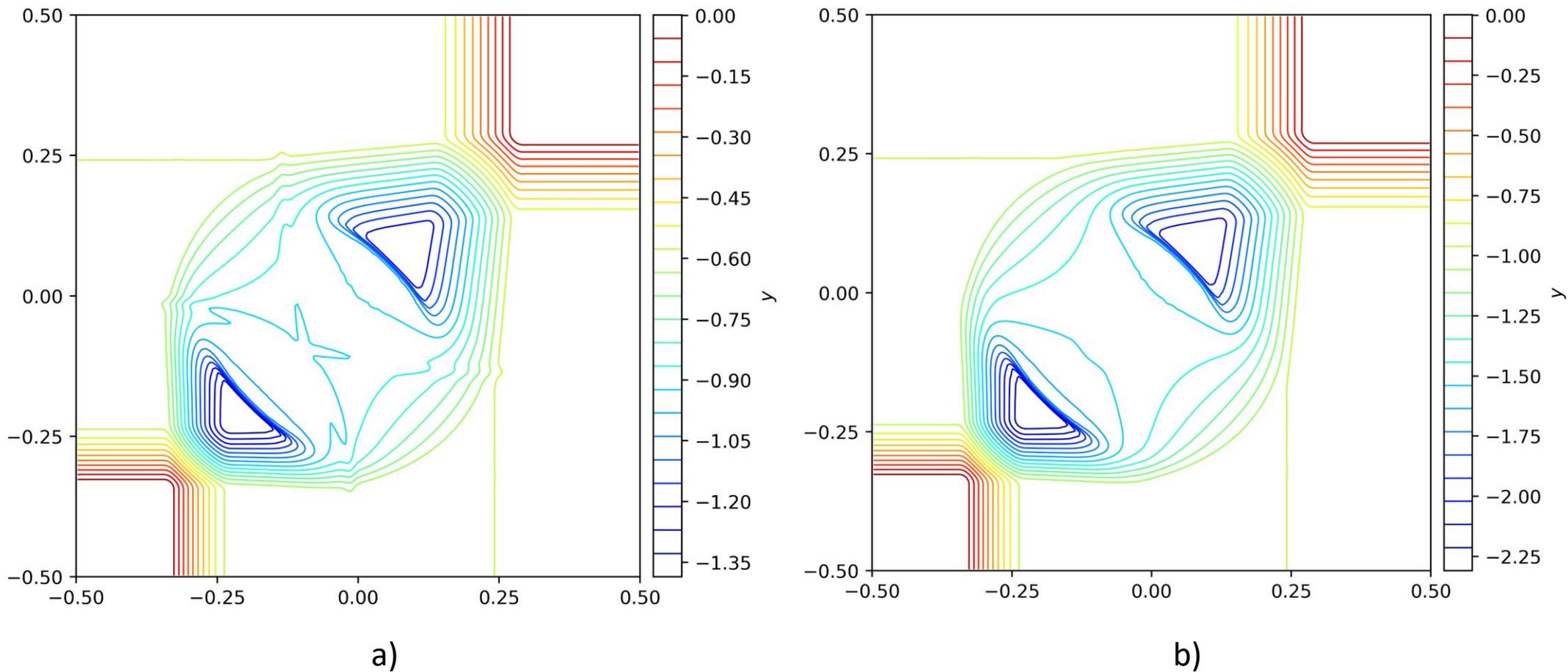

*Fig. 20) Relativistic Hydrodynamic Flow: 2DRP-2 (2D Riemann Problem-2). Fig. 20a shows the density logarithm and Fig. 20b shows the pressure logarithm at time t=0.4 obtained using the 7th order LLF-based AFD-WENO scheme with 200 ×200 zones. 25 contours were fit in between the range of minimum and maximum value. The 5th and 9th order schemes also perform well on this problem and are not shown here.*

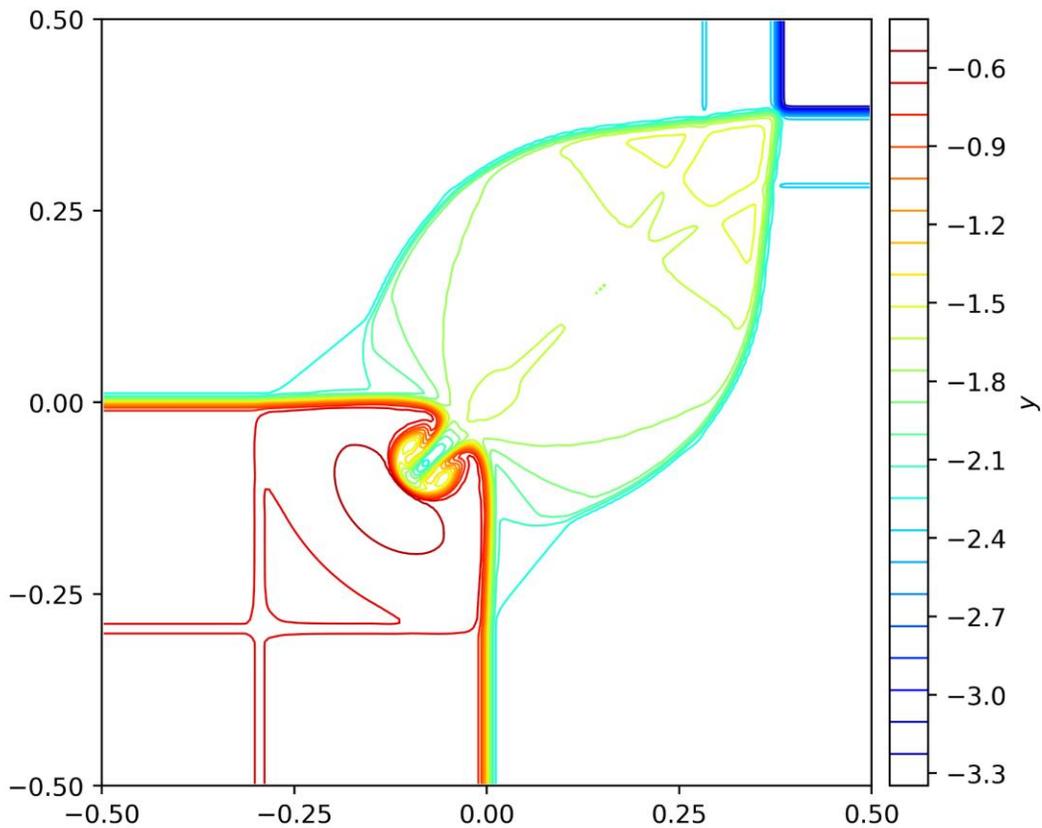 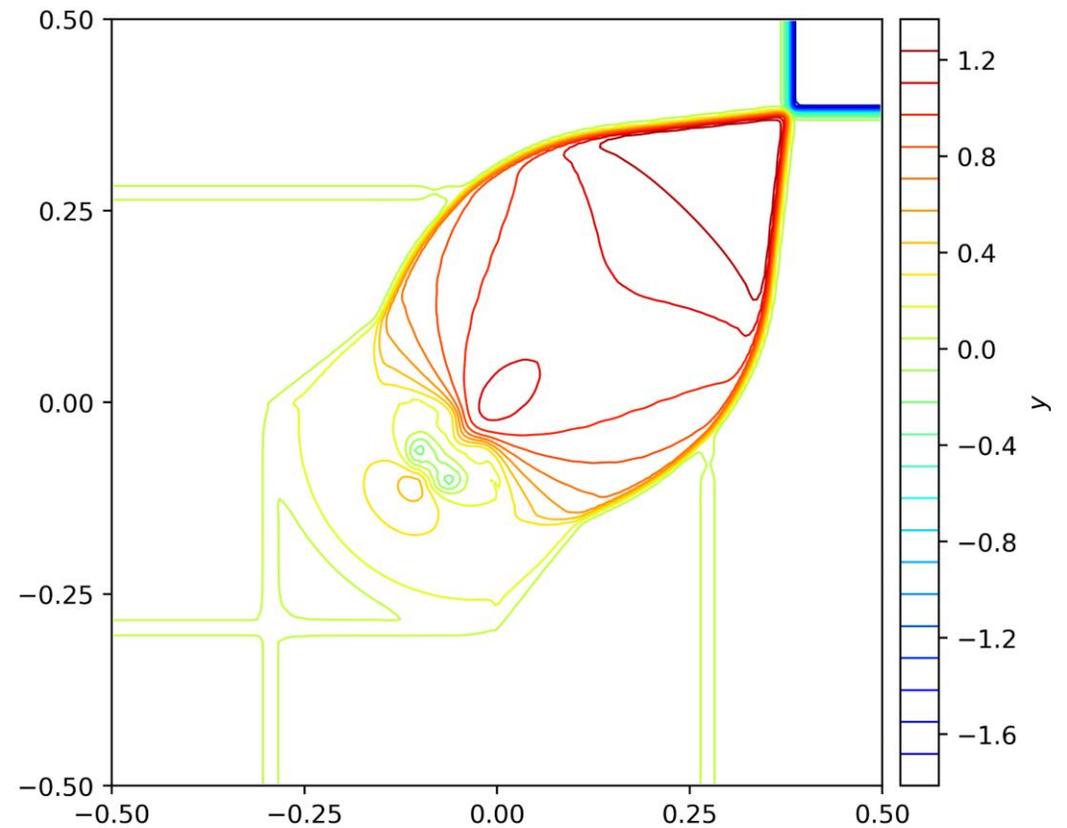

*Fig. 21) Relativistic Hydrodynamic Flow: 2DRP-3 (2D Riemann Problem-3). Fig. 21a shows the density logarithm and Fig. 21b shows the pressure logarithm at time t=0.4 obtained using the 9th order LLF-based AFD-WENO scheme with 200 ×200 zones. 25 contours were fit in between the range of minimum and maximum value. The 5th and 7th order schemes also perform well on this problem and are not shown here.*

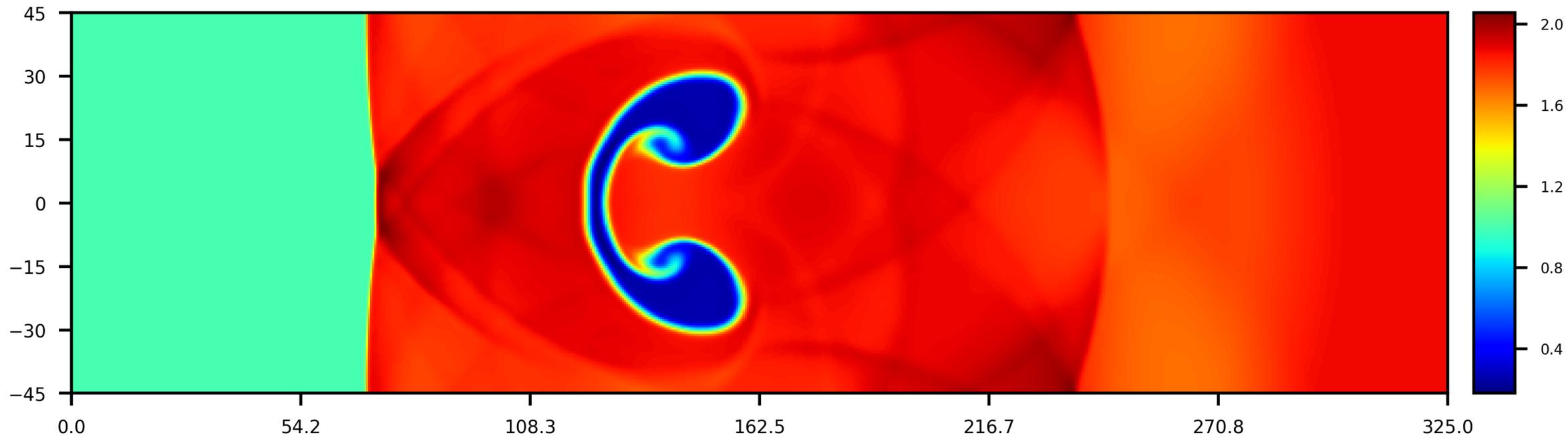

*Fig. 22a) Relativistic Hydrodynamic Flow: SB-1 (Shock-Bubble interaction-1). Panel shows the density profile at time t=450 obtained using the 3rd order LLF-based AFD-WENO scheme with 650×180 zones.*

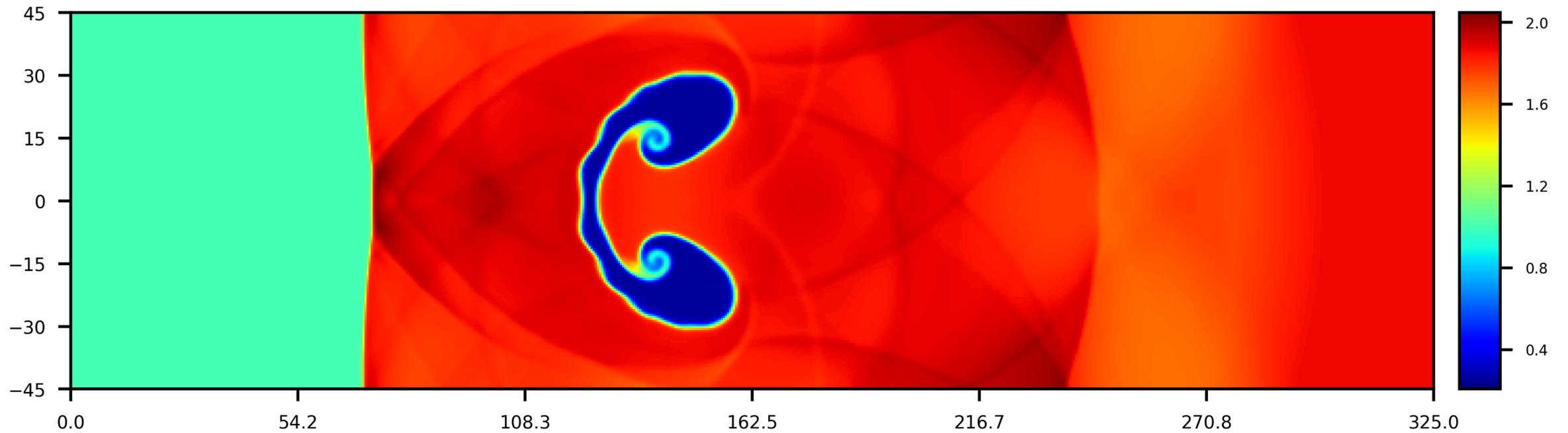

*Fig. 22b) Relativistic Hydrodynamic Flow: SB-1 (Shock-Bubble interaction-1). Panel shows the density profile at time t=450 obtained using the 5th order LLF-based AFD-WENO scheme with 650×180 zones.*

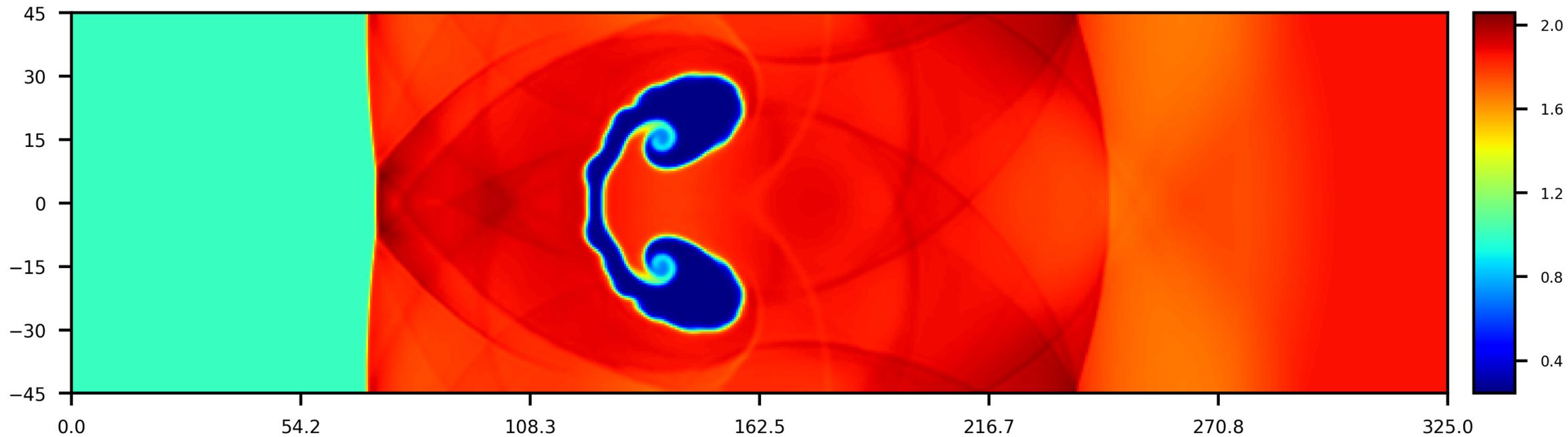

*Fig. 22c) Relativistic Hydrodynamic Flow: SB-1 (Shock-Bubble interaction-1). Panel shows the density profile at time t=450 obtained using the 7th order LLF-based AFD-WENO scheme with 650×180 zones.*

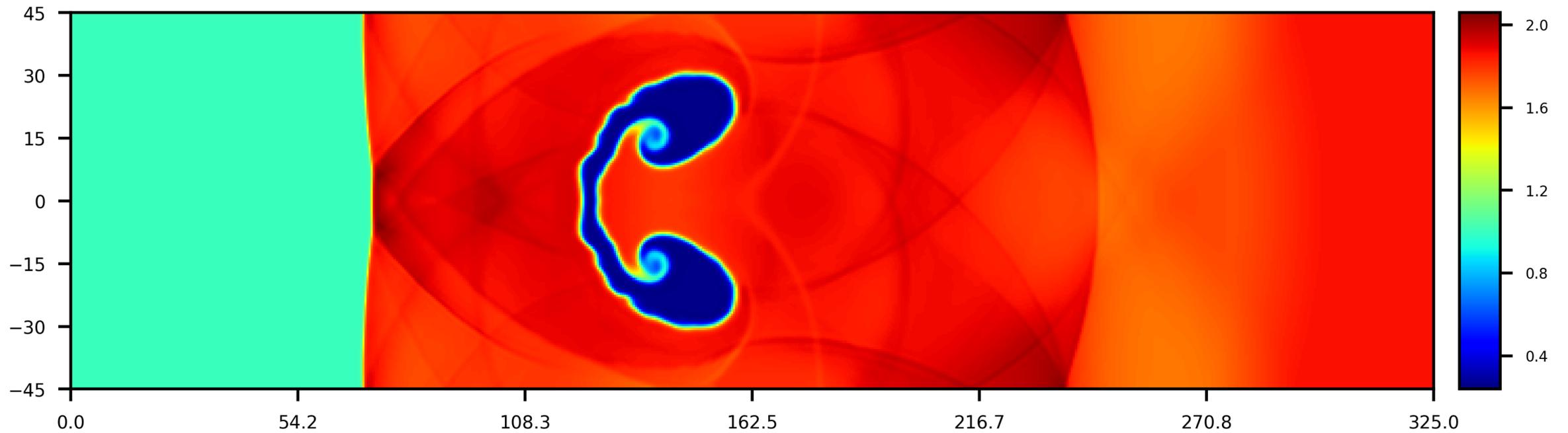

*Fig. 22d) Relativistic Hydrodynamic Flow: SB-1 (Shock-Bubble interaction-1). Panel shows the density profile at time t=450 obtained using the 9$^{th}$ order LLF-based AFD-WENO scheme with 650×180 zones.*

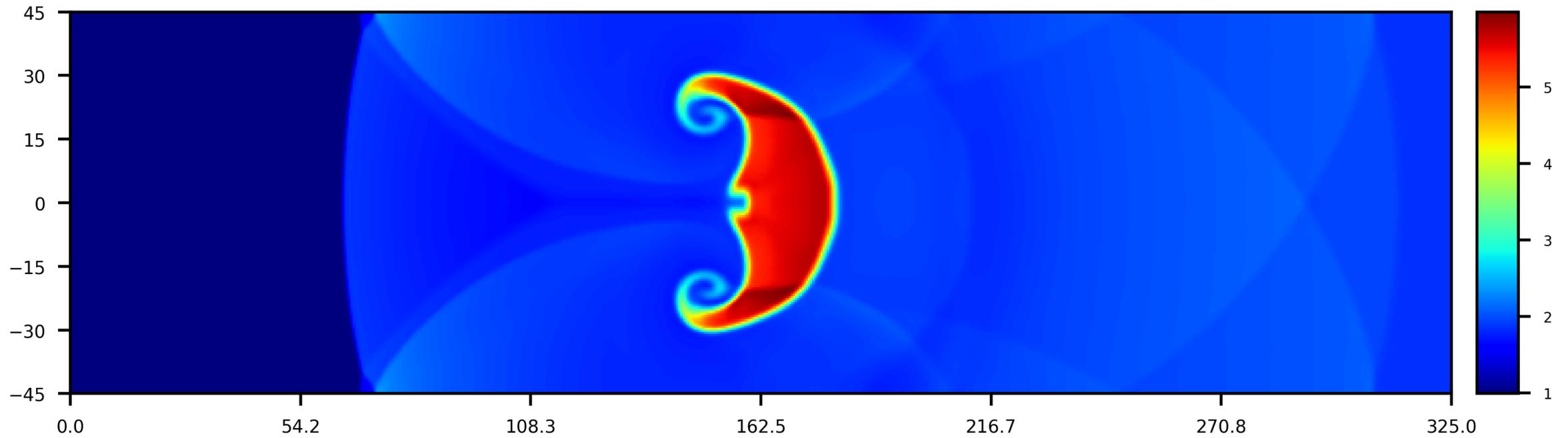

*Fig. 23a) Relativistic Hydrodynamic Flow: SB-2 (Shock-Bubble interaction-2). Panel shows the density profile at time t=500 obtained using the 3rd order LLF-based AFD-WENO scheme with 650×180 zones.*

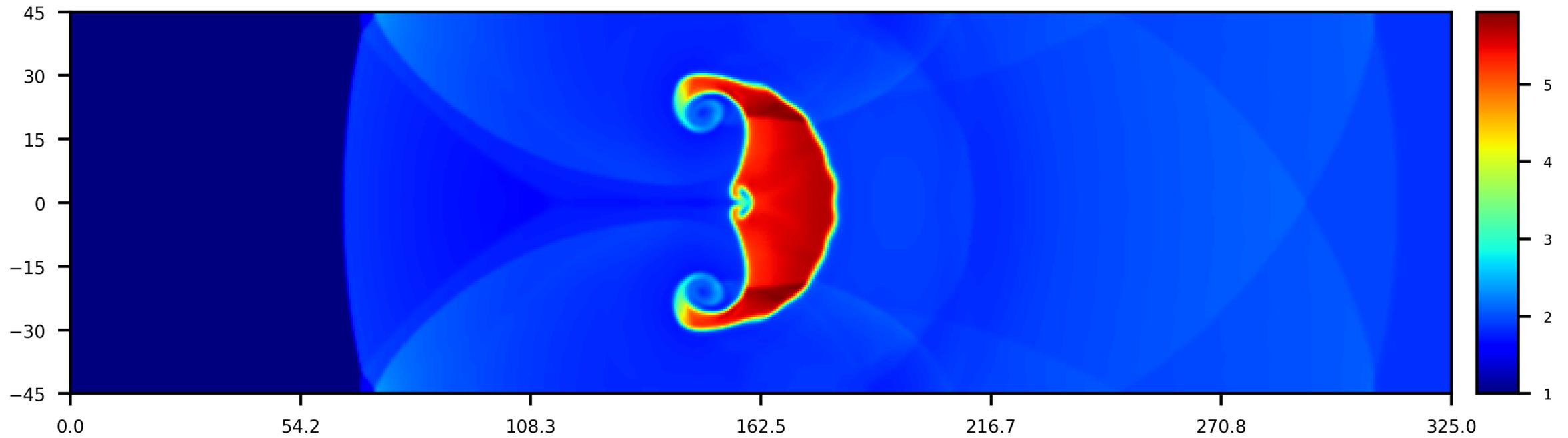

Fig. 23b) Relativistic Hydrodynamic Flow: SB-2 (Shock-Bubble interaction-2). Panel shows the density profile at time t=500 obtained using the 5$^{th}$ order LLF-based AFD-WENO scheme with 650×180 zones.

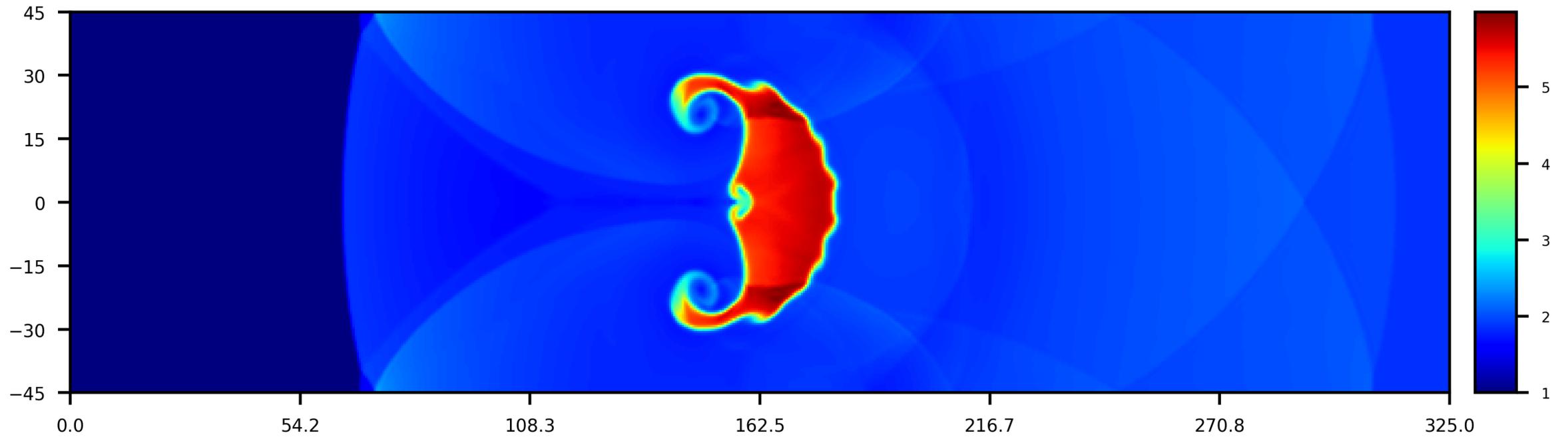

*Fig. 23c) Relativistic Hydrodynamic Flow: SB-2 (Shock-Bubble interaction-2). Panel shows the density profile at time t=500 obtained using the 7th order LLF-based AFD-WENO scheme with 650×180 zones.*

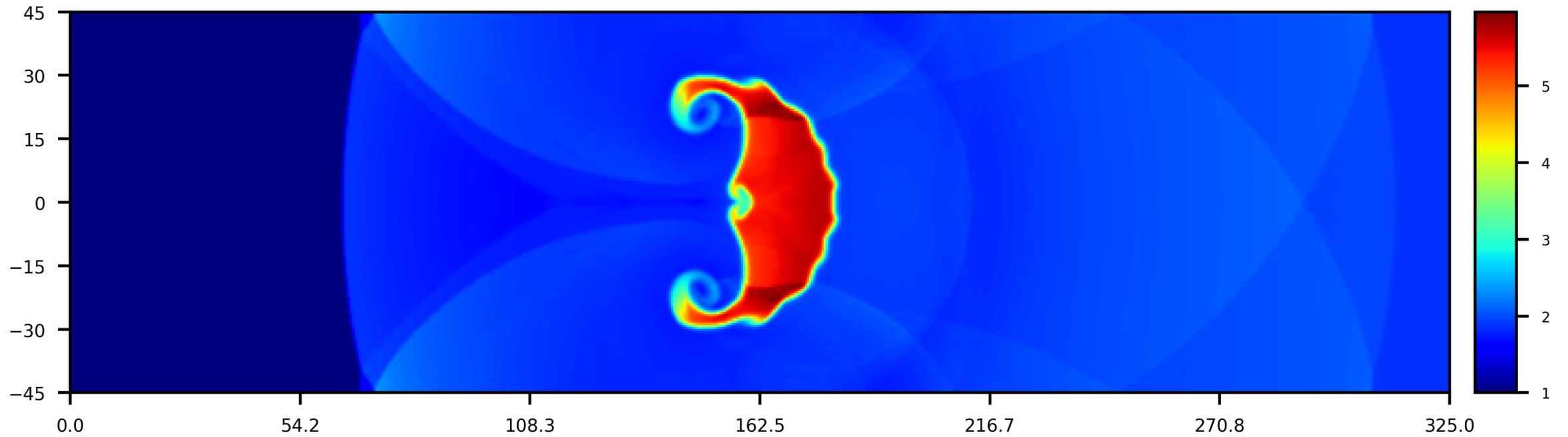

Fig. 23d) Relativistic Hydrodynamic Flow: SB-2 (Shock-Bubble interaction-2). Panel shows the density profile at time t=500 obtained using the 9$^{th}$ order LLF-based AFD-WENO scheme with 650×180 zones.

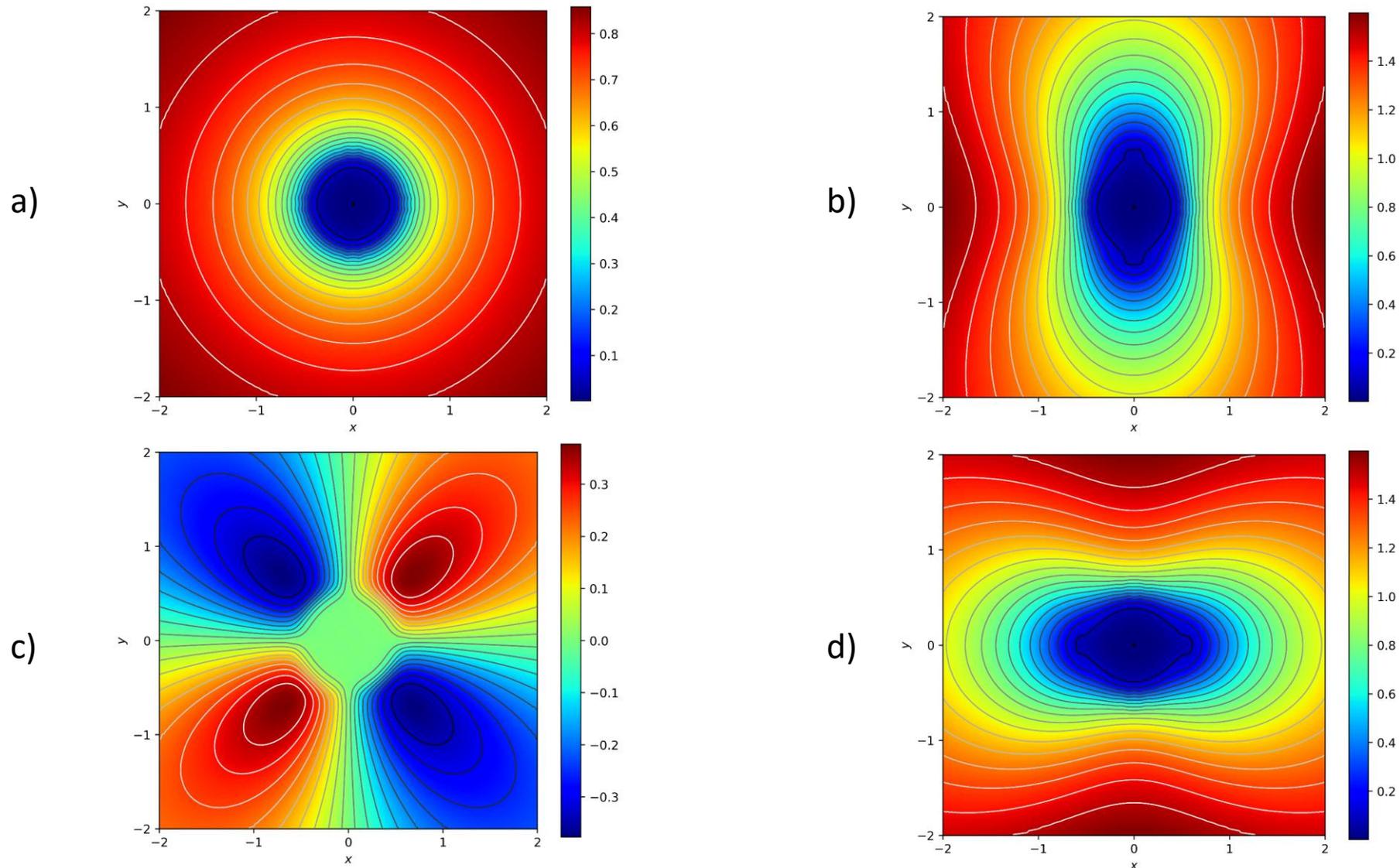

*Fig. 24) Ten-Moment Rarefied Gas flow: 2D near vacuum test problem. Fig. 24a shows the density, Fig. 24b shows the $p_{xx}$ component, Fig. 24c Shows the $p_{xy}$ component and Fig. 24d shows the $p_{yy}$ component at time t=0.05 obtained using the 7th order LLF-based AFD-WENO scheme with 200 zones. The 5th and 9th order schemes also perform well on this problem and are not shown here.*